%% file: Jiang.tex
\documentstyle{amltd}
\begin{document}
\annalsline{157}{2003}
\received{February 5, 2001}
\revised{December 9, 2001}
\startingpage{743}
\def\bye{\input{Jiang.ref} \end{document}}
 \font\tenrm=cmr10
\def\ritem#1{\item[{\rm #1}]}
\input amssym.def
\input amssym.tex
\catcode`\@=11
\font\twelvemsb=msbm10 scaled 1100
\font\tenmsb=msbm10
\font\ninemsb=msbm10 scaled 800
\newfam\msbfam
\textfont\msbfam=\twelvemsb  \scriptfont\msbfam=\ninemsb
  \scriptscriptfont\msbfam=\ninemsb
\def\msb@{\hexnumber@\msbfam}
\def\Bbb{\relax\ifmmode\let\next\Bbb@\else
 \def\next{\errmessage{Use \string\Bbb\space only in math
mode}}\fi\next}
\def\Bbb@#1{{\Bbb@@{#1}}}
\def\Bbb@@#1{\fam\msbfam#1}
\catcode`\@=12

 \catcode`\@=11
\font\twelveeuf=eufm10 scaled 1100
\font\teneuf=eufm10
\font\nineeuf=eufm7 scaled 1100
\newfam\euffam
\textfont\euffam=\twelveeuf  \scriptfont\euffam=\teneuf
  \scriptscriptfont\euffam=\nineeuf
\def\euf@{\hexnumber@\euffam}
\def\frak{\relax\ifmmode\let\next\frak@\else
 \def\next{\errmessage{Use \string\frak\space only in math
mode}}\fi\next}
\def\frak@#1{{\frak@@{#1}}}
\def\frak@@#1{\fam\euffam#1}
\catcode`\@=12

\def\alp{{\alpha}}	
\def\bet{{\beta}}	
\def\gam{{\gamma}}	
\def\Gam{{\Gamma}}	
\def\lam{{\lambda}}	
\def\Lam{{\Lambda}}	
\def\del{{\delta}}		
\def\eps{{\varepsilon}}		
\def\sig{{\sigma}}
\def\Sig{{\Sigma}}
\def\bks{{\backslash}}
\def\tet{{\theta}}
\def\zet{{\zeta}}
\def\Tet{{\Theta}}
\def\ome{{\omega}}
\def\Ome{{\Omega}}
\def\Sp{\mathop{\rm Sp}\nolimits}
\def\diag{{\rm diag}}
\def\Bil{{\rm Bil}}
\def\Gr{{\rm Gr}}

\def\tilsp{\mathop{\widetilde{\rm Sp}}\nolimits}
\def\tilspn{{\widetilde{\rm Sp}_{2n}}}
\def\spn{{{\rm Sp}_{2n}}}
\def\SO{{\rm SO}}
\def\GL{{\rm GL}}
\def\GSp{{\rm GSp}}
\def\U{{\rm U}}
\def\O{{\rm O}}
\def\Hom{{\rm Hom}}
\def\c-Ind{{\rm c-Ind}}
\def\Ind{{\rm Ind}}
\def\Im{{\rm Im}}
\def\SL{{\rm SL}} 
\def\diag{{\rm diag}}
\def\max{{\rm max}}

\def\span{\mathop{\rm span}\nolimits}
\def\Span{\mathop{\rm Span}\nolimits}
\def\tilg{{\widetilde g}}
\def\tilh{{\widetilde h}}
\def\tilP{{\widetilde P}}
\def\tilm{{\widetilde m}} 
\def\tilM{{\widetilde M}}
\def\tilN{{\widetilde N}}
\def\tilU{{\widetilde U}}
\def\oh{{\overline h}}
\def\hatrho{{\widehat\rho}}
\def\hattau{{\widehat\tau}}
\def\hata{{\widehat a}}
\def\calB{{\cal B}}
\def\tilcalb{{\widetilde{\cal B}}}
\def\calH{{\cal H}}
\def\calL{{\cal L}}
\def\nek{,\ldots,}
\def\AA{{\Bbb A}}
\def\CC{{\Bbb C}}
\font\titr=cmr10 scaled 1700
\font\emi= cmmi10 scaled 1700 
\font\titr=cmr10 scaled 1700
 \title{The local converse theorem\\ for {\titr SO(2}{\emi n}{\titr+1)} and
applications} 
\shorttitle{The local converse theorem} 

 \acknowledgements{During the work of this paper, the first named
author was partly supported by NSF Grants DMS-9896257 and DMS-0098003, by
the Sloan Research Fellowship, McKnight Professorship at University of Minnesota, 
and by NSF Grant DMS-9729992 through the Institute for
Advanced Study, Princeton, in the fall, 2000.
The second named author was supported by a grant from
the Israel-USA Binational Science Foundation.\hfill\break
\hglue19pt 1991 {\it Mathematics Subject Classification}: 11F, 22E.} 
 \twoauthors{Dihua Jiang}{David Soudry}
\institutions{University of
Minnesota, Minneapolis, MN\\
{\eightpoint {\it E-mail address\/}: dhjiang@@math.umn.edu}\\
\vglue6pt
Tel Aviv University, Tel Aviv, Israel\\
{\eightpoint {\it E-mail address\/}: soudry@@math.tau.ac.il}}

\centerline{\bf Abstract}
\vglue12pt 
In this paper we characterize irreducible generic representations of\break $\SO_{2n+1}(k)$ 
(where $k$ is a $p$\/{\rm -}\/adic field) by means of twisted local gamma factors 
(the {\it Local Converse Theorem}). As applications, we prove that two irreducible 
generic cuspidal automorphic representations of $\SO_{2n+1}({\Bbb A})$ (where 
${\Bbb A}$ is the ring of adeles of a number field) are equivalent if their local 
components are equivalent at almost all local places (the {\it Rigidity Theorem}); 
and prove the\break {\it Local Langlands Reciprocity Conjecture} for generic supercuspidal 
representations of $\SO_{2n+1}(k)$. 
 \vglue8pt 
\section{Introduction}

In the theory of admissible representations of $p$-adic reductive groups, 
one of the basic problems is to characterize an irreducible admissible 
representation up to isomorphism. 
Keeping in mind the link of the theory of admissible representations 
of $p$-adic reductive groups to the modern theory of automorphic forms, 
we consider in this paper the characterization of irreducible admissible 
representations by the local gamma factors and their twisted versions. 
Such a characterization is traditionally called the 
{\it Local Converse Theorem}, and is the local analogue of the (global) 
{\it Converse Theorem} for $\GL(n)$. We refer to \cite{CP-S1} and \cite{CP-S2} for 
detailed explanation of converse theorems.  

The local converse theorem for the general linear group, $\GL(n)$, was first 
formulated by I. Piatetski-Shapiro in his unpublished Maryland notes (1976) 
with his idea of deducing the local converse theorem from his (global) 
converse theorem. 
It was first proved by G. Henniart in \cite{Hn2} using a local approach. 
The local converse theorem is a basic ingredient in the recent proof 
of the local Langlands conjecture for $\GL(n)$ by M. Harris and R. Taylor 
\cite{HT} and by\break G. Henniart \cite{Hn3}. 

The formulation of the local converse theorem in this case is as follows. 
Let $\tau$ and $\tau'$ be irreducible admissible generic representations 
of $\GL_n(k)$, where $k$ is a $p$-adic field 
(non-archimedean local field of characteristics zero). 
Following \cite{JP-SS}, one defines the twisted local gamma factors 
$\gamma(\tau\times\varrho,s,\psi)$ and\break $\gamma(\tau'\times\varrho,s,\psi)$, 
where $\varrho$ is an irreducible admissible generic representation of 
$\GL_l(k)$ and $\psi$ is a given nontrivial additive
character of $k$. 

\proclaimtitle{Henniart, \cite{Hn2}}
\proclaim{Theorem}  \hskip-8pt
Let $\tau$ and $\tau'$ be irreducible admissible generic representations of 
$\GL_n(k)$ with the same central character. If the twisted local gamma factors 
are the same{\rm ,} i.e.\ $$
\gamma(\tau\times\varrho,s,\psi)=
\gamma(\tau'\times\varrho,s,\psi)
$$
for all irreducible supercuspidal representations $\varrho$ of $\GL_l(k)$ with 
$l=1,2,\cdots,\break n-1${\rm ,} then the representation $\tau$ is isomorphic to the 
representation $\tau'$. 
\endproclaim 

This theorem has been refined by J. Chen in \cite{Ch} (unpublished) 
so that the twisting condition on $l$ reduces from $n-1$ to $n-2$ 
(using a  local approach) and 
by J. Cogdell and I. Piatetski-Shapiro in \cite{CP-S1} (using  a global approach 
and assuming both $\tau$ and $\tau'$ are supercuspidal). 
It is expected (as a conjecture of H. Jacquet, \S 8 in \cite{CP-S1}) that the 
twisting condition on $l$ should be reduced from $n-1$ to $[{n\over2}]$. 
We note also that the local converse theorem for generic representations of 
$\U(2,1)$ and for $\GSp(4)$ was established by E. M. Baruch in \cite{B1} and 
\cite{B2}. 

The objective of this paper is to prove the local converse theorem for 
irreducible admissible generic representations of $\SO_{2n+1}(k)$. 

\proclaimtitle{The Local Converse Theorem}
\proclaim{Theorem} 
Let $\sigma$ and $\sigma'$ be irreducible admissible generic representations 
of $\SO_{2n+1}(k)$. If the twisted local gamma factors 
$\gamma(\sigma\times\varrho,s,\psi)$ and 
$\gamma(\sigma'\times\varrho,s,\psi)$ are the same{\rm ,} i.e.\ $$
\gamma(\sigma\times\varrho,s,\psi)=\gamma(\sigma'\times\varrho,s,\psi)
$$
for all irreducible supercuspidal representations $\varrho$ of $\GL_l(k)$ with 
$l=1,2,\cdots,\break 2n-1${\rm ,} then the representations $\sigma$ and $\sigma'$ are 
isomorphic. 
\endproclaim

Note that the twisted local gamma factors used here are the ones studied 
either by F. Shahidi in \cite{Sh1} and \cite{Sh2} or by D. Soudry in 
\cite{S1} and \cite{S2}. It was proved by Soudry that the twisted local 
gamma factors defined by these two different methods are in fact the same. 
It is expected that the local converse theorem (Theorem 1.2) should be 
refined so that it is enough to twist the local gamma factors in 
Theorem 1.2 by the irreducible supercuspidal representations $\varrho$ of 
$\GL_l(k)$ for $l=1,2,\cdots,n$. 
This is compatible with the conjecture of Jacquet 
as mentioned above. In a forthcoming paper of the authors, 
we shall prove the finite field analogue of Jacquet's conjecture 
and provide strong evidence for the refined local converse theorem. 

The local converse theorem for $\SO(2n+1)$ has many significant applications 
to both the local and global theory of representations of $\SO(2n+1)$. 
For the global theory, we can prove that the weak Langlands functorial lift 
from irreducible generic cuspidal automorphic representations of $\SO(2n+1)$ 
to irreducible automorphic representations of $\GL(2n)$ is injective up to 
isomorphism (Theorem 5.2) 
(The weak Langlands functorial lift in this case was recently 
established in \cite{CKP-SS}.); 
that the image of the backward lift from irreducible generic self-dual 
automorphic representations of $\GL(2n)$ to $\SO(2n+1)$ is irreducible, 
which was conjectured in \cite{GRS1} 
(The details of this application will be given in \cite{GRS5}.); and that
the {\it Rigidity Theorem} holds for irreducible generic cuspidal automorphic 
representations of $\SO(2n+1)$ (Theorem 5.3). 

Two important applications of the local converse theorem to the theory 
of admissible representations of $\SO_{2n+1}(k)$ are included in this paper. 
The first one is the explicit local Langlands functorial lifting taking irreducible 
generic supercuspidal representations of $\SO_{2n+1}(k)$ to $\GL_{2n}(k)$ 
(Theorem 6.1). Since the Langlands dual group of $\SO_{2n+1}(k)$ is $\Sp_{2n}(\CC)$, 
the Langlands functorial lift conjecture asserts that the natural embedding of 
$\Sp_{2n}(\CC)$ into $\GL_{2n}(\CC)$ yields a lift of irreducible admissible 
representations of $\SO_{2n+1}(k)$ to $\GL_{2n}(k)$. 
Let ${\cal GL}_{2n}^{\rm ifl}(k)$ (`if\hskip.5pt l' denotes the image of the functorial lifting) 
be the set of all equivalence classes of irreducible admissible generic 
representations of $\GL_{2n}(k)$ of the form  
$$
\tau=\eta_1\times\eta_2\times\cdots\times\eta_t, 
$$
where $\eta_i$ are irreducible unitary supercuspidal self-dual 
representations of\break $\GL_{2n_j}(k)$ with $j=1,2,\cdots,t$ and $\sum_{j=1}^tn_i=n$, 
such that 
\begin{itemize}
\item[(1)] $\eta_i\not\cong\eta_j$ if $i\neq j$, and 
\item[(2)] the local $L$-function $L(\eta_j,\Lambda^2,s)$ has a pole at $s=0$ 
for $j=1,2,\cdots,t$. 
\end{itemize} 
We denote by ${\cal SO}_{2n+1}^{\rm igsc}(k)$ the set of all equivalence 
classes of irreducible generic supercuspidal representations of $\SO_{2n+1}(k)$. 
We prove the local Langlands functorial conjecture for ${\cal SO}_{2n+1}^{\rm igsc}(k)$ in 
this paper. 

\proclaim{Theorem} \proclaimtitle{Local Langlands Functoriality (Theorem 6.1)}
There exists a unique bijective map 
$$
\ell\ :\ \sigma\mapsto \tau=\ell(\sigma)
$$
from ${\cal SO}_{2n+1}^{\rm igsc}(k)$ to ${\cal GL}_{2n}^{\rm ifl}(k)${\rm ,} which 
preserves the twisted local $L$\/{\rm -}\/factors{\rm ,}\break $\epsilon$\/{\rm -}\/factors and gamma \pagebreak
factors{\rm ,}  i.e.\ \begin{eqnarray*}
L(\sigma\times\varrho,s)
&=&
L(\tau\times\varrho,s),
\\
\epsilon(\sigma\times\varrho,s,\psi)
&=&
\epsilon(\tau\times\varrho,s,\psi)
\\
\noalign{\noindent 
and}
\gamma(\sigma\times\varrho,s,\psi)
&=&
\gamma(\tau\times\varrho,s,\psi)
\end{eqnarray*}
for all irreducible supercuspidal representations $\varrho$ of $\GL_l(k)$ with 
$l$ being any positive integer.
\endproclaim

The second application is the local Langlands reciprocity conjecture for irreducible 
generic supercuspidal representations of $\SO_{2n+1}(k)$ (Theorem 6.4). 
Let $W_k$ be the Weil group associated to the   local field $k$. We take  
$$
W_k\times \SL_2({\Bbb C})
$$
as the Weil-Deligne group (\cite{M} and \cite{Kn}). 
Let ${\cal G}_{2n}^{\rm ah}(k)$ be the set of conjugacy classes of admissible 
homomorphisms $\rho$ from $W_k\times{\rm SL}_2({\Bbb C})$ to $\Sp_{2n}({\Bbb C})$. 
If we write 
$$
\rho=\oplus_i\rho^0_i\otimes\lam^0_i,
$$
then the admissibility of $\rho$ means that $\rho^0_i$'s are continuous complex 
representations of $W_k$ with $\rho^0_i(W_k)$ semi-simple and $\lam^0_i$'s 
are algebraic complex representations of $\SL_2(\CC)$. The local Langlands 
reciprocity conjecture for $\SO_{2n+1}(k)$ asserts that for each local Langlands 
parameter $\rho$ in ${\cal G}_{2n}^{\rm ah}(k)$, there is a finite set $\Pi(\rho)$ 
(called the local $L$-packet associated to $\rho$) of  
equivalence classes of irreducible admissible representations of $\SO_{2n+1}(k)$, 
such that the union $\cup_{\rho}\Pi(\rho)$ gives a partition of the set of 
equivalence classes of irreducible admissible representations of $\SO_{2n+1}(k)$ 
and the reciprocity map taking $\rho$ to $\Pi(\rho)$ is compatible with various 
local factors attached to $\rho$ and $\Pi(\rho)$, respectively. 

Let ${\cal G}^0_{2n}(k)$ be the set of conjugacy classes of all $2n$-dimensional, 
admissible, completely reducible, multiplicity-free, symplectic 
complex representations $\rho^0$ of the Weil group $W_k$. Then we prove the following 
theorem. 

\proclaimtitle{Local Langlands Reciprocity Law (Theorem 6.4)}
\proclaim{Theorem} 
There exists a unique bijection 
$$
{\frak R}_{2n}\ :\ \rho^0_{2n}\mapsto {\frak R}_{2n}(\rho_{2n}^0)
$$
from the set ${\cal G}_{2n}^0(k)$ onto the set 
${\cal SO}_{2n+1}^{\rm igsc}(k)$ such that 
\begin{itemize}
\item[{\rm (L)}] $L(\rho^0_{2n}\otimes\rho^0_l,s,)
=
L({\frak R}_{2n}(\rho_{2n})\times {\frak r}_l(\rho^0_l),s)${\rm ,} 
\item[($\epsilon$)] $\epsilon(\rho^0_{2n}\otimes\rho^0_l,s,\psi)
=
\epsilon({\frak R}_{2n}(\rho_{2n})\times {\frak r}_l(\rho^0_l),s,\psi)${\rm ,} and 
\item[($\gamma$)] $\gamma(\rho^0_{2n}\otimes\rho^0_l,s,\psi)
=
\gamma({\frak R}_{2n}(\rho_{2n})\times {\frak r}_l(\rho^0_l),s,\psi)$ 
\end{itemize}
for all irreducible continuous representations $\rho^0_l$ of $W_k$ of   dimension $l$.  Here $\tau$ is the reciprocity
map to ${\rm GL}_l(k)$, obtained by {\rm [HT], [Hn3] (}\/see Theorem {\rm 6.2).} 
\endproclaim

Note that by Theorem 1.2, each local $L$-packet $\Pi(\rho)$ has at most one generic 
member. Theorem 1.4 establishes the  Langlands conjecture in this case up to 
the explicit construction of the relevant $L$-packets, which is a very interesting 
and difficult problem. We shall consider the local Langlands conjectures for 
general generic representations of $\SO_{2n+1}(k)$ and other related problems in 
a forthcoming work (\cite{JngS}).

Our proof of the local converse theorem goes as follows. 
Based on the basic properties of twisted local gamma factors established 
by D. Soudry in \cite{S1} and \cite{S2} and by F. Shahidi \cite{Sh1} 
and \cite{Sh2}, we study the existence of poles of twisted local gamma 
factors and related properties. This leads us to reduce
the proof of Theorem 1.2 to the case where both $\sigma$ and $\sigma'$ are 
supercuspidal (Theorem 5.1). 
To prove the local converse theorem for the case of supercuspidal 
representations (Theorem 4.1), we must combine the local method with the 
global method. 
More precisely, we first develop the explicit local Howe duality for 
irreducible generic supercuspidal representations of $\SO_{2n+1}(k)$ and 
$\tilspn(k)$, the metaplectic (double) cover of $\Sp_{2n}(k)$ 
(Theorem 2.2), which is more or less the local version 
of the global results of M. Furusawa \cite{F}.  Then, using the global weak 
Langlands functorial lifting from $\SO(2n+1)$ to $\GL(2n)$ \cite{CKP-SS} and 
the local backward lifting from $\GL_{2n}(k)$ to $\tilspn(k)$ \cite{GRS2} and 
\cite{GRS6}, we can basically relate our local converse theorem for $\SO(2n+1)$ 
to that for $\GL(2n)$. 
See the proof of Theorem 4.1 for details. The point here is to use 
preservation properties of twisted local gamma factors under 
various liftings (Propositions 3.3 and 3.4). 
It is worthwhile to mention here that the ideas and the methods used in this 
paper are applicable to other classical groups. 

This paper is organized as follows. 
In Section~2, we work out some explicit properties of local Howe duality for 
irreducible generic supercuspidal representations of $\SO_{2n+1}(k)$ 
and $\tilspn(k)$. 
The preservation property of (the pole at $s=1$ of) twisted local gamma 
factors under various liftings will be discussed in Section~3. 
In Section~4, we prove the local converse theorem for supercuspidal 
representations and in Section~5, we prove 
the theorem in the general case. The global applications mentioned above 
will be discussed at the end of Section~5. 
We determine in Section~6 the explicit structure of the 
image of the local Langlands functorial lifting from irreducible generic 
supercuspidal representations of $\SO_{2n+1}(k)$ to $\GL_{2n}(k)$ and 
prove the local Langlands reciprocity law for irreducible generic 
supercuspidal representations of $\SO_{2n+1}(k)$. 

Since $\SO(3)\cong{\rm PGL}(2)$, the main theorems in this paper are known in the 
case of $n=1$. Note that the theories of twisted local gamma factors for 
$\SO(3)\times \GL(r)$ via [S1,2], or via Shahidi's method, or via [JP-SS], 
for ${\rm PGL}(2)\times \GL(r)$ are all the same.
The reason for this is the multiplicativity property of gamma factors (which
is known in all cases above). This reduces comparison of gamma factors to
supercuspidal representations. Such representations can be embedded as 
components at one place of (irreducible)\break automorphic cuspidal representations,
unramified at all remaining finite places. Since gamma factors are "globally 1" 
(this is a restatement of the functional equation for the global L function),
we get the identity of the gamma factors for supercuspidal representations. From 
now on we assume that $n\geq 2$ (this will be helpful for one technical 
reason concerning the theta lifting). 
\vglue4pt

Our project on this topic was started when we attended the conference on 
Automorphic Forms and Representations at Oberwolfach (March 2000) 
organized by Professors S. Kudla and J. Schwermer. 
The main results of this paper were obtained when we participated at 
the Automorphic Forms Semester at Institut Henri Poincar{\' e} (Paris, Spring 2000)
organized by Professors H. Carayol, M. Harris,
J. Tilouine, and M.-F.  Vign{\' e}ras. This paper was finished
when the first named author was a member of the Institute
for Advanced Study (Princeton, Fall 2000). 
We would like to thank all the organizers of the above two 
research activities and the Institutes for providing a stimulating 
research environment. 
We would like to thank D. Ginzburg and S. Rallis for their 
encouragement during our work on this project. 
Our discussion with G. Henniart was very important for the proof of 
Theorem 6.4. We are grateful to him for providing us the proof of Theorem 6.3 
\cite{Hn1}. We thank the referee for his careful reading, and for his valuable
comments, questions and suggestions.

\section{Howe duality for $\SO(2n+1)$ and
$\widetilde{\Sp}(2n)$}

In this section, we prove certain properties of the local
Howe duality between $\SO_{2n+1}(k)$  and
$\widetilde{\Sp}_{2n}(k)$, applied to irreducible,
generic, supercuspidal representations, and then we
discuss relevant aspects of the global theta
correspondence for irreducible, automorphic, cuspidal
representations of $\SO_{2n+1}(\AA)$  and
$\widetilde{\Sp}_{2n}(\AA)$. Here $\widetilde{\Sp}_{2n}$ denotes 
the metaplectic (double) cover of $\Sp_{2n}$ over both the local field $k$ and the 
ring of adeles $\AA$ (\cite{Mt}).

\demo{{\rm 2.1.} Local Howe duality}
Let $k$  be a non-archimedean local field of
characteristic zero.  Let $V$ be a $(2n+1)$-dimensional vector 
space over $k$, equipped with a nondegenerate symmetric
form $(\cdot,\cdot)_V$  of Witt index $n$.  Let $W$  be a\break
$2m$-dimensional vector space over $k$,  equipped with a
nondegenerate symplectic form $(\cdot,\cdot)_W$.  We fix
a basis 
$$\{ e_1\nek e_n,e,e_{-n}\nek e_{-1}\}$$
of $V$  over $k$, such that $(e_i,e_j)_V=(e_{-i},
e_{-j})_V=0$, $(e_i,e_{-j})_V=\del_{ij}$, for
$i,j=1\nek n$,  and  we may assume that $(e,e)_V=1$.  Thus
$$V^+=\Span_k\{e_1\nek e_n\},\quad
V^-=\Span_k\{e_{-1}\nek e_{-n}\}$$ 
are dual maximal totally isotropic subspaces of $V$, and we get
a polarization of $V$,
$$V=V^++ke+V^-\ .$$
Similarly, we fix a basis 
$$\{f_1\nek f_m,f_{-m}\nek f_{-1}\}$$
of $W$ over $k$, such that
$(f_i,f_j)_W=(f_{-i},f_{-j})_W=0$  and
$(f_i,f_{-j})_W=\del_{ij}$, for $i,j=1\nek m$. Thus,
$$W^+=\Span_k\{ f_1\nek f_m\}\quad ,\qquad W^-=\Span_k
\{ f_{-1}\nek f_{-m}\}$$
are dual maximal isotropic subspaces of $W$, and we get
the polarization 
$$W=W^++W^-\ .$$

Consider the tensor product $V\otimes W$ of $V$ and $W$. It is a symplectic 
space of dimension $2m(2n+1)$, equipped with the symplectic form
$(,)_V\otimes (,)_W$. With the chosen bases, we may identify $V$ with $k^{2n+1}$ 
(column vectors) and $W$ with $k_{2m}$ (row vectors). Then we have   
$\O_{2n+1}(k)\cong \O(V)$, acting from the left on $V$, and $\Sp_{2m}(k)\cong\Sp(W)$, 
acting from the right on $W$. We let $\Sp(V\otimes W)\cong\Sp_{2m(2n+1)}(k)$ act 
from the right on $V\otimes W$. Then $\O(V)\times\Sp(W)$ is naturally embedded in
$\Sp(V\otimes W)$ by means of the following action 
$$
(v\otimes w)(g,h)=g^{-1}\cdot v\otimes w\cdot h.
$$

Let $\psi$  be a nontrivial character of $k$. The 
Weil representation $\ome_\psi$  of the metaplectic
group $\widetilde{\Sp}_{2m(2n+1)}(k)$  can be realized
in the space of Bruhat-Schwartz functions ${\cal S}(V^m)$, 
where $V^m=V\times\cdots\times V$  ($m$ copies).  We
restrict $\ome_\psi$  to the image of the natural
embedding of $\O_{2n+1}(k)\times\widetilde{\Sp}_{2m}(k)$
inside $\widetilde{\Sp}_{2m(2n+1)}(k)$, in order to study the local Howe duality 
between representations of $\O_{2n+1}(k)$ and $\widetilde{\Sp}_{2m}(k)$. 

In the following we identify
\begin{equation}
V^m\cong V\otimes W^+=V\otimes f_1\oplus\cdots\oplus V\otimes f_m. 
\end{equation}
We restrict $\ome_\psi$  to the image of the embedding
of $\O_{2n+1}(k)\times\widetilde{\Sp}_{2m}(k)$  inside
$\widetilde{\Sp}_{2m(2n+1)}(k)$.  Here are some
formulae.  Let $\varphi\in {\cal S}(V^m)$.  Then
$$\ome_\psi (g,1)\varphi (v_1\nek
v_m)=\varphi(g^{-1}\cdot v_1\nek g^{-1}\cdot v_m)$$
for $g\in \O_{2n+1}(k)$  and $(v_1\nek v_m)\in V^m$.
Next, let $\tilP_m=\tilM_m\tilN_m$  be the inverse image
in $\widetilde{\Sp}_{2m}(k)$  of the Siegel parabolic
subgroup $P_m$  of $\Sp_{2m}(k)$.  Thus,
$$\tilM_m=\left\{ (\tilm(a),\eps):\tilm(a)=
\left( \begin{array}{cc}
a&0\\
0&a^*
\end{array}\right)
\in\Sp_{2m}(k), a\in \GL_m(k),\eps =\pm 1\right\}$$
which is a semi-direct product of $\GL_m(k)$  and $\{\pm
1\}$. Note that   $\tilN_m$  is the direct product of $N_m$  and
$\{\pm 1\}$,  since the double cover splits over unipotent subgroups.
(See \cite{Mt}.)
Here
$$N_m=\left\{n(X)=\left( \begin{array}{cc}
I_m&X\\
0&I_m
\end{array}\right)
\in\Sp_{2m}(k)\right\}\ .$$
We will identify $N_m$ with $N_m\times \{1\}$. 

From the definition of the Weil (or Oscillator) representation $\ome_\psi$, 
we have that for $(m(a),\eps)\in\tilM_m$,
\begin{equation}
\ome_\psi(1,(\tilm(a),\eps))\varphi(v_1\nek v_m)
=
\chi_\psi((\det a)^m)|\det a|^{m\over 2}\varphi((v_1\nek v_m)a) \quad
\end{equation}
where $\chi_\psi$ is the character of the two-fold cover of $k$ associated 
to $\psi$ (through the Weil factor); and for $n(X)\in N_m$, 
\begin{equation}
\ome_\psi(1,n(X))\varphi(v_1\nek v_m)
=
\psi\left({1\over2}{\rm tr}[\Gr(v_1\nek v_m)Xw_m]\right)\varphi(v_1\nek v_m)
\end{equation}
where ${\rm tr}(\cdot)$ is the usual trace of a matrix, $w_m$
is the $m\times m$  matrix, whose entries are all zero except these 
along the second diagonal, which are all one, and finally 
\begin{equation}
\Gr(v_1\nek v_m)=\Big((v_i,v_j)_V\Big)_{m\times m}, 
\end{equation}
is the Gram matrix. (See (2.9) in \cite{GRS4} for more formulas.)

Let $\sig$  be an irreducible admissible representation
of $\O_{2n+1}(k)$, acting on a space $V_\sig$.  Consider,
as in p. 47 of \cite{MVW},
$$S(\sig):=\bigcap_\alp\ker(\alp)\ ,$$
where $\alp$  runs over all elements of
$\Hom_{\O_{2n+1(k)}}(S,V_\sig)$, $S=S(V^m)$.
Define
\begin{equation}
S[\sig]:=S/S(\sig)
\end{equation}
It is clear that $S[\sig]$  affords a representation of
$\O_{2n+1}(k)\times\widetilde{\Sp}_{2m}(k)$.  According to
page~47 of \cite{MVW}, $S[\sig]$  has the form
$$\sig\otimes\Tet_\psi^{n,m}(\sig)$$ 
where $\Tet^{n,m}_\psi (\sig)$ is a smooth representation
of $\widetilde{\Sp}_{2m}(k)$.  Assume that
$$
\Hom_{\O_{2n+1(k)}}(S,V_\sig )\not= 0.
$$ 
Then the Howe duality conjecture states that $\Tet_\psi^{n,m}(\sig)$
has a unique sub-representation $\Tet^{n,m}_\psi
(\sig)^0$, such that the quotient representation
\begin{equation}
\tet^{n,m}_\psi(\sig):=\Tet^{n,m}_\psi(\sig)\big/
\Tet^{n,m}_\psi (\sig)^0
\end{equation}
is irreducible.  The map taking $\sig$ to
$\tet^{n,m}_\psi(\sig)$  is called the local $\psi$-Howe
lift from $\O_{2n+1}(k)$ to $\widetilde{\Sp}_{2m}(k)$.
Similarly, in the reverse direction, given an
irreducible, admissible representation $\pi$ of
$\widetilde{\Sp}_{2m}(k)$, such that
$\Hom_{\widetilde{\Sp}_{2m}(k)}$ $(S,V_\pi)\not= 0$, we
have the spaces $S(\pi),S[\pi]$, $\Tet^\psi_{m,n}(\pi)$,
such that 
$$
S[\pi]\cong\Tet^\psi_{m,n}(\pi)\otimes\pi
$$
over $\O_{2n+1}(k)\times\widetilde{\Sp}_{2m}(k)$.  The
Howe duality conjecture states that
$\Tet^\psi_{m,n}(\pi)$ has a unique sub-representation  
$\Tet^\psi_{m,n}(\pi)^0$, such that the quotient
$$
\tet^\psi_{m,n}(\pi):=\Tet^\psi_{m,n}(\pi)\big/
\Tet^\psi_{m,n}(\pi)^0
$$ 
is irreducible.  We will say in such a case
that $\tet^\psi_{m,n}(\pi)$ is the local $\psi$-Howe
lift of $\pi$ to $\O_{2n+1}(k)$.  

In general, if $\sig$ and $\pi$  are irreducible admissible representations of
$\O_{2n+1}(k)$  and $\widetilde{\Sp}_{2m}(k)$ respectively such that
$$\Hom_{\O_{2n+1}(k)\times\widetilde{\Sp}_{2m}(k)}
(\ome_{\psi},\sig\otimes\pi)\not= 0\ ,$$
then we say that $\pi$ is a local $\psi$-Howe lift of
$\sig$, and $\sig$  is a local $\psi$-Howe lift of
$\pi$ (without assuming the existence of the local Howe duality conjecture).  
The local Howe duality conjecture was proved by
Waldspurger \cite{W}, when the residual characteristic
of $k$ is odd.  In particular, in such a case, if $\pi$
is a $\psi$-local Howe lift of $\sig$ (notations as above)
then $\pi=\tet^{n,m}_\psi(\sig)$ and
$\sig=\tet^\psi_{m,n}(\pi)$.  The following theorem of
Kudla, concerning local Howe duality for supercuspidal
representations is free from the restriction on the
residual characteristic. 

\proclaimtitle{[K1, Th.~2.1] or [MVW, \S  VI.4, Chap.~3]}
\proclaim{Theorem} 
Let $\sig$  and $\pi$  be irreducible{\rm ,} supercuspidal
representations of $\O_{2n+1}(k)$  and
$\widetilde{\Sp}_{2m}(k)$  respectively.

$(1)$  There is a positive integer $m_0=m_0(\sig)${\rm ,} such
that for any integer $1\le m<m_0${\rm ,} $\Hom_{\O_{2n+1}(k)}(S,V_\sig)=0${\rm ,}
and for any integer $m\ge m_0${\rm ,} $(\Hom_{\O_{2n+1}(k)}(S,V_\sig)\not= 0${\rm ,} 
and hence\/{\rm )} $\Tet^{n,m}_\psi(\sig)\not= 0$.  Moreover{\rm ,} 
if $m=m_0$ then $\Tet^{n,m}_\psi(\sig)$ is irreducible and supercuspidal. 
In particular{\rm ,} $\phantom{\sum^1}$
$$
\Tet^{n,m_0}_\psi(\sig)=\tet^{n,m_0}_\psi(\sig)\ .$$
If $m>m_0${\rm ,} then $\Tet^{n,m}_\psi(\sig)$  is of finite
length and is not supercuspidal.

Similar results hold for $\pi$ {\rm (}\/denote $n_0=n_0(\pi)${\rm )}.

$(2)$  We have{\rm ,}
$$\tet^{n_0,m}_\psi\left(\tet^\psi_{m,n_0}(\pi)\right)
=\pi$$
and
$$\tet^\psi_{m_0,n}\left(\tet^{n,m_0}_\psi(\sig)\right)=\sig\ .$$
{\rm (}\/We use the Weil representation as in Remark {\rm 2.3} of
{\rm {\rm \cite{K1}}.)}
\endproclaim

\demo{{R}emark {\rm 2.1}}
Since $\O_{2n+1}(k)=\{\pm
I_{2n+1}\}\times \SO_{2n+1}(k)$, every irreducible  
representation of $\O_{2n+1}(k)$  remains irreducible
upon restriction to $\SO_{2n+1}(k)$. Conversely, let
$\bet =-I_{2n+1}$.  Then for every irreducible
representation $\sig$ of $\SO_{2n+1}(k)$, $\sig=\sig^\bet$, so that $\sig$  
extends to an irreducible
representation of $\O_{2n+1}(k)$. It extends in two ways, $\sig^{+}$
and $\sig^{-}$, to $\O_{2n+1}(k)$, where $\sig^{+}(\bet)=id_{V_{\sig}}$
and $\sig^{-}(\bet)=-id_{V_{\sig}}$. Clearly,
\begin{eqnarray*}
& &\Hom_{\SO_{2n+1}(k)\times\widetilde{\Sp}_{2m}(k)}(\ome_{\psi},\sig\otimes\pi)\\
&=&
\Hom_{\O_{2n+1}(k)\times\widetilde{\Sp}_{2m}(k)}(\ome_{\psi},\sig^{+}\otimes\pi)
+
\Hom_{\O_{2n+1}(k)\times\widetilde{\Sp}_{2m}(k)}(\ome_{\psi},\sig^{-}\otimes\pi).
\end{eqnarray*}

In cases where the Howe duality conjecture holds (e.g.\ when $k$ has odd residual
characteristic) if
$$
\Hom_{\SO_{2n+1}(k)\times\widetilde{\Sp}_{2m}(k)}(\ome_{\psi},\sig\otimes\pi)\not=0,
$$ 
then exactly one of the spaces 
$$
\Hom_{\O_{2n+1}(k)\times\widetilde{\Sp}_{2m}(k)}(\ome_{\psi},\sig^{\pm}\otimes\pi)
$$ 
is nonzero.

Let $\sig$ and $\pi$ be irreducible admissible representations of $\SO_{2n+1}(k)$ 
and $\widetilde{\Sp}_{2m}(k)$ respectively. Assume that
$$
\Hom_{\SO_{2n+1}(k)\times\widetilde{\Sp}_{2m}(k)}(\ome_{\psi},\sig\otimes\pi)\not=0.
$$
Then we say that $\sig$ is a local $\psi$-Howe lift of $\pi$, and that
$\pi$ is  a local $\psi$-Howe lift of $\sig$. There shouldn't be confusion
with the similar notion for $\O_{2n+1}(k)\times\widetilde{\Sp}_{2m}(k)$. (The
groups are different.) Again, if the last condition holds and the Howe duality
conjecture is valid, then the local $\psi$-Howe lift of $\pi$ to  $\O_{2n+1}(k)$
is one of the representations $\sig^{\pm}$, denote it by $\sig^{\eps}$, and
then the local $\psi$-Howe lift of $\sig^{\eps}$ to $\widetilde{\Sp}_{2m}(k)$
is $\pi$. In general, if  $\pi$ is  a local $\psi$-Howe lift of $\sig$, then
we can assert that at least one of the representations $\sig^{\pm}$ is a local 
$\psi$-Howe lift of $\pi$.  

One of our main goals in this section is to show, for
irreducible, generic, supercuspidal representations
$\sig,\pi$ of $\SO_{2n+1}(k)$ and
$\widetilde{\Sp}_{2m}(k)$  respectively, that $n_0(\pi)=m$ and for exactly one of 
the representations $\sig^{\pm}$, denote it by $\sig^{\eps}$, $m_0(\sig^{\eps})=n$. 
(In the first case $\pi$ has to have a Whittaker model compatible with
$\psi$.)
\enddemo

Let $U_n$  (resp.\ $\tilU_m$) be the standard maximal
unipotent subgroup of\break $\SO_{2n+1}(k)$ (resp.\ $\widetilde{\Sp}_{2m}(k)$); here $\tilU_m$ is the image
of the embedding of the standard maximal unipotent
subgroup of $\Sp_{2m}(k)$ inside
$\widetilde{\Sp}_{2m}(k)$.  Let $Z_l$ be the standard
maximal unipotent subgroup of $\GL_l(k)$.  Then, since the covering of 
$\widetilde{\Sp}_{2m}(k)$ splits over unipotent subgroups (\cite{Mt}),  
$$
U_n=m(Z_n)\cdot V_n,\ \ \ \ \tilU_m=\tilm(Z_m)N_m\times 1
$$
where
\begin{eqnarray*}
 m(Z_n)&=&\left\{m(z)=\left( \begin{array}{ccc}
z&&\\
&1&\\
&&z^*
\end{array}\right)
\in \SO_{2n+1}(k):z\in Z_n\right\}\\[0.5em]
 \tilm(Z_m)&=&\left\{\tilm(z)=\left( \begin{array}{cc}
z&\\
&z^*
\end{array}\right)
\in\Sp_{2m}(k):z\in Z_m\right\}\\[0.5em]
 V_n&=&\left\{v(y,z)=\left( \begin{array}{ccc}
I_n&y&z\\
&1&y'\\
&&I_n
\end{array}\right)
\in \SO_{2n+1}(k)\right\}.
\end{eqnarray*}
We will identify $\tilm(Z_m)N_m$ with $\tilU_m$. 

Let $\psi$  be a nontrivial character of $k$. Denote by
$\psi_U$  the following nondegenerate character of
$U_n$:
\begin{equation}
\psi_U(m(z)v(y,e)):=\psi(z_{12}+\cdots +
z_{n-1,n})\psi(y_n):=\psi_n(z)\psi_U(v(y,e)) \qquad 
\end{equation}
where $m(z)v(y,e)\in m(Z_n)\cdot V_n=U_n$. 
For $\lam\in k^*$, denote by $\psi_{\tilU,\lam}$ the
nondegenerate character of $\tilU_m$, which corresponds to
$\psi$ and $\lam$:
\vglue-9pt
\begin{equation}
\psi_{\tilU,\lam}(\tilm(z)n(X)):=\psi(z_{12}+\cdots
+z_{m-1,m})\psi\left({\lam\over 2}X_{m1}\right):
=\psi_m(z)\psi_{\tilU,\lam}(n(X))
\end{equation}
where $\tilm(z)n(X)\in \tilm(Z_m)N_m=\tilU_m$. 
\vglue4pt 
An irreducible admissible representation $\sig$ (resp.\ $\pi$) of $\SO_{2n+1}(k)$ (resp.\ $\tilsp_{2m}(k)$) is
called $\psi_U$-generic (resp.\ $\psi_{\tilU,\lam}$-generic) if $\sig$ (resp.\ $\pi$)
admits a nonzero $\psi_U$ (resp.\ $\psi_{\tilU,\lam}$)
Whittaker functional, i.e.\ a nonzero element of\break
$\Hom_{U_n}(\sig,\psi_U)$ (resp.\ $\Hom_{\tilU_m}(\pi,
\psi_{\tilU,\lam})$). Note that if a
representation of\break $\SO_{2n+1}(k)$ has a Whittaker model
with respect to one nondegenerate character, then it has
a Whittaker model with respect to any nondegenerate
character, since the maximal split torus of $\SO_{2n+1}(k)$ acts 
transitively on the set of all generic characters of $U_n$. 
This is not necessarily the case for
representations of $\tilsp_{2m}(k)$.

\specialnumber{2.1}\proclaim{Proposition}
Let $\sig$  be an irreducible generic representation of
$\SO_{2n+1}(k)$.  Let $1\le m<n$  be an integer.  Then
$\sig$ has no nonzero local $\psi$\/{\rm -}\/Howe lifts to
$\tilsp_{2m}^{\phantom{1}}(k)$ {\rm (}\/and thus{\rm ,} each of the representations $\sig^{\pm}$
has no nonzero local $\psi$\/{\rm -}\/Howe lifts to
$\tilsp^{\phantom{1}}_{2m}(k)${\rm .)}
\endproclaim

\demo{Proof}
This is the local version of Proposition 2 in \cite{F}.
The proof is the appropriate analog of the proof
in \cite{F}.  
Let $m<n$, and assume that there is an irreducible
admissible representation $\pi_m$ of $\tilsp_{2m}(k)$,
acting in a (nontrivial) space $V_{\pi_m}$, which is a
local $\psi$-Howe lift of $\sig$. This means that there
is a nontrivial $\SO_{2n+1}(k)$-intertwining and
$\tilsp_{2m}(k)$-equivariant map  
$$
\rho :{\cal S}(V^m)\otimes V_{\pi^\vee_m}\longrightarrow
V_\sig.
$$
Here $\pi^\vee_m$  denotes the representation
contragredient to $\pi_m$ (acting in $V_{\pi^\vee_m}$.)
Let $\eta_{\psi_U}$  be a (nontrivial) $\psi_U$-Whittaker
functional on $V_\sig$. Consider
$$
b_{\psi_U}:=\eta_{\psi_U}\circ\rho :{\cal S}(V^m)\otimes
V_{\pi^\vee_m}\longrightarrow {\CC}
$$
which is a nontrivial bilinear form satisfying 
\begin{equation}
b_{\psi_U}(\ome_\psi(u,h)\varphi,\pi^\vee_m(h)\xi)=
\psi_U(u)b_{\psi_U}(\varphi,\xi)
\end{equation}
for $u\in U_n$, $h\in\tilsp_{2m}(k)$, $\varphi\in {\cal S}(V^m)$,
$\xi\in V_{\pi^\vee_m}$.  We will show that, for $m<n$,
the space of bilinear forms, satisfying the equivariance
property (2.9), is zero, and this will be a
contradiction.  To this end, we pass to a realization of
$\ome_\psi$ in a mixed model
$$
{\cal S}(W^n\times W^+)\cong {\cal S}(V^m)
$$
where $W^n\times W^+$ is the direct product of the spaces $W^n$ 
and $W^+$ (\S II.7, Chapter 2 in \cite{MVW}).
More precisely,
$$
[V^+\otimes W+e\otimes W^-]+[V^-\otimes W+e\otimes W^+]
$$
is a polarization of $V\otimes W$  (with respect to the 
symplectic form $(,)_V\otimes (,)_W$).  We may realize
$\ome_\psi$ in $S[V^-\otimes W+e\otimes W^+]\cong
{\cal S}(W^n\times W^+)\cong {\cal S}(W^n)\otimes {\cal S}(W^+)$.  We identify 
\begin{eqnarray*}
 (y_1\nek y_n)&\longleftrightarrow& e_{-n}\otimes
y_1 +\cdots +e_{-1}\otimes y_n\quad ,
\quad y_i\in W\\[0.5em]
    y^+&\longleftrightarrow& e\otimes
y^+\quad ,\quad y^+\in W^+\ .
\end{eqnarray*}
We keep denoting the Weil representation by $\ome_\psi$
(in the mixed model as well).  Let $\varphi\in
{\cal S}(W^n\times W^+)$, and consider an element $v(0,z)$  in
the center of $V_n$.  We have, from the definition of the mixed model of the 
Weil representations (\S II.7, Chapter 2 in \cite{MVW}), 
\begin{eqnarray}
&&\hskip-24pt(\ome_\psi(v(0,z),1)\varphi)(y_1\nek y_n;y^+)\\
&&\qquad =
\psi\left({1\over 2}{\rm tr}(\Gr(y_1\nek y_n)w_nz)\right)\varphi(y_1\nek y_n;y^+)\nonumber
\end{eqnarray}
where $\Gr(y_1\nek y_n)=((y_i,y_j)_W)_{n\times n}$.
Let $V_n(0,Z)=\{v(0,z)\in V_n\subset\SO_{2n+1}(k)\}$.  
Denote by $J_{V_n(0,Z)}$ the Jacquet functor along $V_n(0,Z)$
(with respect to the trivial character.  We view
$b_{\psi_U}$ first as a bilinear form on
$J_{V_n(0,Z)}({\cal S}(W^n\times W^+))\times V_{\pi^\vee_m}$,
satisfying (2.9). Let 
\begin{eqnarray*}
C_0&=&\{ (y_1\nek y_n;y^+)\in W^n\times W^+|(y_i,y_j)_W=0,\forall\ i,j\le n\}\\
C&=&W^n\times W^+\bks C_0. 
\end{eqnarray*}
It is clear that $C_0$ is closed in $W^n\times W^+$, and as the complement of $C_0$, 
$C$ is open in $W^n\times W^+$. We claim that
\begin{equation}
J_{V_n(0,Z)}({\cal S}(W^n\times W^+))\cong {\cal S}(C_0).
\end{equation}
Indeed, by \cite{BZ}, we have an exact sequence 
$$
0\rightarrow {\cal S}(C)\stackrel{i}{\longrightarrow}{\cal S}(W^n\times W^+)\stackrel{r}{\longrightarrow}{\cal
S}(C_0)\rightarrow 0,
$$
where $r$  is the restriction to $C_0$, and $i$ is the
embedding which takes a function supported in $C$,  and
extends it by zero to the whole of $W^n\times W^+$.
Using the exactness of Jacquet functors, we get  
$$
0\rightarrow J_{V_n(0,Z)}({\cal S}(C))\stackrel{i^*}{\longrightarrow}
J_{V_n(0,Z)}({\cal S}(W^n\times W^+))\stackrel{r^*}{\longrightarrow}
J_{V_n(0,Z)}({\cal S}(C_0))\rightarrow 0. 
$$ 
From (2.10), it follows that $J_{V_n(0,Z)}({\cal S}(C))=0$.
Note that for $(y_1\nek y_n; y^+)\break\in C$,  the
character
$$
v(0,z)\mapsto \psi \left({1\over 2}{\rm tr}(\Gr(y_1\nek y_n)w_nz)\right)
$$
is nontrivial, and hence, for $\varphi \in {\cal S}(C)$, there is
a large enough compact subgroup $\Ome_\varphi$ of
$V_n(0,Z)$, such that
\begin{eqnarray*}
& &\hskip-48pt \int_{\Ome_\varphi}(\ome_\psi(1,v(0,z))\varphi)(y_1\nek y_n;y^+)dz\\
&&\qquad = \int_{\Ome_\varphi}\psi({1\over 2}{\rm tr}(\Gr(y_1\nek y_n)w_nz))
\varphi(y_1\nek y_n;y^+)dz\\
&&\qquad = 0\  
\end{eqnarray*}
(by (2.10)). We conclude that
$$
J_{V_n(0,Z)}({\cal S}(W^n\times W^+))\cong
J_{V_n(0,Z)}({\cal S}(C_0))\cong {\cal S}(C_0)\ .
$$
With ${\cal S}(C_0)$ as a $U_n\times\tilsp_{2m}(k)$-module,
we now view $b_{\psi_U}$ as a bilinear form on
$J_{m(Z_n),\psi_n}({\cal S}(C_0))\times V_{\pi^\vee_m}$,
satisfying (2.9), where $J_{m(Z_n),\psi_n}$ denotes the
Jacquet functor along $m(Z_n)$, with respect to the
nondegenerate character $\psi_n=\psi_U\Big|_{m(Z_n)}$. 
Note that from the definition of the Weil representations on the mixed model 
(\S II.7, Chapter 2 in \cite{MVW}), 
the action of 
$$
m(z)=\left( \begin{array}{ccc}
z&&\\
&1&\\
&&z^*
\end{array}\right)\in m(Z_n)\ ,
$$
induced by $\ome_\psi$, in ${\cal S}(C_0)$  is given by
\begin{equation}
\ome_\psi (m(z),1)\varphi(y_1\nek y_n;y^+)=\varphi((y_1\nek y_n)w_nzw_n;y^+)\hskip.4in
\end{equation}
where we still use $\ome_\psi$ to denote action on
${\cal S}(C_0)$.  Consider the orbits of the action of
$w_nZ_nw_n$ on $\{(y_1\nek y_n)\in W^n\mid \Gr(y_1\nek y_n)=0\}$.  
They have the form
\begin{equation}
(0\cdots 0 x_10\cdots 0x_20\cdots 0x_{j-1}0\cdots 0
x_j0\cdots 0)w_nZ_nw_n
\end{equation}
where $\{x_1\nek x_j\}$ are linearly independent, span
a totally isotropic subspace of $W$, and the spaces of zeros in
(2.13) are of given sizes.  Let us write
$C_0=\bigcup_{0\le j}C_0(j)$, where
$$
C_0(j)=\{(y_1\nek y_n;y^+)\in C_0\mid\dim \Span \{ y_1
\nek y_n\}\le j\}.
$$
Note that $C_0(j)=C_0$, if $j\geq n$. We let $C_0(-1)$ be the empty set. 

By \cite{BZ}, we have the exact sequences
\begin{eqnarray}
&&\\
0\rightarrow J_{m(Z_n),\psi_n}({\cal S}(C_0\bks C_0(j)))
&\rightarrow& J_{m(Z_n),\psi_n}({\cal S}(C_0\bks C_0(j-1)))\nonumber\\
&\rightarrow& J_{m(Z_n),\psi_n}({\cal S}(C_0(j)\bks C_0(j-1)))
\rightarrow 0\nonumber
\end{eqnarray}
for $j=0,1\nek n$. We define the following subsets 
\begin{equation}
\Ome_{j,e}=\bigcup
\left[(0\cdots 0x_10\cdots 0x_20\cdots 0x_{j-1}0\cdots
0x_j0\cdots 0)w_nZ_nw_n\right]\times W^+\ ,\enspace
\end{equation}
where the union is taken over all the representative sets (as in (2.13))\break
$\{x_1\nek x_j\}$ 
which span a $j$-dimensional, totally isotropic subspace 
of $W$;   $e$~stands for an injective map from the index set 
$\{1\nek j\}$ of $y_1\nek y_j$ into the whole index set $\{1\nek n\}$. Then we have 
\begin{equation}
C_0(j)\bks C_0(j-1)=\bigcup_e\Ome_{j,e}. 
\end{equation}
It is clear that in (2.16), there are ${n\choose j}$ different terms. We order them 
as: $e_1,e_2\nek e_{n\choose j}$. Since the orbits $\Ome_{j,e_i}$  have the same 
dimension, they are both open and closed in $\bigcup_{i=1}^{n\choose j}\Ome_{j,e_i}$, 
and hence, we have 
\begin{equation}
J_{m(Z_n),\psi_n}({\cal S}(\cup_{i\ge k}\Ome_{j,e_i}))
=
J_{m(Z_n),\psi_n}({\cal S}(\cup_{i\ge k+1}\Ome_{j,e_i}))
\oplus
J_{m(Z_n),\psi_n}({\cal S}(\Ome_{j,e_k})).
\end{equation}

Let us now show that, for $j<n$, and any $k$,
\begin{equation}
\Hom_{m(Z_n)\times\tilsp_{2m}(k)}\Big(J_{m(Z_n),\psi_n}
\Big({\cal S}(\Ome_{j, e_k})\Big)\otimes V_{\pi^\vee_m},
\CC\Big)=0\hskip.5in
\end{equation}
where $m(Z_n)$  acts on $\CC$  according to $\psi_n$.
For this, write ${\cal S}(\Ome_{j,e_k})={\cal S}(\Ome'_{j,e_k})\otimes
{\cal S}(W^+)$, where $\Ome'_{j,e_k}$ is the first factor of 
$\Ome_{j,e_k}$ in (2.15).  Denote by $\ome'_\psi$  the
standard Weil representation of $\tilsp_{2m}(k)$ in
${\cal S}(W^+)$.  Let $\phi_1\in {\cal S}(\Ome'_{j,e_k})$, and
$\phi_2\in {\cal S}(W^+)$.  Then the action of $m(Z_n)\times
\tilsp_{2m}(k)$, induced by $\ome_\psi$, on
${\cal S}(\Ome'_{j,e_k})\otimes {\cal S}(W^+)$  is given by
$$
\ome_\psi(m(z),h)(\phi_1\otimes\phi_2)(y_1\nek y_n;y^+)
=
\phi_1((y_1\oh\nek y_n\oh)\cdot w_nzw_n)(\ome'_\psi(h)\phi_2)(y^+)
$$
where $\oh$ is the projection of $h\in\tilsp_{2m}(k)$
to $\Sp_{2m}(k)$.  Let $R_{j,e_k}$  be the stabilizer,
in $[m(Z_n)\times\tilsp_{2m}](k)$, of $(0\cdots 0
f_{-j}0\cdots 0 f_{-(j-1)}0\cdots 0 f_{-2}0\cdots
0f_{-1}0\cdots 0)$, where the ordering is given by
$e_k$, and the action is given by 
$$(y_1\nek y_n)\cdot (m(z),\oh)=(y_1\oh\nek y_n\oh)\cdot
w_nzw_n\ .$$
Then by Witt's theorem, it follows that $\Ome'_{j,e_k}$ is an 
$[m(Z_n)\times\tilsp_{2m}](k)$-orbit and hence ${\cal S}(\Ome'_{j,e_k})$ 
is isomorphic to the compactly induced representation 
$$
\c-Ind^{m(Z_n)\times\tilsp_{2m}}_{R_{j,e_k}}(1)
$$ 
as a representation of $m(Z_n)\times\Sp_{2m}(k)$.  Thus, the left-hand side of (2.18)
is isomorphic to 
$$
\Hom_{m(Z_n)\times\tilsp_{2m}(k)}
\left(\c-Ind^{m(Z_n)\times\tilsp_{2m}}_{R_{j,e_k}}(1)
\otimes\ome'_\psi\otimes\pi^\vee_m,\psi_n\right)
$$
which is, by Frobenius reciprocity, isomorphic to 
\begin{equation}
\Hom_{R_{j,e_k}}\left(\psi^{-1}_n\otimes\ome'_\psi\otimes\pi^\vee_m,1\right).
\end{equation}
The space (2.19) is zero for $j<n$, since then $R_{j,e_k}$
contains a subgroup of the form $L\times 1$,  where 
$L\subset m(Z_n)$  contains a simple root subgroup, and
hence $\psi_n|_L\not= 1$.  This proves (2.18).
It follows, from (2.9), (2.11), (2.14)--(2.18), that the
space of $b_{\psi_U}$ in (2.9) is isomorphic to
\begin{equation}
\Hom_{U_n\times\tilsp_{2m}(k)}\Big(J_{U,\psi_U}
({\cal S}(\Ome_n))\otimes V_{\pi^\vee_m},\CC\Big)\ .
\end{equation}
This space is zero  if $m<n$, since then $\Ome_n$ is
clearly empty.  This completes the proof of Proposition 2.1. 
\enddemo

Let us continue the line of argument of the proof for Proposition 2.1 
in case $m=n$.  Now re-denote $\pi=\pi_n$.  Let $T$  be
an element of the space (2.20), which we view now as
$$
\Hom_{U_n\times\tilsp_{2n}(k)}({\cal S}(\Ome'_n)\otimes
{\cal S}(W^+)\otimes V_{\pi^\vee},\CC)
$$ 
where $U_n$ acts on $\CC$ according to $\psi_U$.  Thus, we may think of $T$ as a
trilinear form $T(\phi_1,\phi_2,\xi)$.  Fixing $\phi_2\in {\cal S}(W^+)$ 
and $\xi\in V_{\pi^\vee}$, we obtain a map 
$$
\phi_1\mapsto T(\phi_1,\phi_2,\xi)
$$ 
which is a smooth distribution on the $m(Z_n)\times\Sp_{2n}(k)$ orbit $\Ome'_n$, 
and hence can be written uniquely in the form
\begin{equation}
T(\phi_1,\phi_2,\xi)
=\int_{R\bks Z_n\times\Sp_{2n}(k)}
\phi_1((f_{-n}h\nek f_{-1}h)w_nzw_n)\Phi_{\phi_2,\xi}(z,h)d(z,h)
\end{equation}
where $d(z,h)$  is a right $Z_n\times\Sp_{2n}(k)$-invariant measure on 
$R\bks Z_n\times\Sp_{2n}(k)$, and $\Phi_{\phi_2,\xi}$ is a (right) smooth function on
$Z_n\times\Sp_{2n}(k)$ and is left $R$-invariant. Here $R$ is the stabilizer 
$(R_n)$ in $Z_n\times\Sp_{2n}(k)$ of $(f_{-n}\nek f_{-1})$ under
the action (on $\Ome'_n$)
$$
(x_1\nek x_n)\cdot (z,h)=(x_1h\nek x_nh)w_nzw_n.
$$
Using the equivariance properties (with respect to
$m(Z_n)\times \Sp_{2n}(k))$
we find that
\begin{equation}
\Phi_{\phi_2,\xi}(z,h)=\psi^{-1}_n(z)\Phi_{\ome'_\psi(\tilh)
\phi_2,\pi^\vee(\tilh)\xi}(1,1)
\end{equation}
where $\tilh$  is a pre-image in $\tilsp_{2n}(k)$ of $h$.  Note again that
\begin{eqnarray*}
&&\hskip-48pt
\phi_1(f_{-n}\cdot h\nek f_{-1}\cdot h)w_nzw_n)\ome'_\psi(\tilh)\phi_2(y^+)
\\&&\qquad=
\ome_\psi(m(z),\tilh)(\phi_1\otimes\phi_2)(f_{-n}\nek f_{-1}; y^+). 
\end{eqnarray*}
Thus, the integrand of (2.21) is (using (2.22))
$$
\psi^{-1}_n(z)\Phi_{[\phi_1(f_{-n}\cdot h\nek f_{-1}\cdot h)w_nzw_n)
\ome'_\psi(\tilh)\phi_2], \pi^\vee(\tilh)\xi}(1,1).
$$
The function in square brackets in the first index of $\Phi$ is (by the last
equality)
$$y^+\rightarrow\ome_\psi
(m(z),\tilh)(\phi_1\otimes\phi_2)(f_{-n}\nek f_{-1}; y^+).$$
Substituting in (2.21), we get
\begin{eqnarray}
\quad && T(\phi_1,\phi_2,\xi)\\
&&=
\int_{R\bks Z_n\times\Sp_{2n}(k)}
\psi^{-1}_n(z)\Phi_{\ome_\psi(m(z),\tilh)
(\phi_1\otimes\phi_2)(f_{-n}\nek f_{-1};*),\pi^\vee(\tilh)\xi}(1,1)d(z,h).\nonumber
\end{eqnarray}
We have not yet used the property
\begin{equation}
T(\ome_\psi(1,v_n(t,x))(\phi_1\otimes\phi_2),\xi)=
\psi(t_n)T(\phi_1,\phi_2,\xi).
\end{equation}
It follows from the definition of the mixed model of the Weil representations 
(\S II.7, Chapter 2 in \cite{MVW}) that
\begin{equation}
\ome_\psi(1,v(t,x))\varphi(f_{-n}\nek f_{-1};y^+)
=\psi\left(\sum^n_{i=1}t_i(y^+,f_{-i})_W\right)
\varphi(f_{-n}\nek f_{-1}; y^+)\ .
\end{equation}
Using (2.22)--(2.25), we conclude that
$$\psi(t_n)\Phi_{\phi_2,\xi}(z,h)
=
\Phi_{\psi(\sum^n_{i=1}(zt)_i(*,f_{-i})w)\phi_2(*),\xi}(z,h)$$
for all $z\in Z_n$, $h\in\Sp_{2m}(k)$,  and in particular, for $z=1$,
\begin{equation}
\Phi_{\psi(\sum^n_{i=1}t_i(*,f_{-i})_W)\phi_2(*),\xi}(1,1)
=\psi(t_n)\Phi_{\phi_2,\xi}(1,1)
\end{equation}
Regarding, for fixed $\xi$, $\phi_2\mapsto \Phi_{\phi_2,\xi}(1,1)$ 
as a distribution 
on $W^+$, (2.26) implies that
it is supported at $y^+=f_n$.  Thus
\begin{equation}
\Phi_{\phi,\xi}(1,1)=W(\xi)\phi_2(f_n)
\end{equation}
for some $W(\xi)\in \CC$. This implies that the trilinear functional 
$T(\phi_1,\phi_2,\xi)$ is in fact given by the following integral 
\begin{equation}
\int^{\phantom{1}}_{R\bks Z_n\times\Sp_{2n}(k)}\ome_\psi(m(z),\tilh)
(\phi_1\otimes\phi_2)(f_{-n}\nek f_{-1}; f_n)
\psi^{-1}_n(z)W(\pi^\vee(\tilh)\xi)d(z,k)
\end{equation}
where $\tilh$ is any pre-image (in $\tilsp_{2n}(k)$) of
$h$.  Finally, note that 
$$R=\left\{\left(z,\left( \begin{array}{cc}
z&x\\
0&z^*
\end{array}\right)\right)\in Z_n\times\Sp_{2n}(k)\right\}.$$
Using the left $R$-invariance of $\Phi_{\phi_2,\xi}(z,h)$,
(2.22) and (2.27) we find that
$$W\left(\pi^\vee\left( \begin{array}{cc}
z&x\\
0&z^*
\end{array}\right)
\xi\right)=\psi_n(z)\psi\left({\textstyle {1\over 2}}x_{n1}\right)W(\xi);$$
i.e., $W$ is a (necessarily nontrivial, if $T$  is
nontrivial) $\psi_{\tilU,1}$-Whittaker functional on
$V_{\pi^\vee}$. Put
$$W^{\psi_{\tilU,1}}_{\pi^\vee}(\xi)(\tilh)=
W(\pi^\vee(\tilh)\xi)\ .$$
This is the corresponding $\psi_{\tilU,1}$-Whittaker
function.  Now we can rewrite (2.28), for
$\varphi=\phi_1\otimes\phi_2$,  as
\begin{equation}
\int_{\mu_2\cdot\tilN_n\bks\tilsp_{2n}(k)}\
\ome_\psi(1,\tilh)\varphi(f_{-n}\nek
f_{-1};f_n)W^{\psi_{\tilU,1}}_{\pi^\vee}(\xi)(\tilh)d\tilh.\hskip.5in
\end{equation}
Here $\mu_2=\{\pm 1\}$  is the kernel of the projection $\tilsp_{2n}
(k)\to\Sp_{2n}(k)$.  Note that this is the local version
of \cite[formula (12)]{F}. 
Note that the integral (2.29) converges absolutely. Indeed the integrand has 
compact support. To see this, we may assume that $\varphi=\phi_1\otimes\phi_2$,
as before, and due to the Iwasawa decomposition, it is enough to note that
$\phi_1(f_n\cdot za'\nek f_1\cdot za')W^{\psi_{\tilU,1}}_{\pi^\vee}(\xi)(a',1)$,
has compact support,
where $z\in \tilde{m}(Z_n)(k)$, $a=\diag(a_1\nek a_n)$, and $a'=\diag(a,a^*)$. 
Recall that a Whittaker function, restricted to the diagonal subgroup is "vanishing
at infinity", meaning that if $\max_{1\leq i\leq n}\{|a_i|\}$ is large, then
$W^{\psi_{\tilU,1}}_{\pi^\vee}(\xi)(a',1)$ vanishes. 
Clearly $\phi_1(f_n\cdot za'\nek f_1\cdot za')$ vanishes if 
$\max_{1\leq i\leq n}\{|a_i|^{-1}\}$ is large. We conclude that $a'$ has to lie in
a compact set of the diagonal subgroup, and hence also $z$ has to lie in 
a compact set of $\tilde{m}(Z_n)(k)$.

Let us summarize what we have shown in case $m=n$. 

\specialnumber{2.1}\proclaim{{C}orollary}
$(1)$ Let $\sig$  be an irreducible generic representation
of $\SO_{2n+1}(k)$.  Assume that $\pi$  is an irreducible
representation of $\tilsp_{2n}(k)${\rm ,} which is a local
$\psi$\/{\rm -}\/Howe lift of $\sig$.  Then $\pi$ is
$\psi^{-1}_{\tilU,1}$\/{\rm -}\/generic.
Moreover{\rm ,} the functional $b_{\psi_U}${\rm ,} viewed as a
bilinear form on $\ome^{\phantom{|}}_\psi\otimes\pi^\vee${\rm ,}  equals{\rm ,} up
to scalars{\rm ,} to $(2.29)$ {\rm (}\/with a fixed  
$\psi^{\phantom{|}}_{\tilU,1}$\/{\rm -}\/Whittaker model on $\pi^\vee${\rm ).}  The
$\psi_U$\/{\rm -}\/Whittaker model of $\sig$  is spanned by the
functions   
\begin{equation}
g\mapsto\int_{\mu_2\cdot\tilN_n\bks\tilsp_{2n}(k)}\
\ome_\psi(g,\tilh)\varphi(f_{-n}\nek f_{-1};f_n)
W^{\psi_{\tilU,1}}_{\pi^\vee}(\xi)(\tilh)d\tilh.\hskip.5in
\end{equation}

$(2)$ Let $\pi$  be an irreducible{\rm ,} supercuspidal{\rm ,}
$\psi^{-1}_{\tilU,1}$\/{\rm -}\/generic representation of
$\tilsp_{2n}(k)$.  Then $\pi$ has a nontrivial 
local $\psi$\/{\rm -}\/Howe lift to $\SO_{2n+1}(k)$.  Moreover{\rm ,}
there is a nontrivial space $t_\psi(\pi)$  of
$\psi$\/{\rm -}\/Whittaker functions on $\SO_{2n+1}(k)${\rm ,} invariant
to right translations{\rm ,} such that
\begin{equation}
\Hom_{\SO_{2n+1}(k)\tilsp_{2n}(k)}(\ome_\psi\otimes
\pi^\vee, t_\psi(\pi))\not= 0
\end{equation} and 
$t_\psi(\pi)$ is spanned by the functions $(2.30)$.
\endproclaim

\demo{Proof}
We have already shown part (1).  We now prove part (2).
Since $\pi^\vee$ is $\psi_{\tilU ,1}$-generic, we may
define the integrals (2.30), which are absolutely convergent (explained
 just before the statement of Cor.~2.1).
It is easily seen that 
these integrals are not identically zero as $(\varphi,\xi)$ varies.  
Let $t_\psi(\pi)$ be
the space of functions on $\SO_{2n+1}(k)$, spanned by
the integrals (2.30).  Note that these are
$\psi_U$-Whittaker functions on $\SO_{2n+1}(k)$, and that 
$t_\psi(\pi)$ affords a smooth representation, by right
translations, of $\SO_{2n+1}(k)$.  By construction, we
clearly have (2.31). We may, of course, substitute in (2.30)
any $g$ in $\O_{2n+1}(k)$. Denote by $t'_\psi(\pi)$ the space of
functions on $\O_{2n+1}(k)$ thus obtained; it affords, as before, 
a smooth representation by right translations of $\O_{2n+1}(k)$. We have
$$\Hom_{\O_{2n+1}(k)\tilsp_{2n}(k)}(\ome_\psi\otimes
\pi^\vee, t'_\psi(\pi))\not= 0.$$
This implies that 
   
$$\Hom_{\O_{2n+1}(k)\tilsp_{2n}(k)}(\Tet^\psi_{n,n}(\pi)\otimes\pi,t'_\psi(\pi)
\otimes\pi)\not= 0\ .$$
In particular, $\Tet^\psi_{n,n}(\pi)\not= 0$.  Since
$\pi$ is supercuspidal, $\Tet^\psi_{n,n}(\pi)$  is of
finite length as a representation of $\O_{2n+1}(k)$ (Theorem 2.1) 
and hence has an irreducible
quotient; call it $\sig'$.  We have nontrivial maps 
$$\ome_\psi\longrightarrow S[\pi]=\Tet^\psi_{n,n}(\pi)
\otimes\pi\longrightarrow \sig'\otimes\pi$$
and hence $\sig'$  is a local $\psi$-Howe lift of $\pi$ to $\O_{2n+1}(k)$.
Let $\sig$ be the restriction of $\sig'$ to $\SO_{2n+1}(k)$. Then $\sig$ is
a local $\psi$-Howe lift of $\pi$ to $\SO_{2n+1}(k)$.
\enddemo

To continue, we introduce the notion of a Bessel model of
special type for representations of $\SO_{2n+1}(k)$.
Bessel models for representations of orthogonal groups
are discussed in general in \cite{GP-SR}.

Let $Q_{n-1}=M_{n-1}V_{n-1}$ be the standard
maximal parabolic subgroup of $\SO_{2n+1}(k)$,  with
Levi subgroup isomorphic to
$\GL_{n-1}(k)\times \SO_3(k)$, and unipotent radical
$$
V_{n-1}=\left\{v'(y,z)=\left( \begin{array}{ccc}
I_{n-1}&y&z\\
&I_3&y'\\
&&I_{n-1}
\end{array}\right)
\in \SO_{2n+1}(k)\right\}.
$$
Let $V_3=\Span_k\{e_n,e,e_{-n}\}$.  We choose, for
$\lam\in k^*$, a vector $e_\lam\in V_3$, such that
$(e_\lam,e_\lam)_V=\lam$.  If $\lam$ is a square
$\alp^2$, we choose $e_{\alp^2}=\alp e$.  Define a
character $\chi_\lam$  of $V_{n-1}$ by 
$$\chi_\lam(v'(y,z))=\psi\left((y\cdot e_\lam,e_{-(n-1)})_V\right)$$
where we view $y$ as the linear map which takes
$x_1e_n+x_2e+x_3e_{-n}$ in $V_3$  to
$\sum^{n-1}_{i=1}(y_{i1}x_1+y_{i2}x_2+y_{i3}x_3)e_i$.  It follows that the 
connected component of the stabilizer of $\chi_\lam$ in $M_{n-1}$ is the subgroup
\begin{eqnarray*}
S_{M_{n-1}}(\chi_\lam)&=&\left\{s(p,d)=\left( \begin{array}{ccc}
p&&\\
&d&\\
&&p^*
\end{array}\right) \right.\\
&&\hskip-.6in\left.\phantom{s(p,d)=\left( \begin{array}{c}
p \\
d\\
 p^*
\end{array}\right)}
\in \SO_{2n+1}(k):p\in P_{n-1},d\in \SO_3(k),d\cdot e_\lam
=e_\lam\right\}\ ,
\end{eqnarray*}
where $P_{n-1}=\left\{\left( \begin{array}{ccc}
a&b\\
0&1
\end{array}\right)
\in \GL_{n-1}(k)\right\}$. Let
$$D_\lam:=\{s(I_{n-1},d)\in S_{M_{n-1}}(\lam)\}\ .$$
Note that $D_\lam$ is abelian and 
$$D_1=\left\{ s\Big(I_{n-1},\left( \begin{array}{ccc}
a&&\\
&1&\\
&&a^{-1}
\end{array}\right)
\Big):a\in k^*\right\}\ .$$
For $g\in \GL_{n-1}(k)$,  denote $m'(g)=\left( \begin{array}{ccc}
g&&\\
&I_3&\\
&&g^*
\end{array}\right)
\in M_{n-1}$.  Put
$$R_\lam:=D_\lam m'(Z_{n-1})V_{n-1}\ .$$
Let $\nu$ be a character of $D_\lam$.  $D_\lambda$ is isomorphic to the special
orthogonal group of the orthocomplement of $e_\lam$  in   
$V_3$.  Define a character of $R_\lam$  by
$$b_{(\nu,\psi,\lam)}(d\cdot
m'(z)v'(y,x))=\nu(d)\psi_{n-1}(z)\chi_\lam(v'(y,x))$$
for $d\in D_\lam$, $z\in Z_{n-1}$, $v'(y,x)\in
V_{n-1}$.  We say that an irreducible, admissible
representation $\sig$ of $\SO_{2n+1}(k)$  has a
(nontrivial) 
{\it Bessel \pagebreak model of type $(R_\lam,\nu)$\/}, if
$$\Hom_{R_\lam}(\sig, b_{(\nu,\psi,\lam)})\not= 0\ .$$
If $\nu =1$,  we say that the Bessel model (of type
$(R_\lam,1)$) is {\it special\/}.

\specialnumber{2.2}\proclaim{Proposition}
If $\sig$ is an irreducible{\rm ,} supercuspidal{\rm ,} generic
representation of $\SO_{2n+1}(k)${\rm ,}  then $\sig$  has a
nontrivial Bessel model of special type $(R_1,1)$.
\endproclaim

\demo{Proof}
Let $\sig$  be an irreducible, supercuspidal, generic
representation of $\SO_{2n+1}(k)$,  acting in a space
$V_\sig$.  Let $\eta_{\psi_U}$ be a Whittaker functional
on $V_\sig$, with respect to $(U_n,\psi_U)$, i.e.\ $$
\eta_{\psi_U}(\sig
(u)\xi)=\psi_U(u)\eta_{\psi_U}(\xi)
$$
for $u\in U_n$, $\xi\in V_\sig$.  For $\xi\in V_\sig$, let
$W_\xi(g)=\eta_{\psi_U}(\sig(g)\xi)$ be the
corresponding Whittaker function. Since $\sig$  is
supercuspidal, $W_\xi$ is compactly supported modulo $U_n$
(on the left).  Now consider    
$$\bet(\xi):=\int_{k^{n-1}\times
k^*}W_\xi\left(m\left(\left( \begin{array}{cc}
&1\\
I_{n-1}&
\end{array}\right)
\left( \begin{array}{cc}
I_{n-1}&y\\
&t
\end{array}\right)
\right)\right)dy|t|^{-n+1}d^*t$$
where for $a\in \GL_n(k)$, $m(a)=\left( \begin{array}{ccc} 
a&&\\
&1&\\
&&a^*
\end{array}\right)_{\phantom{|}}
\in \SO_{2n+1}(k)$.  
It follows from the definition of $\bet(\xi)$ and the supercuspidality of $\sig$ that 
the integral (which is a Mellin transform) is absolutely convergent, and we 
can choose $\xi$ so that $\bet(\xi)\not= 0$. By direct verification, one can
check that 
$\bet$  is a Bessel functional of special type $(R_1,1)$ attached to $\sig$. 
\enddemo

\specialnumber{2.3}\proclaim{Proposition}
Let $\sig$  be an irreducible admissible representation
of $\SO_{2n+1}(k)$.  Let $\pi$  be an irreducible
admissible $\psi_{\tilU,\lam}$\/{\rm -}\/generic representation  
of $\tilsp_{2n}(k)${\rm ,} such that $\pi$  is a local
$\psi$\/{\rm -}\/Howe lift of $\sig$.  Then $\sig$  has a
nontrivial Bessel model of special type $(R_\lam,1)$.
\endproclaim

\demo{Proof}
The idea of the proof is similar to that of the
corresponding global statement (Prop. 1 in \cite{F}).
For later needs, we consider a slightly more general
situation.  Let $\sig$  be an irreducible, admissible
representation of $\SO_{2r+1}(k)$, where $r\le n$.
 Let $\pi$  be an irreducible, admissible
$\psi_{\tilU,\lam}$-generic representations of
$\tilsp_{2n}(k)$  acting in a space $V_\pi$.  Assume
that $\pi$ is a local $\psi$-Howe lift of $\sig$.  Then
there is a nontrivial $\tilsp_{2n}(k)$-intertwining and
$\SO_{2r+1}(k)$-equivariant map 
$$
\rho : {\cal S}(V^n)\otimes V_{\sig^\vee}\longrightarrow V_\pi
$$
($V$, as before, is the vector space, of dimension
$2r+1$, over $k$, on which $\SO_{2r+1}(k)$ acts from the
left, preserving $(\cdot,\cdot)_V$.  Also,  $V_{\sig^\vee}$  is a
realization  of~$\sig^\vee$.)  Let
$\eta_{\psi_{\tilU,\lam}}$  be a (nontrivial) Whittaker
functional on $V_\pi$, with respect to
$(\tilU_n,\psi_{\tilU,\lam})$.  As in the proof of Proposition
2.1, consider the composition
$b_{\psi_{\tilU,\lam}}=\eta_{\psi_{\tilU,\lam}}\circ\rho$.
We view $b_{\psi_{\tilU,\lam}}$, as a (nontrivial)
bilinear form on ${\cal S}(V^n)\times V^{\phantom{|}}_{\sig^\vee}$  satisfying
the quasi-invariance property
\begin{equation}
b_{\psi_{\tilU,\lam}}(\ome_\psi(g,u)\varphi,\sig^\vee(g)
\xi) =\psi_{\tilU,\lam}(u)b_{\psi_{U,\lam}}(\varphi,\xi)
\end{equation}
for $u\in\tilU_n$, $g\in \SO_{2r+1}(k)$, $\varphi\in
{\cal S}(V^n)$, $\xi\in V_{\sig^\vee}$. Let
$J_{N_n,\psi_{\tilU,\lam}}$ 
denote the Jacquet functor with respect to $N_n$ and
$\psi_{\tilU,\lam}|_{N_n}$.  Then we may first view
$b_{\psi_{\tilU,\lam}}$  as a bilinear form on
$J_{N_n,\psi_{\tilU,\lam}}({\cal S}(V^n))\times
V^{\phantom{|}}_{\sig^\vee}$, satisfying
\begin{equation}
b_{\psi_{\tilU,\lam}}(\ome_\psi(g,(\tilm(z),1))\varphi,
\sig^\vee(g)\xi)=\psi_n(z)b_{\psi_{\tilU,\lam}}
(\varphi,\xi)
\end{equation}
for $z\in Z_n$. We continue denoting by $\ome_\psi$  the
action of $\tilm(Z_n)\times \SO_{2r+1}(k)$ on $\varphi$
in $J_{N_n,\psi_{\tilU,\lam}}({\cal S}(V^n))$. Let   
$$
V^n_\lam 
=\left\{(v_1\nek v_n)\in V^n:\Gr(v_1\nek v_n)=\left( \begin{array}{cc}
0_{n-1}&0\\
0&\lam
\end{array}\right)\right\}
$$
where $\Gr(v_1\nek v_n)=((v_1,v_j)_V)_{n\times n}$.  Now,  
\begin{equation}
J_{N_n,\psi_{\tilU,\lam}}({\cal S}(V^n))\cong {\cal S}(V^n_\lam)\ .
\end{equation}
This follows as in (2.11). We have the exact sequence
$$
0\rightarrow J_{N_n,\psi_{\tilU,\lam}}({\cal S}(V^n\bks V^n_\lam))
\stackrel{i}{\longrightarrow}
J_{N_n,\psi_{\tilU,\lam}}({\cal S}(V^n))\stackrel{r}{\longrightarrow}
J_{N_n,\psi_{\tilU,\lam}}({\cal S}(V^n_\lam))\rightarrow 0
$$
where $r$  is induced by restriction of functions on
$V^n$  to $V^n_\lam$, and  $i$  is induced by extending
functions on $V^n\bks V^n_\lam$  by zero.  By (2.32), we
find, as in the proof of Proposition~2.1, that
$J_{N_n,\psi_{\tilU,\lam}}({\cal S}(V^n\bks V^n_\lam))=0$.  This
proves (2.34).  Note that $N_n$  acts on ${\cal S}(V^n_\lam)$
by $\psi_{\tilU,\lam}$.  Thus, the space of bilinear
forms (2.33) is isomorphic to
\begin{equation}
\Hom_{\SO_{2r+1}(k)}\left(J_{\tilm(Z_n),\psi_n}
({\cal S}(V^n_\lam))\otimes V_{\sig^\vee},\CC\right)
\end{equation}
where $\SO_{2r+1}(k)$  acts trivially on $\CC$.  By Witt's theorem, the
orbits of\break $\SO_{2r+1}(k)\times\tilm(Z_n)$  in $V^n_\lam$
have the form
\begin{equation}
\SO_{2r+1}(k)\cdot (0\cdots 0e_10\cdots 0e_20\cdots
0e_j0\cdots 0e_\lam)Z_n\ .
\end{equation}
Here $j$  and the location of $e_1\nek e_j$  among the
zeroes determine the orbit.  As in the proof of Proposition 
2.1, the orbit (2.36) contributes zero to (2.35), as long as there are zeroes in the
representative of (2.36).  In particular, the space
(2.35) is zero, if $r<n-1$,  and hence $\sig$  cannot
have a (nontrivial) $\psi$-Howe lift to $\tilsp_{2n}(k)$, 
which is $\psi_{\tilU,\lam}$-generic.

We go back to the case of the proposition, $r=n$. As we
just explained, the space (2.35) is isomorphic to 
\begin{equation}
\Hom_{\SO_{2n+1}(k)\times\tilm(Z_n)}({\cal S}(\Ome_n)\otimes V_{\sig^\vee},\CC)
\end{equation}
where $\tilm(Z_n)$  acts on $\CC$ by $\psi_n$,
and
$$\Ome_n=\SO_{2n+1}(k)(e_1,e_2\nek e_{n-1},e_\lam)Z_n\ .$$
Note that the space ${\cal S}(\Ome_n)$ is isomorphic to the compactly induced 
representation 
$$
\c-Ind^{\SO_{2n+1}(k)\times \tilm(Z_n)}_{R'_\lam}(1),
$$
where $R'_\lam$  is the stabilizer in $\SO_{2n+1}(k)\times\tilm(Z_n)$ of
$(e_1\nek e_{n-1},e_\lam)$, consisting of elements of following type: 
\begin{equation}
\left\{\!\left(\!\left( \begin{array}{ccc}
\zet&y&b\\ 
&d&y'\\
&&\zet^*
\end{array}\right)
,\tilm\left( \begin{array}{cc}
\zet&x\\
0&1
\end{array}\right)\!\right)\!\in\! \SO_{2n+1}(k)\times\tilm(Z_n)
 \mid{\zet\in Z_{n-1},\ de_\lam{=}e_\lam\atop y\cdot e_\lam
{=}\sum^{n-1}_{i=1}x_ie_i}\right\}. 
\end{equation}
Here, we view $d$  as an element of $\SO(V_3)$  and $y$
as an element of 
$$
\Hom_k(V_3,\Span_k\{e_1\nek e_{n-1}\}).
$$  
What we proved, so far, is that the space of bilinear forms (2.32) is isomorphic to 
$$
\Hom_{\SO_{2n+1}(k)\times\tilm(Z_n)}
\left(\c-Ind_{R'_\lam}^{\SO_{2n+1}(k)\times\tilm(Z_n)}(1)\otimes\sig^\vee,
\psi_n\right),
$$
which, by Frobenius reciprocity is isomorphic to
$\Hom_{R'_\lam}({\rm res}_{R'_\lam}(\psi^{-1}_n\otimes
\sig^\vee), 1)$.  When we consider (2.38), it is easy to see
that the last space is isomorphic to 
\begin{equation}
\Hom_{R_\lam}({\rm res}_{R_\lam}(\sig),b_{(1,\psi,\lam)}).
\end{equation}
We used the fact that $\sig$  is self-dual
(see \cite[p.~91]{MVW}).  This means that $\sig$ has a
(nontrivial) Bessel model of type $(R_\lam,1)$. 
\enddemo

Let us continue the line of proof of Proposition 2.3 and consider
the case $r=n-1$.  We will keep the same notation. Since in this case $\lam$ must be 
a square (so take $\lam =1$), the space (2.35) (with
$r=n-1$) is isomorphic to
$$\Hom_{\SO_{2n-1}(k)\times\tilm(Zn)}({\cal S}(\Ome'_n)\otimes
V_{\sig^\vee},\CC_{\psi_n)})$$
where
$$\Ome'_n= \SO_{2n-1}(k)(e_1,e_2\nek e_{n-1},e)Z_n\ .$$
Again the space ${\cal S}(\Ome'_n)$ can written as a compactly induced representation 
$$
\c-Ind^{\SO_{2n-1}(k)\times\tilm(Z_n)}_{S_1}(1),
$$
where
$$
S_1=\left\{(\left( \begin{array}{ccc}
\zet&y&b\\
&1&y'\\
&&\zet^*
\end{array}\right),\tilm\left( \begin{array}{cc}
\zet&y\\
0&1
\end{array}\right))
\in \SO_{2n-1}(k)\times\tilm(Z_n)\right\}\ .$$
Thus, the space of bilinear forms (2.32) is isomorphic to
$$
\Hom_{\SO_{2n-1}(k)\times\tilm(Zn)}
\Big(\c-Ind^{\SO_{2n-1}(k)\times\tilm(Z_n)}_{S_1}(1)\otimes\sig^\vee,\psi_n\Big),
$$
which is, by Frobenius reciprocity, isomorphic to 
\begin{equation}
\Hom_{S_1}({\rm res}_{S_1}(\psi^{-1}_n\otimes\sig^\vee),1)
\cong
\Hom_{U'}(\sig^\vee,\psi_{U'})
\cong\Hom_{U'}(\sig,\psi_{U'}).\hskip.2in
\end{equation}
Here $U'$  is the standard maximal unipotent subgroup of
$\SO_{2n-1}(k)$ and $\psi_{U'}$ is its standard
nondegenerate character defined by $\psi$. Since the last
space is nontrivial, we conclude that $\sig$ is generic.

As in the proof of Proposition 2.1, where we obtained (2.29),
we may view $b_{\psi_{\tilU,\lam}}(\varphi,\xi)$,
satisfying (2.32), as a distribution on $V^n$, for fixed
$\xi$; then the content of the proof of the
isomorphism of the space (2.32) with (2.39) is that 
$b_{\psi_{\tilU,\lam}}(\varphi,\xi)$ has the form, for
$r=n$,
\begin{equation}
b_{\psi_{\tilU,\lam}}(\varphi,\xi)=\int_{S_\lam\bks
\SO_{2n+1}(k)}\ome_\psi(g,1)\varphi(e_1,e_2\nek e_{n-1},
e_\lam)\bet(\xi)(g)dg\quad
\end{equation}
where $\bet$ is a nonzero Bessel functional on $V_\sig$,
of type $(R_{\lam,1})$, and $S_\lam$ is the stabilizer in $SO_{2n+1}(k)$
of $(e_1\nek e_{n-1},e_\lam)$. Similarly, in case $r=n-1,
\lam =1$, the content of the isomorphism of the space
(2.32) and the space (2.40)) is that
$b_{\psi_{\tilU,1}}(\varphi,\xi)$ has the form
\begin{equation}
b_{\psi_{\tilU,1}}(\varphi,\xi)=\int_{C'\bks \SO_{2n-1}
(k)}\ome_\psi(g,1)\varphi(e_1,e_2\nek e_{n-1},e)
W(\sig(g)\xi)dg\quad
\end{equation}
where $W$  is a $\psi_{U'}$-Whittaker functional on
$V_\sig$, and $C'$ is the stabilizer in $\SO_{2n-1}(k)$ of $(e_1\nek e_{n-1},e)$.

Note that the integrals in (2.41) and (2.42) converge absolutely. This is shown as for 
the integral in (2.29). For example, let us sketch the convergence of (2.41) in
case $\lam=1$. By the Iwasawa decomposition, it is enough to show that
the integration along $S_1\bks B$ is carried in a compact support, 
where $B$ denotes the Borel subgroup of $\SO_{2n+1}(k)$. Thus it is enough to show 
that the following function has compact support
$$
\varphi\left(z^{-1}a^{-1}e_1\nek z^{-1}a^{-1}e_{n-1},e+\sum^n_{i=1}x_ie_i\right)
\bet(\xi)(\diag(a,u(x_n),a^*)
$$
where $a=\diag(a_1\nek a_{n-1})$, $z\in Z_{n-1}$, and 
$$
u(x_n)=\left( \begin{array}{ccc}
1&-x_n&-1/2x^2_n\\
&1&x_n\\
&&1
\end{array}\right).
$$
Since $\varphi\in S(V^n)$, looking at its last coordinate, we see that the
support in $(x_1\nek x_n)\in k^n$ is compact. Next, the function 
$\bet(\xi)(\diag(a,u(x_n),a^*))$ vanishes for $a\in (k^*)^{n-1}$, if 
$\max_{1\leq i\leq n-1}\{|a_i|\}$ is large, and 
$x_n$ remains in a compact set of $k$. The proof for this is the same as
for Whittaker functions. Denote the last function by $f(a,x_n)$. We can find
a unipotent element $u$ in $\SO_{2n+1}(k)$ close to the identity, such that
it fixes $\bet(\xi)$. We then get $f(a,x_n)=b_{(1,\psi,1)}(u^a)f(a,x_n)$ (see
the paragraph before Prop.\ 2.2 for notation), where $u^a$ denotes conjugation of
$u$ by $\diag(a,I_3,a^*)$. From the last equality for any $u$ close enough
to the identity, we conclude that if $f(a,x_n)$ is nonzero, then the 
coordinates of $a$ are bounded (above). 
Finally, if 
$\varphi(z^{-1}a^{-1}e_1\nek z^{-1}a^{-1}e_{n-1},e+\sum^n_{i=1}x_ie_i)$
is nonzero, then $\max_{1\leq i\leq n-1}\{|a_i|^{-1}\}$ is bounded, and then
$z$ must lie in a compact set as well.
 
We summarize.

\specialnumber{2.2}\proclaim{{C}orollary}
Let $\pi$  be an irreducible $\psi_{\tilU,\lam}$\/{\rm -}\/generic
representation of $\tilsp_{2n}(k)$.
\vglue4pt
$(1)$ Assume that $\sig$ is an irreducible representation of
$\SO_{2n+1}(k)${\rm ,} which is a local $\psi$\/{\rm -}\/Howe lift  
of $\pi$.  Then $\sig$  has a Bessel model of special
type $(R_{\lam,1})$.  Moreover{\rm ,} the functional
$b_{\psi_{\tilU,\lam}}${\rm ,} viewed as a bilinear form on
$\ome_\psi\otimes\sig$ $(\cong\ome_\psi\otimes\sig^\vee)$
has the form $(2.41)${\rm ,} where $\bet$ is a Bessel functional
on $V_\sig$, of type $(R_\lam,1)$.  The
$\psi_{\tilU,\lam}$\/{\rm -}\/Whittaker model of $\pi$  is spanned
by the functions
\begin{equation}
h\mapsto\int_{S_\lam\bks \SO_{2n+1}(k)}\ome_\psi(g,h)
\varphi(e_1,e_2\nek e_{n-1},e_\lam)\bet(\xi)(g)dg\ .\hskip.5in
\end{equation}
\vglue4pt
$(2)$  The representation $\pi$ has no nontrivial local
$\psi$\/{\rm -}\/Howe lifts to $\SO_{2r+1}(k)${\rm ,} for $r<n-1$.
\vglue4pt
$(3)$  Assume that $\sig$  is an irreducible representation
of $\SO_{2n-1}(k)${\rm ,} which is a local $\psi$\/{\rm -}\/Howe lift of
$\pi$.  Then $\lam$  is a square {\rm (}\/take $\lam =1${\rm )}
and $\sig $  is generic.  Moreover{\rm ,} the functional
$b_{\psi_{\tilU,1}}$ has the form $(2.42)$.  The
$\psi_{\tilU,1}$\/{\rm -}\/Whittaker model of $\pi$ is spanned  
by the functions
\begin{equation}
h\mapsto\int_{C'\bks \SO_{2n-1}(k)}\ome_\psi(g,h)\varphi
(e_1,e_2\nek e_{n-1},e)W(\sig(g)\xi)dg.\hskip.5in
\end{equation}
\endproclaim

\specialnumber{2.4}\proclaim{Proposition}
Let $\sig$  be an irreducible{\rm ,} generic{\rm ,} supercuspidal
representation of $\SO_{2n+1}(k)$.  Then $\sig$  has a
nontrivial local $\psi$\/{\rm -}\/Howe lift to $\tilsp_{2n}(k)$.
Moreover{\rm ,} there is a nontrivial space $t_\psi(\sig)$ of 
$\psi_{\tilU,1}$\/{\rm -}\/Whittaker functions on
$\tilsp_{2n}(k)${\rm ,} which is invariant to right translations and 
is spanned by the functions $(2.43)$ with $\bet$ a Bessel functional 
on $V_\sig$ of special type $(R_1,1)${\rm ,} such
that
\begin{equation}
\Hom_{\SO_{2n+1}(k)\times\tilsp_{2n}(k)}(\ome_\psi\otimes\sig,
t_\psi(\sig))\not= 0.
\end{equation}
\endproclaim

\demo{Proof}
By Proposition~2.2, $\sig$  has a nontrivial Bessel functional
$\bet$ of special type $(R_1,1)$.  Consider the
integrals (2.43) (with this $\bet$ and $\lam =1$).  They
converge absolutely (as shown before Cor.~2.2). 
Since the space consisting of the functions $\varphi(g^{-1}(e_1\nek e_{n-1}, e_n))$ 
contains the space ${\cal S}(S_1\bks SO_{2n+1}(k))$, 
it follows that integrals in (2.43) as the $(\varphi,\xi)$ vary 
cannot be identically 
zero for any given $\bet(\xi)(g)$, by means of the usual density argument.

Let $t_\psi(\sig)$ be the space of functions on
$\tilsp_{2n}(k)$  spanned by the integrals (2.43).
Then $t_\psi(\sig)$  consists of $\psi_{\tilU,1}$-Whittaker
functions and affords a smooth representation, by right
translations of $\tilsp_{2n}(k)$.  By construction, we
clearly have (2.45). Write (2.43) as a sum of two terms
$W^+_{\varphi,\xi}(h) + W^-_{\varphi,\xi}(h)$, where
\begin{eqnarray*}
 &&\hskip-20pt W^{\pm}_{\varphi,\xi}(h)
\\&&  ={1\over2}\int_{S_\lam\bks \SO_{2n+1}(k)}(\ome_\psi(g,h)
\varphi\pm\ome_\psi(-g,h)\varphi)(e_1,e_2\nek e_{n-1},e_\lam)\bet(\xi)(g)dg .\end{eqnarray*}
Denote by $t^{\pm}_\psi(\sig)$ the space spanned by the functions 
$W^{\pm}_{\varphi,\xi}$ as $ \varphi$ varies in $S$ and $ \xi$ varies
in $V_\sig $.
Since $t_\psi(\sig)$ is nontrivial, one of the spaces $t^{\pm}_\psi(\sig)$
say $t^{\eps}_\psi(\sig)$, is nontrivial. Of course  $t^{\eps}_\psi(\sig)$,
is a nontrivial  $\tilsp_{2n}(k)$ - submodule of $t_\psi(\sig)$ and
$$\Hom_{\tilsp_{2n}(k)\O_{2n+1}(k)}(\ome_\psi\otimes\sig^\eps,t^{\eps}_\psi(\sig))\not=0$$
and hence    
$$\Hom_{\tilsp_{2n}(k)\O_{2n+1}(k)}\Big(
\Tet^{n,n}_\psi(\sig^\eps)\otimes\sig^\eps,
t_\psi^\eps(\sig)\otimes\sig^\eps\Big)\not= 0\ .$$
In particular, $\Tet^{n,n}_\psi(\sig^\eps)\not= 0$, and since
$\sig^\eps$  is supercuspidal, we conclude, as in the proof
of Corollary~2.1(2), that $\sig^\eps$ has a nontrivial local
$\psi$-Howe lift to $\tilsp_{2n}(k)$. In particular, $\sig$ has a nontrivial
local $\psi$-Howe lift to $\tilsp_{2n}(k)$.
\enddemo

The main theorem of this section is:

\specialnumber{2.2}\proclaim{Theorem}
Let $\sig$  and $\pi$  be irreducible, supercuspidal
representations of $\SO_{2n+1}(k)$ and $\tilsp_{2n}(k)$
respectively. 
Assume that $\sig$  is generic and that $\pi$  is
$\psi_{\tilU,1}$\/{\rm -}\/generic.  Then
\begin{itemize}
\item [{\rm (1)}] $\sig$  has a unique nontrivial local
$\psi$-Howe lift to $\tilsp_{2n}(k)$.  This lift is
supercuspidal and $\psi_{\tilU,1}$\/{\rm -}\/generic.
\item [{\rm (2)}] For $n\ge 2${\rm ,} $\pi$  has a unique nontrivial
local $\psi$\/{\rm -}\/Howe lift to $\SO_{2n+1}(k)$.
This lift is supercuspidal and generic. 
\end{itemize}

\endproclaim

\demo{Proof}
Let $\sig$ be an irreducible, supercuspidal, generic
representation of $\SO_{2n+1}(k)$.  By Proposition~2.4, $\sig$
has a nontrivial local $\psi$-Howe lift to
$\tilsp_{2n}(k)$. Let $ \sig^\eps$ be as in the proof of Prop. 2.4.
By Proposition~2.1, $n=n_0(\sig^\eps)$,  and
hence, by Theorem 2.1,
$\Tet^{n,n}_\psi(\sig^\eps)=\tet^{n,n}_\psi(\sig^\eps)$ is
irreducible and supercuspidal.  By the proof of Proposition~2.4, we have
$$
\Hom_{\tilsp_{2n}(k)\O_{2n+1}(k)}
\Big(\tet^{n,n}_\psi(\sig^\eps)\otimes\sig^\eps,t^\eps_\psi(\sig)\otimes
\sig^\eps\Big)\not= 0 .
$$
This implies that $\tet^{n,n}_\psi(\sig^\eps)$ is a
sub-representation of $t_\psi^\eps(\sig)$.  Since
$t_\psi^\eps(\sig)$  is realized in a space of
$\psi_{\tilU,1}$-Whittaker functions,
$\tet^{n,n}_\psi(\sig^\eps)$  is also
$\psi_{\tilU,1}$-generic. Put $ \pi_\eps=\tet^{n,n}_\psi(\sig^\eps)$.
Of course  $ \pi_\eps $ is a nontrivial local $\psi$-Howe lift of $\sig$ to
$\tilsp_{2n}(k)$. Note now that it is impossible to have both $\pi_+$
and $\pi_-$   nontrivial, since in such a case, we will get that 
both $\pi_\pm$ are $\psi_{\tilU,1}$-generic local $\psi$-Howe lifts of $\sig.$
By Corollary~2.2(1), both $\psi_{\tilU,1}$-Whittaker models of $\pi_\pm$ are given
by the spans of the integrals (2.43). This implies that $\pi_+ \cong \pi_-$,
and hence, by Theorem 2.1, $\sig^+ \cong \tet^\psi_{n,n}(\pi_+)
\cong \tet^\psi_{n,n}(\pi_-) \cong \sig^-.$
This is impossible since $\sig^+$ and $\sig^-$ are not isomorphic.
Thus, if $\pi$ is a  local $\psi$-Howe lift of $\sig$ to
$\tilsp_{2n}(k)$, then $ \pi$ is a  local $\psi$-Howe lift of one of the
representations $\sig^\pm$, say $\sig^{\eps'}$, and then it follows, by the 
above and Theorem 2.1 that $\pi_{\eps'} \cong \tet^{n,n}_\psi(\sig^{\eps'}) \cong \pi.$
This forces $\eps'=\eps$ and $\pi \cong \pi_\eps.$
 This proves part (1).    Note that
$\tet^{n,n}_\psi(\sig^\eps)\cong t_\psi(\sig)$. This follows from
Corollary~2.2(1).  Indeed, since $\tet^{n,n}_\psi(\sig^\eps)$  is
a local $\psi$-Howe lift of $\sig$, and it is
$\psi_{\tilU,1}$-generic, then its
$\psi_{\tilU,1}$-Whittaker model is spanned by the
functions (2.43), i.e.\ the spanning set of $t_\psi(\sig)$.

Let $\pi$  be an irreducible, supercuspidal,
$\psi_{\tilU,\lam}$-generic representation of
$\tilsp_{2n}(k)$.  Assume that $n\ge 2$.  We claim that
$\pi$  has no nontrivial local $\psi$-Howe lift to  
$\SO_{2n-1}(k)$.  Otherwise, if $\sig'$  is an
irreducible representation of $\SO_{2n-1}(k)$, which is a
local $\psi$-Howe lift of $\pi$, then by Corollary~2.2, parts (2), (3), 
$\sig'$  is supercuspidal and generic.  By part (1)
of Theorem 2.2   (just proved), $\sig'$  has a
nontrivial, supercuspidal ($\psi_\tilU$-generic)
$\psi$-Howe lift to $\tilsp_{2n-2}(k)$. This contradicts 
the tower principle of Theorem 2.1.  (The $\psi$-Howe
lifts of the supercuspidal representation $\sig'$
of $\SO_{2n-1}(k)$ to both $\tilsp_{2n-2}(k)$  and
$\tilsp_{2n}(k)$  are nontrivial and supercuspidal.)
Note that supercuspidal Weil representations of
$\widetilde{\rm SL}_2(k)$ do lift to $\SO_1(k)$. We
conclude from this, Corollary~2.2(2) and Cor.~2.1(2) that
$n_0(\pi)=n$, and hence, by Theorem 2.1,
$\Tet^\psi_{n,n}(\pi)=\tet^\psi_{n,n}(\pi)$ is
irreducible and supercuspidal.  From the proof of Corollary~2.1(2), 
using the same notation, it follows that 
$$
\Hom_{\tilsp_{2n}(k)\O_{2n+1}(k)}
\Big(\tet^\psi_{n,n}(\pi)\otimes\pi,t'_\psi(\pi)
\otimes\pi\Big)\not= 0\ .
$$
This implies that $\tet^\psi_{n,n}(\pi)$ is a
sub-representation of $t'_\psi(\pi)$.  Since $t'_\psi(\pi)$
is realized in a space of $\psi_{\tilU,1}$-Whittaker
functions, $\tet^\psi_{n,n}(\pi)$  is also
$\psi_{\tilU,1}$-generic. (Note again that
$\tet^\psi_{n,n}(\pi)\cong t'_\psi(\pi)$,  as follows
from Corollary~2.1(1).)  This completes the proof of
Theorem 2.2. 
\enddemo

\demo{{\rm 2.2.} Relation to global Howe duality}
In this subsection, we realize an irreducible, generic,
supercuspidal representation as a local component of an
irreducible, automorphic, cuspidal, generic
representation and discuss the local-global relation. 

Let $F$  be a number field and $\nu_0$ be a finite 
place of $F$, such that $F_{\nu_0}=k$.  Let $\sig$  be
an irreducible, supercuspidal, generic representation of
$\SO_{2n+1}(k)$, and let $\pi$ be the (unique) local $\psi$-Howe lift
of $\sig$ to $\tilsp_{2n}^{\phantom{l}}(k)$. By
Theorem 2.2, $\pi$  is an irreducible, supercuspidal,
$\psi_{\tilU,1}$-generic representation of
$\tilsp_{2n}(k)$. There is an element, which is a square in $k^*$, 
$\alp^2$  (and we may even take it to lie in $1+{\cal P}_{\nu_0}^N$,
where  ${\cal P}_{\nu_0}$ is the maximal ideal in the ring of integers
of $k$), such that there is a nontrivial character $\psi_0$ of
 $\AA/F$ ($\AA$  is the adele ring of $F$), such that 
$\psi_{0,\nu_0}=\psi^{\alp^2}$, i.e.\ $\psi_{0,\nu_0}(x)=\psi(\alp^2x)$, for all
$x$ in $F_{\nu_0}$. Modifying $\psi$ by $\psi^{\alp^2}$ is not harmful, since
$\ome_\psi \cong \ome_{\psi^{\alp^2}}$ and $\pi$ is $\psi_{\tilU,1}$-generic,
if and only if it is $\psi^{\alp^2}_{\tilU,1}$-generic. Note that $\pi$ is generic 
with respect to the character 
$$
u \rightarrow \psi(u_{12}+\nek+u_{n-1,n}+u_{n,n+1}), 
$$
if and only if it is generic with respect to the character
$$
u \rightarrow \psi(c_1u_{12}+\nek+\, c_{n-1}u_{n-1,n}+c_n^2u_{n,n+1}),
$$ 
for any $c_1,\nek,c_n$ in $k^*$ (this is clear from the action by conjugation of the
diagonal subgroup on $\tilU$). Thus, we may replace $\psi$ by $\psi^{\alp^2}$,
and hence we may just assume that $\alp=1$, so that there is a nontrivial 
character $\psi_0$ of $\AA/F$, satisfying  $\psi_{0,\nu_0}=\psi$. 

Let $S_0$ be the (finite) set of finite places $\nu$ of $F$, which satisfy (at least) one of 
the following conditions.
\begin{itemize}
\item [1)] The residual characteristic of $F_\nu$ is the same as that 
of  $F_{\nu_0}$.
\item [2)] The residual characteristic of $F_\nu$ is two.
\item [3)] $\psi_{0,\nu}$ is not normalized (i.e.\ its conductor is not the ring of 
integers at $\nu$).
\item [4)] The residual characteristic at $\nu$ is equal to that of a place
$\nu'$, which satisfies the previous condition.
\end{itemize}

Choose for each place $\nu\in S_0$ an irreducible, supercuspidal, 
$(\psi_{0,\nu})_{\tilU_\nu,1}$-generic representation $\pi_\nu$ of 
$\tilsp_{2n}(F_\nu)$, such that $\pi_{\nu_0}=\pi$. As in Theorem 2.2 of [V],
there is an irreducible, automorphic, cuspidal, 
$(\psi_0)_{\tilU_\AA,1}$-generic representation $\Pi$  of $\tilsp_{2n}(\AA)$,
such that $\Pi_\nu\cong\pi_\nu$, for all $\nu\in S_0$. (Here we abuse
notation and denote by $\tilU_\AA$ the adele analogue \pagebreak of $\tilU$.)  

Consider the global theta lift $\Tet(\Pi,\psi_0)$
of $\Pi$ from $\tilsp_{2n}(\AA)$ to $\SO_{2n+1}(\AA)$.
To be consistent with our local set up,
$\Tet(\Pi,\psi_0)$  is spanned by
$$ 
g\mapsto\int_{\Sp_{2n}(F)\bks\tilsp_{2n}(\AA)}
\tet^\phi_{\psi_0}(g,h)\overline{\xi(h)}dh.
$$
Here $\tet^\phi_{\psi_0}(g,h)$ is the theta series for the
dual pair $\SO_{2n+1}(\AA)\times\tilsp_{2n}(\AA)$
associated to $\psi_0$ and the Schwartz function $\phi$, and 
$\xi$ varies in the space of $\Pi$. By Proposition~3 in [F],
$\Tet(\Pi,\psi_0)$  is nontrivial and generic (in the
sense that the $\psi_0$-Whittaker coefficient (and hence any other
Whittaker coefficient) is
nontrivial on $\Tet(\Pi,\psi_0)$).  We claim that $\Tet(\Pi,\psi_0)$ is cuspidal.  
Otherwise, there is an
integer $m<n$, such that $\Tet_m(\Pi,\psi_0)$, the theta
lift of $\Pi$  to $\SO_{2m+1}(\AA)$, $\Tet_m(\Pi,\psi_0)$,  is nontrivial.
Take the first such $m$.  Then $\Tet_m(\Pi,\psi_0)$  is
cuspidal.  Let $\Sig$ be an irreducible summand of
$\Tet_m(\Pi,\psi_0)$. We clearly have, at $\nu_0$,
$$
\Hom_{\SO_{2m+1}(k)\times\tilsp_{2n}(k)}
\Big(\ome^{(m,n)}_\psi\otimes\pi^\vee,\Sig_{\nu_0}\Big)
\not= 0
$$
where $\ome^{(m,n)}_\psi$  is the Weil representation
for the dual pair $\SO_{2m+1}(k)\times\tilsp_{2n}(k)$.
Thus $\pi$  has a local $\psi$-Howe lift to
$\SO_{2m+1}(k)$.  This contradicts the tower principle
of Theorem 2.1, since $\pi$  already has a local $\psi$-Howe lift
to the supercuspidal representation $\sig$  on
$\SO_{2n+1}(k)$. Since $\Tet(\Pi,\psi_0)$  is cuspidal,
we can decompose $\Tet(\Pi,\psi_0)=\bigoplus\Sig_i\ $
into a direct sum of irreducible summands.  Note that 
each summand is (irreducible and) cuspidal. Since
$\Tet(\Pi,\psi_0)$  is generic, there is a summand, call
it $\Sig$, which is generic.  Since $\Sig_{\nu_0}$  is a
local $\psi$-Howe lift of $\pi$, we find by Theorem 2.2,
that $\Sig_{\nu_0}\cong\sig$. Similarly, for any other $\nu\in S_0$, 
$\Sig_\nu$  is a local $\psi$-Howe lift of $\pi_\nu$, and again, since
$\pi_\nu$ is a supercuspidal, $(\psi_{0,\nu})_{\tilU_\nu,1}$-generic 
representation $\pi_\nu$ of
$\tilsp_{2n}(F_\nu)$, we get (by Theorem 2.2) that $\Sig_\nu$
is a supercuspidal and generic. This proves:
\enddemo

\specialnumber{2.5}\proclaim{Proposition}
Let $\sig$ be an irreducible{\rm ,} supercuspidal{\rm ,} generic
representation of $\SO_{2n+1}(k)$.  Let
$\pi$ be the local $\psi$\/{\rm -}\/Howe
lift of $\sig$  to $\tilsp_{2n}(k)$. Let $F$ be a number field{\rm ,}
and let $\nu_0$ be a place of $F${\rm ,} such that $F_{\nu_0}\cong k$. Assume that there is $\psi_0$ as before. Then
there exists an irreducible{\rm ,} automorphic{\rm ,} cuspidal{\rm ,}
$(\psi_0)_{\tilU_\AA,1}$\/{\rm -}\/generic representation $\Pi$ of 
$\tilsp_{2n}(\AA)$  such that
\begin{itemize}
\item[{\rm 1)}] $\Pi_{\nu_0}\cong\pi$.
\item[{\rm 2)}] There exists an irreducible{\rm ,} automorphic{\rm ,}
cuspidal{\rm ,} generic representation $\Sig$ of
$\SO_{2n+1}(\AA)${\rm ,} such that
$\Sig\subset\Tet(\Pi,\psi_0)$ and
$\Sig_{\nu_0}\cong\sig$.
\end{itemize}
Moreover{\rm ,} let{\rm ,} for each $\nu\in S_0, \nu\not=\nu_0$ {\rm (}$S_0$ as before\/{\rm ),}
$\pi_\nu$ be an irreducible{\rm ,} supercuspidal{\rm ,}
$(\psi_{0,\nu})_{\tilU_\nu,1}$-generic representation of
$\tilsp_{2n}(F_\nu)$. Then we may take $\Pi$ and $\Sig$ as above{\rm ,}
such that $\Pi_\nu\cong \pi_\nu$ {\rm (}\/and hence $\Sig_\nu${\rm ,} being the 
local $\psi_{0,\nu}$ Howe lift of $\pi_\nu$, is supercuspidal and generic\/{\rm ),} 
for all $\nu\in S_0, \nu\not=\nu_0$. 
\endproclaim

\section{Local gamma factors} 

The basic theory of local gamma factors for $\SO(2n+1)$ and the twisted ones 
has been established by F. Shahidi (\cite{Sh1} and \cite{Sh2}) and D. Soudry 
(\cite{S1}, \cite{S2}, and \cite{S3}), by different methods. 
We shall discuss the basic properties of local gamma factors related to 
various lifting problems of representations of $p$-adic groups. We will also 
discuss local gamma factors for metaplectic groups. 

\demo{{\rm 3.1.} Basic facts on local gamma factors}
Let $k$ be a non-archimedean local field of characteristic zero with 
residual field consisting of $q$ elements. We shall recall mainly from 
\cite{S1} some basic facts on local gamma factors, for $\SO(2n+1)$
and later on we  recall, mainly from [GRS2,3], local gamma factors for 
metaplectic groups. 

Let $\sigma$ be an irreducible admissible generic representation of 
$\SO_{2n+1}(k)$ with ${\cal W}(\sigma,\psi)$  the associated standard 
Whittaker model with respect to the additive character $\psi$. 
Let $\varrho$ be an irreducible 
admissible generic representation of $\GL_l(k)$ with 
${\cal W}(\varrho,\psi^{-1})$ 
 the associated standard Whittaker model with respect to the additive 
character $\psi^{-1}$. 
Let $P_l=M_lN_l$ be the standard maximal parabolic
subgroup of (split) $\SO_{2l}(k)$ with the Levi subgroup
$M_l\cong \GL_l(k)$.  Let $I(\varrho,s)$ be the unitarily
induced representation of $\SO_{2l}(k)$ from $P_l$, which
is realized in the space of all smooth functions: 
$$
\Phi_{\varrho,s}\ :\ 
\SO_{2l}(k)\rightarrow {\cal W}(\varrho,\psi^{-1})
$$
satisfying the following condition:
$$
\Phi_{\varrho,s}(m(a)ny)(x)=|\det a|^{s+{l-1\over 2}}\Phi_{\varrho,s}(y)(xa)
$$
where $y\in \SO_{2l}(k)$, $m(a)\in M_l(k)$, $n\in N_l(k)$, and 
$x\in \GL_l(k)$. For the sake of convenience, 
we shall write $\Phi_{\varrho,s}$ as a $\CC$-valued function with two variables: 
$$
\Phi_{\varrho,s}(y)(x)=\xi_{\varrho,s}(y;x),
$$
where $y\in \SO_{2l}(k)$ and $x\in \GL_l(k)$. 

The local Rankin-Selberg convolution ${\cal A}(W_\sigma;\xi_{\varrho,s})$ 
for $\SO_{2n+1}(k)\times \GL_l(k)$ is defined by formulas (1.2.3) and (1.3.1) in 
\cite{S1} for $n\leq l$ and $n\geq l$, respectively. 
To illustrate the construction, 
we recall the definition of the integral ${\cal A}(W_\sigma;\xi_{\varrho,s})$ 
for the case $l\leq n$. For $d\in M_{(n-l)\times l}(k)$, we set 
$$
{\overline{x}}(d):=
\left(\begin{array}{ccccc}
I_l&0      & &       &    \\
d  &I_{n-l}& &       &    \\
   &       &1&       &    \\
   &       & &I_{n-l}&    \\
   &       & &d'     &I_l\end{array}\right)\in
\SO_{2n+1}(k). 
$$
Set ${\overline{X}}_{(l,n)}
:=\{{\overline{x}}(d)\ |\ d\in M_{(n-l)\times l}(k)\}$. For 
$W_\sigma\in {\cal W}(\sigma, \psi)$ and $\xi_{\varrho,s}\in I(\varrho,s)$, 
the local Rankin-Selberg convolution integral 
${\cal A}(W_\sigma;\xi_{\varrho,s})$ 
is defined as formula (1.3.1) in \cite{S1} by 
\begin{equation}
{\cal A}(W_\sigma;\xi_{\varrho,s})
:=
\int_{V_l\backslash \SO_{2l}(k)}
\int_{{\overline{X}}_{(l,n)}}
W_\sigma({\overline{x}}(d)j_{l,n}(h))\xi_{\varrho,s}(h,I_l)
d{\overline{x}}(d)dh,\hskip.25in
\end{equation}
where $j_{l,n}(h)$ is the embedding of $\SO_{2l}$ into $\SO_{2n+1}$ given by 
$$
h=\left(\begin{array}{ccc}
A&B\\
C&D\end{array}\right)\mapsto 
\left(\begin{array}{ccc}
A&      &B\\
 &I_{2(n-l)+1}& \\
C& &D\end{array}\right). 
$$
$V_l$ is the standard maximal unipotent subgroup of
$\SO_{2l}(k)$.
As in Chapter 9 of \cite{S1}, to obtain a functional equation for local 
Rankin-Selberg convolution integrals, one applies the intertwining operator 
$$
M(w_l,\cdot)\ :\ \ I(\varrho,s) \rightarrow I(\varrho^{w_l},s), 
$$
which is defined by 
$$
M(w_l,\xi_{\varrho,s})(y;x):=\int_{N_l}\xi_{\varrho,s}(w_ln_ly;x)dn_l
$$
where $w_l$ is the Weyl element in $\SO_{2l}$ as defined in \S 9.1.
Another local Rankin-Selberg convolution integral 
${\tilde{\cal A}}(W_\sigma;\xi_{\varrho,s})$  
for $\SO_{2n+1}(k)\times \GL_l(k)$ 
can be defined as in \S 9.2, \S 9.3, \S 9.4 and \S 9.5 of \cite{S1} 
for various different cases. Now ${\tilde{\cal A}}(W_\sigma;\xi_{\varrho,s})$ is
obtained with slight modification from ${\cal A}(W_\sigma;\xi_{\varrho,s})$
by substitution of $M(w_l,\xi_{\varrho,s})$ instead of $\xi_{\varrho,s}$.  

By Theorem 10.1 of \cite{S1}, there exists a rational function in $q^{-s}$, 
$\Gamma(\sigma\times\varrho,s,\psi)$, such that 
\begin{equation}
\Gamma(\sigma\times\varrho,s,\psi)\cdot {\cal A}(W_\sigma;\xi_{\varrho,s})
=
{\tilde{\cal A}}(W_\sigma;\xi_{\varrho,s}).
\end{equation}
Let $\gamma(\varrho,\Lambda^2,2s-1,\psi)$ be the local coefficient 
(local gamma factor) 
of Shahidi \cite{Sh2}, 
corresponding to the intertwining operator $M(w_l,\cdot)$. 
Then the local Rankin-Selberg convolution gamma factor 
$\gamma(\sigma\times\varrho,s,\psi)$ 
(or simply the local gamma factor of $\sigma$ twisted by $\varrho$) 
is defined as in \S 10.1 in \cite{S1} by the identity 
\begin{equation}
\gamma(\sigma\times\varrho,s,\psi)\cdot {\cal A}(W_\sigma;\xi_{\varrho,s})
=
{\tilde{\cal A}}^*(W_\sigma;\xi_{\varrho,s}),
\end{equation}
where 
$$
{\tilde{\cal A}}^*(W_\sigma;\xi_{\varrho,s})
=
\gamma(\varrho,\Lambda^2,2s-1,\psi)\cdot
{\tilde{\cal A}}(W_\sigma;\xi_{\varrho,s}).
$$
Hence we have the formula
\begin{equation}
\Gamma(\sigma\times\varrho,s,\psi)
=
{\gamma(\sigma\times\varrho,s,\psi)\over
\gamma(\varrho,\Lambda^2,2s-1,\psi)}, 
\end{equation}
which implies that $\gamma(\sigma\times\varrho,s,\psi)$ is a rational 
function in $q^{-s}$. 

We recall from Chapter 11 in \cite{S1} and from \cite{S3}
the theorem on the multiplicativity 
of the local twisted gamma factor $\gamma(\sigma\times\varrho,s,\psi)$. 
See also Shahidi's work in \cite{Sh1}. 
\enddemo

\proclaimtitle{Multiplicativity of gamma factors (\cite{S1} and \cite{S3})}
\proclaim{Theorem} 
$(1)$~Suppose that an irreducible admissible generic representation $\sigma$ of 
$\SO_{2n+1}(k)$ is a subquotient of 
$\Ind^{\SO_{2n+1}(k)}_{P_r}(\tau_r\otimes \sigma_{n-r})${\rm ,} 
the unitarily induced representation from a standard maximal parabolic 
subgroup $P_r$ of $\SO_{2n+1}(k)${\rm ,} 
where $\tau_r$ is an admissible generic representation of $\GL_r(k)$ and 
$\sigma_{n-r}$ is an admissible generic representation of $\SO_{2(n-r)+1}(k)$. 
Then 
$$
\gamma(\sigma\times\varrho,s,\psi)
=
\omega_\tau(-1)^n
\gamma(\tau_r\times\varrho,s,\psi)
\cdot
\gamma(\sigma_{n-r}\times\varrho,s,\psi)
\cdot
\gamma(\tau_r^\vee\times\varrho,s,\psi),
$$
for any irreducible admissible generic representations $\varrho$ of 
$\GL_l(k)$ with $l$ being any positive integer{\rm ,} 
where $\gamma(\tau_r\times\varrho,s,\psi)$ and 
$\gamma(\tau_r^\vee\times\varrho,s,\psi)$ 
are the local gamma factors defined as in {\rm \cite{JP-SS}  (}$\tau^\vee$ is 
the contragredient representation of $\tau${\rm .)}
\vglue4pt
$(2)$\ 
Suppose that an irreducible admissible generic representation $\varrho$ of 
$\GL_l(k)$ is a subquotient of 
$\Ind^{\GL_l(k)}_{P_{r,l-r}(k)}(\tau_r\otimes\tau_{l-r})${\rm ,}
the unitarily induced representation from a standard maximal parabolic 
subgroup $P_{r,l-r}$ of $\GL_l${\rm ,} 
where $\tau_r$ is an admissible generic representation of $\GL_r(k)$ and 
$\tau_{l-r}$ is an admissible generic representation of $\GL_{l-r}(k)$. 
Then 
$$
\gamma(\sigma\times\varrho,s,\psi)
=
\gamma(\sigma\times\tau_r,s,\psi)
\cdot
\gamma(\sigma\times\tau_{l-r},s,\psi),
$$
for any irreducible admissible generic representations $\sigma$ of 
$\SO_{2n+1}(k)$. 
\endproclaim

A similar theory of local gamma factors $\gamma(\pi\times\varrho,s,\psi)$
can be inferred from [GRS3], for a $\psi_{\tilU,1}$-generic representation
$\pi$ of $\tilsp_{2n}(k)$ and a generic representation $\varrho$ of $\GL_l(k)$.
In this paper we  need the case $n<l$ only  (more precisely $l=2n$), 
and this is explained in [GRS3, \S 6.2] for $k$ non-archimedean. 
(The case $n=l$ is covered in [GPS]
and the case $n>l$ can be done similarly; it follows closely the analogous 
case of $\SO(2n+1)\times \GL(l)$, shown in [S1, \S 8.1, 8.2].) The case where
$k$ is archimedean can be done as well, similarly to [S2, \S 3]. However,
in this paper we need less, namely we may assume that $\pi$ and $\varrho$ are
components at one place of globally generic automorphic forms. In this case,
the local functional equation at one place (giving rise to the corresponding
local gamma factor) follows easily from the global functional equation
satisfied by the global integrals. We review this in the appendix (\S 7.1). 
When we return to the non-archimedean local field $k$, $\pi$, $\varrho$, 
$\psi$, $n<l$, as above, the local gamma factor is the proportionality 
factor in a local functional equation of the form
$$
{\gamma(\pi\times\varrho,s,\psi)\over
\gamma(\varrho,s-1/2,\psi)\gamma(\varrho,\Lambda^2,2s-1,\psi)}
{\cal B}(W,\varphi,\xi_{\varrho,s})
=
{\tilde{\cal B}}(W,\varphi,M(w;\xi_{\varrho,s})). 
$$
Here $W$ is a function in the $\psi$ Whittaker model of $\pi$, $\varphi$
is a Schwartz function on $k^n$, and $\xi_{\varrho,s}$ is a 
holomorphic section
for $J(\varrho,s)$, the representation of $\Sp_{2l}(k)$ induced from 
the Siegel parabolic subgroup $Q_l$ and the representation 
$\varrho|{\rm det} |^{s-{1\over2}}$. The section $\xi_{\varrho,s}$ takes values
in an appropriate Whittaker model of~$\varrho$ and  $M(w;\cdot)$ is the local
intertwining operator attached to the Weyl element 
$$
w=\left(\begin{array}{cc}
0&I_l\\
-I_l&0\end{array}\right).
$$
In [GRS2], we defined 
$$
{\cal B}(W,\varphi,\xi_{\varrho,s})
=
\int_{U(k)\backslash\Sp_{2n}(k)}W(g)J_{\psi,n}(\ome_\psi(g)\varphi,j(g)\xi_{\varrho,s})dg
$$
where $U$ is the standard maximal unipotent subgroup of $\Sp(2n)$. 
The precise form of $J_{\psi,n}$ (a Fourier-Jacobi model) is given [GRS2,(1.11)], 
with $\psi^{-1}$ replacing $\psi$. Here $j(g)$ is an appropriate embedding of $g$ 
in $\Sp_{2l}(k)$ and 
$J_{\psi,n}(\varphi,\xi_{\varrho,s})$ is given by an integral 
which stabilizes
on large compact open subgroups of a certain unipotent subgroup.
In particular it
is holomorphic in $s$, and actually, this is a polynomial 
in $q^{-s}$ ($q$ is the number
of elements in the residue field of $k$). The right-hand side of the 
functional equation
has a similar form, with $M(w;\xi_{\varrho,s})$ replacing $\xi_{\varrho,s}$.

\demo{{\rm 3.2.} Poles of local gamma factors} 
We study here the relation between the existence of possible poles of the 
twisted local gamma factors and the structure of the given irreducible 
admissible generic representations. 

By the subquotient theorem of Jacquet \cite{J} or the classification of irreducible 
admissible generic representation of $\GL_{n}(k)$ (\cite{M}, \cite{BZ} and \cite{Z}), 
an irreducible admissible generic representation $\tau$ of $\GL_{n}(k)$ is a 
subquotient of a certain induced representation. More precisely, the representation 
$\tau$ is a subquotient of the unitarily induced representation 
$$
\Ind^{\GL_n(k)}_{P_{n_1,\cdots,n_r}(k)}
(\tau_1|\cdot|^{z_1}\otimes\tau_2|\cdot|^{z_2}\otimes
\cdots\otimes\tau_r|\cdot|^{z_r}).
$$
For the sake of convenience, we use the notation from \cite{BZ} and denote it 
symbolically by 
\begin{equation}
\tau\prec 
\tau_1|\cdot|^{z_1}\times\tau_2|\cdot|^{z_2}\times
\cdots\times
\tau_r|\cdot|^{z_r}
\end{equation}
where $n=\sum_{i=1}^rn_i$ and $\tau_i$'s are irreducible unitary 
supercuspidal representations of $\GL_{n_i}(k)$. 
Without loss of generality, 
we may assume that $z_i$'s are real numbers satisfying the condition:
\begin{equation}
z_1\geq z_2\geq \cdots\geq z_r,
\end{equation}
because of the assumption of the unitarity of the $\tau_i$'s. 
$(P_{n_1,\cdots,n_r}, \tau_1\otimes\cdots\otimes\tau_r)$ is called the 
supercuspidal support of $\tau$, which is determined by $\tau$ uniquely up to 
permutation, and $(z_1,\cdots,z_r)$ is called the exponent of $\tau$. 
\enddemo

\specialnumber{3.1}\proclaim{Proposition}
Assume that an irreducible admissible generic representation $\tau$ of 
$\GL_n(k)$ has the supercuspidal support 
$(P_{n_1,\cdots,n_r}, \tau_1\otimes\cdots\otimes\tau_r)$ 
and the exponent $(z_1,\cdots,z_r)$. 
Let $\varrho$ be an irreducible unitary supercuspidal representation of 
$\GL_l(k)$. 
Then the only possible real poles of the twisted local
gamma factor $\gamma(\tau\times\varrho,s,\psi)$ occur
at $s=1-z_i$ and the only possible real zeros occur at
$s=-z_i$. 
In either case{\rm ,} we have $\varrho\cong \tau^\vee_i$ for some $i\in \{1\nek r\}$.
\endproclaim

\demo{Proof}
By the multiplicativity theorem for the twisted local gamma factors 
(\cite{JP-SS}), 
\begin{equation}
\gamma(\tau\times\varrho,s,\psi)
=
\prod_{i=1}^r\gamma(\tau_i\times\varrho,s+z_i,\psi).
\end{equation}

For each factor $\gamma(\tau_i\times\varrho,s+z_i,\psi)$, 
$$
\gamma(\tau_i\times\varrho,s+z_i,\psi)
=
\epsilon(\tau_i\times\varrho,s+z_i,\psi)\cdot
{L(\tau_i^\vee\times\varrho^\vee,1-(s+z_i))\over 
L(\tau_i\times\varrho,s+z_i)}.
$$
It follows that the local gamma factor $\gamma(\tau_i\times\varrho,s+z_i,\psi)$ 
may have a possible real pole at $s=1-z_i$ and a possible real zero at $s=-z_i$. 
If one of these two things occurs, then we must have 
\vglue12pt
\hfill $
\varrho\cong\tau_i^\vee.
$ 
\enddemo
\vglue12pt

This has the following useful consequence. 

\specialnumber{3.1}\proclaim{{C}orollary}
For any irreducible unitary supercuspidal representation $\varrho$ of 
$\GL_l(k)$, $s=1-z_r$ is the rightmost possible real pole of the product 
$$
\prod_{i=1}^r\gamma(\tau_i\times\varrho,s+z_i,\psi). 
$$
If it occurs in one of the factors in the product{\rm ,} it cannot be cancelled 
by the possible zeros of other factors. In this case{\rm ,} one must have 
that the representation $\varrho$ is equivalent to one of the $\tau_j^\vee$'s 
with $z_j=z_r$. 
\endproclaim

We can prove the same result for $\SO(2n+1)$. First we prove the following 
lemma. 

\specialnumber{3.1}\proclaim{Lemma}
Let $\sigma$ be an irreducible generic supercuspidal representation of 
$\SO_{2n+1}(k)$ and 
$\varrho$ be an irreducible supercuspidal representation of $\GL_l(k)$. 
Assume that the twisted local gamma factor 
$\gamma(\sigma\times\varrho,s,\psi)$ has a pole at 
$s=~1$. Then the pole must be simple and the representation 
$\varrho$ must be self\/{\rm -}\/dual.  Moreover{\rm ,}
$L(\varrho,\Lam^2,s)$ has a pole at $s=0$.
\endproclaim

\demo{Proof}
By the definition of $\gamma(\sigma\times\varrho,s,\psi)$,  
\begin{eqnarray*}
&&\hskip-24pt\gamma(\sigma\times\varrho,s,\psi)\cdot {\cal A}(W_\sigma,\xi_{\varrho,s})
\\ 
&&\quad =\
\gamma(\varrho,\Lambda^2,2s-1,\psi)\cdot 
{\tilde{\cal A}}(W_\sigma,\xi_{\varrho,s})\\
&&\quad =\
\epsilon(\varrho,\Lambda^2,2s-1,\psi)\cdot
{L(\varrho^\vee,\Lambda^2,2(1-s))\over L(\varrho,\Lambda^2,2s-1)}
\cdot
{\tilde{\cal A}}(W_\sigma,\xi_{\varrho,s})\\
&&\quad =\
\epsilon(\varrho,\Lambda^2,2s-1,\psi)\cdot
L(\varrho^\vee,\Lambda^2,2(1-s))\cdot
{{\tilde{\cal A}}(W_\sigma,\xi_{\varrho,s})\over L(\varrho,\Lambda^2,2s-1)}.
\end{eqnarray*}

Since both $\sigma$ and $\varrho$ are irreducible, generic and supercuspidal, 
${\cal A}(W_\sigma,\xi_{\varrho,s})$ is entire (as a function in $s$). 
By Theorem 5.1 in \cite{CS}, the normalized local intertwining operator 
$$
L(\varrho,\Lambda^2,2s-1)^{-1}\cdot M(w_l,\cdot)
$$
is entire and hence 
$$
{{\tilde{\cal A}}(W_\sigma,\xi_{\varrho,s})\over L(\varrho,\Lambda^2,2s-1)}
$$ 
is entire. (Actually, at the point of our interest, $s=1$,
the intertwining integral converges.) Further, we can choose 
certain data such that 
${\cal A}(W_\sigma,\xi_{\varrho,s})$ does not vanish. Hence, if 
$\gamma(\sigma\times\varrho,s,\psi)$ has a pole at $s=1$, 
then we must have that 
$L(\varrho^\vee,\Lambda^2,s)$ has a pole at $s=0$. It is known that 
if $L(\varrho,\Lambda^2,s)$ has a pole at $s=0$, then 
$$
\varrho\cong \varrho^\vee
$$
and the pole is simple. In other words, the gamma factor 
$\gamma(\sigma\times\varrho,s,\psi)$ has at most a simple pole at $s=1$ and 
if the pole occurs, then $\varrho$ is self-dual. 
\enddemo

Here is an analogue   needed on the metaplectic side. 

\specialnumber{3.2}\proclaim{Lemma}
Let $\pi$ be an irreducible{\rm ,} supercuspidal representation 
of $\tilsp_{2n}(k)$. Assume that 
it is $\psi_{\tilU,1}$ generic. Assume that $l>n$ and let $\varrho$ be  
a representation of $\GL_l(k)$
induced from $\varrho_1\otimes\cdots\otimes\varrho_r${\rm ,} where $\varrho_i$ 
is an irreducible unitary supercuspidal
representation of $\GL_{m_i}(k)$ with $m_1+ \cdots +m_r=l$. Assume also that 
$\varrho_i$ is not isomorphic to  $\varrho_j^\vee$ for $i\not=j$. 
Then the order of the pole at $s=1$ of 
$\gamma(\pi\times\varrho,s,\psi)$ is less than{\rm ,} or equal to that of 
$\prod_{i=1}^r L(\varrho_i,\Lambda^2,2(1-s))$. In particular{\rm ,} if 
$L(\varrho_i,\Lambda^2,z)$ is holomorphic at $z=0$ {\rm (}e.g.{\rm ,} $\varrho_i$
is not self\/{\rm -}\/dual\/{\rm ),} for $1\leq i\leq r${\rm ,} 
then $\gamma(\pi\times\varrho,s,\psi)$ is holomorphic at $s=1$.
\endproclaim

The proof follows (as in Lemma 3.1) from the local functional equation, 
and the fact that Whittaker functions
for $\pi$ are compactly supported modulo $\tilU_n$ 
(we assume that $\pi$ is supercuspidal.) 
The details are written in Proposition~1 and Corollary~1 of [GRS6]. Note also that
$$
L(\varrho^\vee,\Lambda^2,z)=\prod_{1\leq i<j\leq r}
L(\varrho_i^\vee\times\varrho_j^\vee,z)
\prod_{i=1}^r L(\varrho_i^\vee,\Lambda^2,z)=
\prod_{i=1}^r L(\varrho_i^\vee,\Lambda^2,z).
$$
We use this for $z=2(1-s)$. 

Returning to the orthogonal group we consider more general cases. 
Let $\sigma$ be an irreducible admissible generic representation 
of $\SO_{2n+1}(k)$. By Jacquet's subquotient theorem \cite{J}, there is a standard 
parabolic subgroup $Q$ with its Levi part isomorphic to 
$$
\GL_{m_1}\times\ldots\times\GL_{m_r}\times\SO_{2m_0+1},
$$
($n=\sum_{i=0}^rm_i$) 
and there are irreducible unitary supercuspidal representations $\tau_i$ of 
$\GL_{m_i}(k)$ ($i=1,2,\cdots,r$) and an irreducible generic supercuspidal 
representation $\sig_0$ of $\SO_{2m_0+1}(k)$, such that the representation $\sig$ is 
a subquotient of the unitarily induced representation 
\begin{equation}
\sigma\prec 
\Ind^{\SO_{2n+1}}_{Q}(\tau_1|\det|^{z_1}\otimes\cdots\otimes\tau_r|\det|^{z_r}
\otimes\sigma_0).
\end{equation}
Without loss of generality, 
we may assume that the parameters $z_i$ are real and 
have the property that $z_1\geq z_2\geq\cdots\geq z_r\geq 0$. 
With this assumption, 
we say that the representation $\sigma$ has supercuspidal support 
$(Q;\tau_1,\tau_2,\cdots,\tau_r;\sigma_0)$ and exponents 
$(z_1,z_2,\cdots,z_r)$. 

\specialnumber{3.2}\proclaim{Proposition}
Let $\sigma$ be an irreducible admissible generic representation of 
$\SO_{2n+1}(k)$ with 
supercuspidal support $(Q;\tau_1,\tau_2,\cdots,\tau_r;\sigma_0)$ and 
exponents $(z_1,z_2,\cdots,z_r)$. 
Then $s=1+z_1$ is the rightmost real point at which the twisted 
local gamma factors $\gamma(\sigma\times\varrho,s,\psi)$ can possibly have 
a pole for any irreducible unitary supercuspidal representation $\varrho$ of 
$\GL_l(k)$ where 
$l$ is any positive integer. If the pole at $s=1+z_1$
occurs for some $(l,\varrho)$\/{\rm ,}\/ then 
$$
\varrho\cong\tau_{i_0}
$$
where $\tau_{i_0}$ is a representation among the $\tau_i$'s such that 
$z_{i_0}=z_1$. 
\endproclaim

\demo{Proof}
By the multiplicativity of the local gamma factors (Theorem 3.1), 
\begin{eqnarray*}
 \gamma(\sigma\times\varrho,s,\psi) 
&\hskip-8pt=\hskip-8pt&
\omega_\varrho(-1)^{rn+m_0}\\ &\hskip-8pt\hskip-8pt& \cdot\
\left[\prod_{i=1}^r\gamma(\tau_i\times\varrho,s+z_i,\psi)
\gamma(\tau^\vee_i\times\varrho,s-z_i,\psi)\right]
\gamma(\sigma_0\times\varrho,s,\psi)
\end{eqnarray*}
where $\varrho$ is any irreducible unitary supercuspidal representation of 
$\GL_l(k)$. 

By the argument in the proof of Proposition 3.1, we get that 
\begin{itemize}
\item[(1)] the factor $\gamma(\tau_i\times\varrho,s+z_i,\psi)$ may contribute 
a possible real pole at\break $s=1-z_i$ and a possible real zero at $s=-z_i$, and 
if one of these two things occurs, then $\varrho\cong\tau_i^\vee$, 
\item[(2)] the factor $\gamma(\tau^\vee_i\times\varrho,s-z_i,\psi)$ may 
contribute 
a possible real pole at $s=1+z_i$ and a possible real zero at $s=z_i$, and 
if one of these two things occurs, one has $\varrho\cong\tau_i$, and 
\item[(3)] the factor $\gamma(\sigma_0\times\varrho,s,\psi)$ has no zero for 
${\rm Re}(s)>0$ (\cite[\S 5]{Sh1} and \cite[Prop.~7.2]{Sh2}) 
and may have a possible simple real pole at $s=1$ (Lemma 3.1). 
If the pole occurs, we must have
$\varrho\cong\varrho^\vee$ (and $L(\varrho,\Lam^2,s)$
has a pole at $s=0$).
\end{itemize}
Hence $s=1+z_1$ is the rightmost possible real pole of the product 
$$
\omega_\varrho(-1)^{rn+m_0}
[\prod_{i=1}^r\gamma(\tau_i\times\varrho,s+z_i,\psi)
\gamma(\tau^\vee_i\times\varrho,s-z_i,\psi)]
\gamma(\sigma_0\times\varrho,s,\psi). 
$$
If the pole occurs at $s=1+z_1$, it cannot be cancelled by any possible 
zero from other factors in the product and $\varrho\cong\tau_{i_0}$,
where $\tau_{i_0}$ is a representation among the $\tau_i$'s, such that 
$z_{i_0}=z_1$. 
\enddemo

\specialnumber{3.2}\proclaim{{C}orollary}
Let $\sigma$ and $\sigma'$ be irreducible admissible generic representations 
of $\SO_{2n+1}(k)$ with supercuspidal supports 
$(Q;\tau_1,\tau_2,\cdots,\tau_l;\sigma_0)$ 
and $(Q';\tau'_1,\tau'_2,\cdots,\tau'_{l'};\sigma'_0)${\rm ,} and exponents 
$(z_1,z_2,\cdots,z_l)$ and $(z'_1,z'_2,\cdots,z'_{l'})${\rm ,} respectively. 
If the twisted gamma factors are the same{\rm ,} i.e.\ 
$$
\gamma(\sigma\times\varrho,s,\psi)=\gamma(\sigma'\times\varrho,s,\psi)
$$
for all irreducible supercuspidal representations $\varrho$ of $\GL_l(k)$ with 
$l=1,2,\cdots,\break 2n-1${\rm ,} then $z_1=z_1'$. 
\endproclaim

\demo{{\rm 3.3.} Gamma factors and functorial lift}
The Langlands functorial lift (or transfer) conjecture describes the relation 
between 
automorphic representations of two different groups as long as their Langlands 
dual groups have an `admissible' relation. One may find more details about the 
Langlands conjectures in \cite{B}. In this paper, we need a special case, which we 
describe below in more detail. 

Let $F$ be a number field. The Langlands dual group of $\SO(2n+1)$ is $\Sp_{2n}(\CC)$ 
and the Langlands dual group of $\GL(2n)$ is $\GL_{2n}(\CC)$. The natural embedding 
of $\Sp_{2n}(\CC)$ into $\GL_{2n}(\CC)$ is `admissible', so that by Langlands 
functorial lift conjecture, any irreducible automorphic representation 
$\Sig$ of $\SO_{2n+1}(\AA)$ can be lifted to an irreducible  
automorphic representation $\cal T$ of $\GL_{2n}(\AA)$, functorially.  
In other words, if we write 
$$
\Sig=\otimes_v\Sig_v,\ \ \ {\cal T}=\otimes_v{\cal T}_v,
$$
(as the restricted tensor product of the local components), then ${\cal T}$ is a 
(global) Langlands functorial lift of $\Sig$ if and only if for each local place $v$, 
the local component ${\cal T}_v$ is a (local) Langlands functorial lift of the 
local component $\Sig_v$. Moreover, a lift from $\Sig$ to $\cal T$ is called 
a weak Langlands functorial lift if the local component ${\cal T}_v$ is a (local) 
Langlands functorial lift of the local component $\Sig_v$ for all archimedean
places and for almost all  
places $v$ of $F$,where ${\cal T}_v$ and $\Sig_v$ are unramified. 

In \cite{CKP-SS} the weak Langlands functorial lift from $\SO_{2n+1}(\AA)$ to 
$\GL_{2n}(\AA)$ was proved to exist for irreducible generic cuspidal automorphic 
representations of $\SO_{2n+1}(\AA)$, by using the converse theorem for
$\GL_{2n}$. 
 
Let $\Sigma$ be an irreducible generic cuspidal automorphic representation of 
$\SO_{2n+1}({\Bbb A})$ and ${\cal T}$ be a weak Langlands functorial lifting 
to $\GL_{2n}({\Bbb A})$ of $\Sigma$. The following theorem determines 
the explicit structure of the image of the weak Langlands functorial lifting 
without any assumption (which improves the results in \cite{CKP-SS}). 
\enddemo

\proclaimtitle{\cite{GRS5}}
\specialnumber{3.2}\proclaim{Theorem} 
Let ${\cal T}$ be an irreducible{\rm ,} automorphic representation of  
$\GL_{2n}({\Bbb A})$ which is the weak Langlands functorial lifting  
of an irreducible generic cuspidal automorphic representation $\Sigma$ of 
$\SO_{2n+1}({\Bbb A})${\rm ,} then ${\cal T}$ is  generic and self\/{\rm -}\/dual. 
Moreover{\rm ,} ${\cal T}$ is isomorphic to 
$$
\Ind^{\GL_{2n}({\Bbb A})}_{P_{2n_1\nek 2n_r}({\Bbb A})}
({\cal T}_1\otimes\cdots\otimes{\cal T}_r)
$$
where $P_{2n_1\nek 2n_r}$ is the standard parabolic subgroup of $\GL_{2n}$ 
corresponding to 
the partition $2n=\sum_{i=1}^r2n_i${\rm ;} ${\cal T}_i$ for $i=1,2,\cdots,r${\rm ,} is 
an irreducible{\rm ,} unitary{\rm ,} self\/{\rm -}\/dual{\rm ,} cuspidal{\rm ,} automorphic representation of 
$\GL_{2n_i}({\Bbb A})$ such that the partial exterior square $L$\/{\rm -}\/function 
$L^{\cal S}({\cal T}_i,\Lambda^2,s)$ has 
a pole at $s=1${\rm ,} and ${\cal T}_i\not\cong {\cal T}_j$ if $i\neq j$. 
In particular\/{\rm ,}\/ ${\cal T}$ is uniquely determined up to isomorphism by 
$\Sigma$. 
\endproclaim

Moreover, one has the following result on local gamma factors. 

\specialnumber{3.3}\proclaim{Proposition}
Let $\Sigma_v$ be the $v$\/{\rm -}\/local component of an irreducible generic cuspidal 
automorphic representation $\Sigma$ of $\SO_{2n+1}({\Bbb A})$ and let ${\cal T}$ be the weak 
Langlands functorial lifting to $\GL_{2n}({\Bbb A})$ of $\Sigma$. Then 
for every supercuspidal representation $\tau_v$ of $\GL_l(F_v)$ where $l$ is any positive integer\/{\rm ,} 
one has 
$$
\gamma(\Sigma_v\times\tau_v,s,\psi_v)
=
\gamma({\cal T}_v\times\tau_v,s,\psi_v).
$$
\endproclaim

It is the result of Corollary 5 in \cite{CKP-SS} that the identity
$$
\gamma(\Sigma_v\times\tau_v,s,\psi_v)
=
\gamma({\cal T}_v\times\tau_v,s,\psi_v)
$$
holds for every supercuspidal representation $\tau_v$ of $\GL_l(F_v)$ with 
$l=1,2,\cdots,\break 2n-1$. The argument in  {\it loc. cit.}\ is valid for any $l$, with
no restriction.

\demo{{\rm 3.4.} Structure of the image of local functorial lifting}
Let $\sigma$ be an irreducible generic supercuspidal representation of 
$\SO_{2n+1}(F_{v_0})$ and let 
$\Sigma$ be an irreducible generic cuspidal automorphic representation of 
$\SO_{2n+1}({\Bbb A})$ as constructed in Proposition 2.5. Then 
$$
\Sigma_{v_0}\cong\sigma.
$$
Let ${\cal T}$ be the image of $\Sigma$ under the weak Langlands functorial 
lifting. 
Then, by Proposition 3.3, 
\begin{equation}
\gamma(\Sigma_{v_0}\times\varrho,s,\psi_{v_0})
=
\gamma({\cal T}_{v_0}\times\varrho,s,\psi_{v_0}),
\end{equation}
for all irreducible supercuspidal representations $\varrho$ of 
$\GL_l(F_{v_0})$ ($l$ any positive integer). 

We shall determine the explicit structure of the $v_0$-local component\break
$\tau:={\cal T}_{v_0}$ in terms of the supercuspidal support by using the 
existence of poles of the local gamma factors, the global version of which 
was given in Theorem 3.2.

By the subquotient theorem or the classification of irreducible  
generic representations of $\GL_{2n}(k)$ ($k=F_{v_0}$) (\cite{M}, \cite{BZ} 
and \cite{Z}), the irreducible admissible generic self-dual representation 
$\tau$ of $\GL_{2n}(k)$ is 
a subquotient of an induced representation 
\begin{equation}
\tau\prec 
\tau_1|\cdot|^{z_1}\times
\cdots\times
\tau_r|\cdot|^{z_r}\times
\eta_1\times\cdots\times\eta_t\times
\tau_r^\vee|\cdot|^{-z_r}\times
\cdots\times
\tau_1^\vee|\cdot|^{-z_1},\;
\end{equation}
where $\tau_i$'s are irreducible unitary supercuspidal representations of 
$\GL_{m_i}(k)$, 
and 
$\eta_j$'s are irreducible unitary supercuspidal self-dual representations of\break
$\GL_{2n_j}(k)$ ($n=\sum_{i=1}^rm_i+\sum_{j=1}^tn_j$). From the given data, 
we may also assume that if $z_i=0$ then $\tau_i$ is not self-dual. 
Without loss of generality, 
we may assume that $z_i$'s are real numbers satisfying the condition:
\begin{equation}
z_1\geq z_2\geq \cdots\geq z_r\geq 0.
\end{equation}

\specialnumber{3.3}\proclaim{Theorem}
Let $\sigma$ be an irreducible generic supercuspidal representation of 
$\SO_{2n+1}(k)$. 
There exists a unique irreducible generic representation $\tau$ of 
$\GL_{2n}(k)$ such that 
\begin{equation}
\gamma(\sigma\times\varrho,s,\psi)
=
\gamma(\tau\times\varrho,s,\psi),
\end{equation}
for all irreducible supercuspidal representations $\varrho$ of $\GL_l(k)$  
{\rm (}$l$ any positive integer\/{\rm ).} Moreover{\rm ,} in the notation of $(3.10)${\rm ,} $\tau$  
must have the following properties\/{\rm :}
\begin{itemize}
\item[{\rm (1)}] $\tau=\eta_1\times\cdots\times\eta_t${\rm ,} where $\eta_j$ are 
irreducible 
unitary supercuspidal self\/{\rm -}\/dual representations of $\GL_{2n_j}(k)$ 
{\rm (}$2n=\sum_{j=1}^t2n_j${\rm )} 
and $\eta_i\neq \eta_j$ if $i\neq j${\rm ;}
\item[{\rm (2)}] the local $L$\/{\rm -}\/function $L(\eta_i,\Lambda^2,s)$ has a pole at $s=0$ 
for $i=1,2,\cdots,t$.
\end{itemize}

\endproclaim

\demo{Proof}
The existence of $\tau$ is given by the weak functorial lifting described 
above ($\tau={\cal T}_{v_0}$ in (3.9)) and the uniqueness of $\tau$ follows 
from the local converse theorem for 
$\GL_n$ (\cite{Hn2}). 
It remains to determine the structure of $\tau$
explicitly.  Suppose that $\tau$ has supercuspidal
support as in (3.10). Let $\varrho$ be an irreducible,
supercuspidal representation of $\GL_l(k)$.
By Lemma 3.1, we know that the gamma factor 
$\gamma(\sigma\times\varrho,s,\psi)$ has at most a simple
pole at $s=1$ and if the pole occurs, then $\varrho$ is
self-dual.  By the identity in (3.12), we know that the
local gamma factor $\gamma(\tau\times\varrho,s,\psi)$ has
at most a simple pole at $s=1$ and if the pole occurs,
then $\varrho$ is self-dual. We are going to show that
this information about the existence of a pole of 
$\gamma(\tau\times\varrho,s,\psi)$ at 
$s=1$ and the order of the pole control the structure of $\tau$. 

By using the multiplicativity theorem for the twisted local gamma factors 
in this case (Theorem 3.1 in \cite{JP-SS}), we have 
\begin{eqnarray*}
& &\hskip-16pt\gamma(\tau\times\varrho,s,\psi)\\
& &=
\left[\prod_{i=1}^r
\gamma(\tau_i\times\varrho,s+z_i,\psi)\cdot
\gamma(\tau_i^\vee\times\varrho,s-z_i,\psi)\right]\cdot
\prod_{j=1}^t
\gamma(\eta_j\times\varrho,s,\psi). 
\end{eqnarray*}
By Corollary 3.1, we know that $1+z_1$ is the rightmost
possible real pole of the twisted local gamma factor
$\gamma(\tau\times\varrho,s,\psi)$ and if it occurs, it
cannot be cancelled by the zeros of other twisted local 
gamma factors in the product.

Now, we take $\varrho=\tau_1$. Then the local gamma factor
$\gamma(\tau\times\tau_1,s,\psi)$ has a pole at
$s=1+z_1$.  By the identity in (3.12), we know that the
local gamma factor $\gamma(\sigma\times\tau_1,s,\psi)$ 
has a pole at $s=1+z_1$. Repetition of the argument of Lemma~3.1  implies that the local exterior square
$L$-function $L(\tau_1,\Lambda^2,s)$ has a pole at 
$s=-2z_1$. Since $\tau_1$  is unitary (supercuspidal), we
get that ${\rm Re}(2z_1)=0$, and hence $z_1=0$. We know by Lemma 3.1
that the  pole is simple, and the representation $\tau_1$
is self-dual.  Now, this implies that  
$$z_1=\cdots =z_r=0$$
and the local gamma factor 
$\gamma(\tau\times\tau_1,s,\psi)$ can be expressed as 
$$
\gamma(\tau\times\tau_1,s,\psi)
=[\gamma(\tau_1\times\tau_1,s,\psi)]^2\times [\cdots].
$$
Hence the local gamma factor $\gamma(\tau\times\tau_1,s,
\psi)$ has a pole at $s=1$ with order two or higher.
This contradicts the simplicity of the pole of the local
gamma factor 
$\gamma(\sigma\times\tau_1,s,\psi)$ because of (3.12). 
Therefore we conclude that 
the supercuspidal representation $\tau_1$ should not 
occur in the supercuspidal support of $\tau$. 

By the same argument, we can conclude that all
supercuspidal representations 
$\tau_i$ with $i=1,2,\cdots,r$, do not occur in the
supercuspidal support of $\tau$. Namely we have that 
$$
\tau\prec 
\eta_1\times\cdots\times\eta_t.
$$
Since the induced representation 
$\eta_1\times\cdots\times\eta_t$ is irreducible, we conclude that 
$$
\tau\cong 
\eta_1\times\cdots\times\eta_t.
$$
Since the pole at $s=1$ of the gamma factor 
$\gamma(\sigma\times\varrho,s,\psi)$ is at most simple, 
the representations $\eta_j$ in the supercuspidal support of $\tau$ 
must be all distinct. Finally the existence of the pole of the 
local $L$-functions $L(\eta_i,\Lambda^2,s)$ at $s=0$ for $i=1,2,\cdots,t$ is 
now clear from the argument above (see Lemma 3.1 again).
\enddemo 

\specialnumber{3.3}\proclaim{{C}orollary}
Let $\Sigma$ be an irreducible generic cuspidal automorphic representation of 
$\SO_{2n+1}({\Bbb A})$ and ${\cal T}$ be the weak Langlands functorial lifting 
of $\Sigma$ from $\SO_{2n+1}({\Bbb A})$ to $\GL_{2n}({\Bbb A})$. 
For a finite local place $v$ 
of the number field~$F${\rm ,} 
if the $v$\/{\rm -}\/local component $\Sigma_v$ is supercuspidal{\rm ,} then 
the corresponding $v$\/{\rm -}\/local component ${\cal T}_v$ has form 
$$
{\cal T}_v=\eta_1\times\cdots\times\eta_{t_v}
$$
where the $\eta_j$\/{\rm '}\/s are irreducible unitary self\/{\rm -}\/dual supercuspidal 
representations of $\GL_{2n_j}(F_v)$ such that 
\begin{itemize}
\item[{\rm (1)}] $\eta_i\neq \eta_j$ if $i\neq j$ {\rm (}$2n=\sum_{j=1}^{t_v}2n_j${\rm );} 
\item[{\rm (2)}] the local $L$\/{\rm -}\/function $L(\eta_i,\Lambda^2,s)$ has a pole at $s=0$ 
for $i=1,2,\cdots,t_v$.
\end{itemize}

\endproclaim
\phantom{dddd}
\vglue-24pt 
3.5. {\it Gamma factors and Howe duality}.
We discuss here the relation of local gamma factors under the local Howe 
duality. Of course, these gamma factors should remain invariant, 
but here we need much less. We need just the preservation of a pole at $s=1$,
as follows. 

\specialnumber{3.4}\proclaim{Proposition}
Let $\sigma$ be an irreducible generic supercuspidal representation of 
$\SO_{2n+1}(k)$ and let $\pi$ be the irreducible
$\psi_{\tilU,1}$\/{\rm -}\/generic supercuspidal representation 
of $\tilspn(k)${\rm ,} which is the local $\psi$\/{\rm -}\/Howe lift of $\sigma$. 
Let $\tau$ be the local lift of $\sigma$ to $\GL_{2n}(F)${\rm ,} i.e.\ the 
representation given by Theorem {\rm 3.3.} Write 
$\tau=\eta_1\times\cdots\times\eta_t$ as in Theorem {\rm 3.3.} 
Then $\gamma(\pi\times\tau,s,\psi)$ has a pole of order $t$ at $s=1$.
\endproclaim

\demo{Proof}
Let $\Sigma$ and $\Pi$ be the representations constructed in Proposition 2.5. 
We keep the notation used prior to Proposition 2.5. (We may and do assume 
that $\psi_{0,\nu_0}=\psi$.) Recall the finite set of places $S_0$, and
the choices of $(\psi_{0,\nu})_{\tilde{U_\nu}}$ generic supercuspidal
representations $\pi_\nu$ of $\tilspn(k_\nu)$, 
such that $\Pi_{\nu_0}\cong\pi$, and $\Pi_\nu\cong\pi_\nu$, for each 
place $\nu$ in $S_0$. Recall also that $\Sigma$ is 
the global theta lift of $\Pi$ with respect to $\psi_0$, and, in particular,
$\Sig_\nu$ is the local $\psi_{0,\nu}$- Howe lift of $\pi_\nu$, for
each $\nu$ in $S_0$ (and hence $\Sig_{\nu_0}\cong\sig$). Recall 
that $\Sig$ is globally generic (and hence 
$\Sig_\nu$ is generic, for each $\nu$ in $S_0$).

Let ${\cal T}$ be an irreducible representation of $\GL_{2n}({\Bbb A})$ of 
the form  
$$
\Ind^{\GL_{2n}({\Bbb A})}_{P_{m_1\nek m_r}({\Bbb A})}
({\cal T}_1\otimes\cdots\otimes{\cal T}_r)
$$
where ${\cal T}_i$ is an irreducible, automorphic, cuspidal, unitary
representation of $\GL_{m_i}({\Bbb A})$, $1\leq i\leq r$, $m_1+\cdots+m_r=2n$.
Note that ${\cal T}$ is irreducible. We realize it automorphically by taking
an Eisenstein series induced from 
${\cal T}_1 |{\rm det}|^{s_1}\otimes\cdots\otimes{\cal T}_r |{\rm det}|^{s_r}$, and 
evaluating it at $(s_1\nek s_r)=(0\nek 0)$.

Let $S$ be a finite set of finite places, containing $S_0$, such that 
at finite places outside $S$, $\Pi$ (hence $\Sig$ ) and ${\cal T}$ are
unramified. For such places $\nu$,
$$
\gamma(\Sig_\nu \times {\cal T}_\nu,s,\psi_{0,\nu})=
\gamma(\Pi_\nu \times {\cal T}_\nu,s,\psi_{0,\nu}).
$$

The functional equation satisfied by the global integrals for 
$\SO_{2n+1}\times \GL_{2n}$ and for $\tilspn\times \GL_{2n}$,  
means that the
product over all places of the corresponding gamma factors is one. From 
the last equalities outside $S$, we conclude
$$
\gamma_\infty(\Sig\times {\cal T},s,\psi_0)
\prod_{\nu\in S}\gamma(\Sig_\nu\times {\cal T}_\nu,s,\psi_{0,\nu})=
\gamma_\infty(\Pi\times {\cal T},s,\psi_0)
\prod_{\nu\in S}\gamma(\Pi_\nu\times {\cal T}_\nu,s,\psi_{0,\nu}). 
$$
Here $\gamma_\infty$ denotes the product of local gamma factors over all 
archimedean places. Write $S= \cup_{i=1}^{i=N} S(i)$, where $S(i)$ is the
set of all places $\nu$ in $S$, such that the residual characteristic 
of $k_\nu$ is $p_i$, and $p_1\nek p_N$ are different prime numbers. Put
$$
A_i(s)= \prod_{\nu\in S(i)}{\gamma(\Pi_\nu\times {\cal T}_\nu,s,\psi_{0,\nu})
\over\gamma(\Sig_\nu\times {\cal T}_\nu,s,\psi_{0,\nu})}.
$$ 
Thus, $A_i(s)$ is in ${\Bbb C}(p_i^{-s})$. Write $A_i(s)= R_i(p_i^{-s})$, 
where $R_i(x)\in{\Bbb C}(x)$. We get
$$
\prod_{i=1}^{i=N}R_i(p_i^{-s})=
{\gamma_\infty(\Sig\times {\cal T},s,\psi_0)\over\gamma_\infty(\Pi\times 
{\cal T},s,\psi_0)}.
$$
We show in Section 7.2 that there is $M>0$, such that for $|\Im(s)|>M$, both
$\gamma_\infty(\Sig\times {\cal T},s,\psi_0)$ and
$\gamma_\infty(\Pi\times {\cal T},s,\psi_0)$ are holomorphic with no zeroes.
Thus, if $s_0$ is a pole (resp.\ a zero) of $R_j(p_j^{_s})$, then
$s_0+{2\pi mi\over \log p_j}$ is a pole (resp.\ a zero) of
$R_j(p_j^{-s})$, for
all integers $m$. Since $p_1\nek p_N$ are different, we see that
$\prod_{i=1}^{i=N}R_i(p_i^{-s})$ has an unbounded sequence of 
poles (resp.\ zeroes) on the line ${\rm Re}(s)=s_0$. This is impossible unless
each $R_i(x)$ is an exponential. Thus there are $a_i $ and $b_i$, such
that
$$
a_ie^{b_is} \prod_{\nu \in S(i)}\gamma(\Sig_\nu\times{\cal T}_\nu,s,\psi_{0,\nu})
=
\prod_{\nu \in S(i)}\gamma(\Pi_\nu\times{\cal T}_\nu,s,\psi_{0,\nu}).
$$
In particular, there is an exponential $\alp(s)$, such that
$$
\alp(s) \prod_{\nu \in S'_0}\gamma(\Sig_\nu\times{\cal T}_\nu,s,\psi_{0,\nu})
=
\prod_{\nu \in S'_0}\gamma(\Pi_\nu\times{\cal T}_\nu,s,\psi_{0,\nu})
$$
where $S'_0$ is the set of places $\nu$ in $S$, where the residual
characteristic of $k_\nu$ is the same as that of $k_{\nu_0}$ (note that
this is a subset of $S_0$).

Consider now the special case where ${\cal T}_{\nu_0}\cong\tau$, where $\tau$
is the local lift of $\sig$ to $\GL_{2n}(k_{\nu_0})$, as in Theorem 3.3. More
precisely, take $m_i=2n_i$ and ${\cal T}_i$- such that it has a trivial central
character, and its $\nu_0$ component is isomorphic to $\eta_i$. For each 
$i\leq t$, let $\chi_i$ be a unitary automorphic character of the adeles
of $k$, such that $\chi_{i,\nu_0}=1$ and for $\nu\not=\nu_0$ in $S'_0$,
${\cal T}_{i,\nu}\chi_{i,\nu}$ is not isomorphic to the dual of
${\cal T}_{j,\nu}\chi_{j,\nu}$, for all $i,j\leq t$. (Such characters exist by
Theorem 5, p.103 of [AT]. For example, it is enough to guarantee that
$\chi_{i,\nu}^{2n_i}\chi_{j,\nu}^{2n_j}\not=1$, for $\nu\not=\nu_0$ in $S'_0$.)  
Let ${\cal T}(\chi)$ be the representation obtained from ${\cal T}$ by twisting
${\cal T}_i$ by $\chi_i$. Repeating the last equality, we have an exponential
$\alp(s)$, such that
\begin {equation}
\alp(s) \prod_{\nu \in S'_0}\gamma(\Sig_\nu\times{\cal T}(\chi)_\nu,s,
\psi_{0,\nu})
=
\prod_{\nu \in S'_0}\gamma(\Pi_\nu\times{\cal T}(\chi)_\nu,s,\psi_{0,\nu}).\hskip.4in
\end{equation}
We have,
$$
\gamma(\Sig_{\nu_0}\times{\cal T}(\chi)_{\nu_0},s,\psi_{0,\nu_0})
=
\gamma(\sigma\times\tau,s,\psi)=\gamma(\tau\times\tau,s,\psi),
$$
which, by the structure of $\tau$, has a pole of order $t$ at $s=1$. The 
remaining factors of the left-hand side of (3.13) are holomorphic and 
nonzero at $s=1$. This follows from Lemma 3.1, the multiplicativity property,
and the choice of $\chi_i$.  We conclude that the left hand side of (3.13) 
has a pole of order $t$ at $s=1$.  By Lemma 3.2, all factors on the
 right-hand side of (3.13), corresponding to $\nu\not=\nu_0$, are holomorphic 
at $s=1$ (again, due to our choice of $\chi_i$). We conclude from (3.13)
that $\gamma(\pi\times\tau,s,\psi)$ has a pole of order $t$ at $s=1$. 
This completes the proof of the proposition.
\enddemo

\section{The local converse theorem for $\SO(2n+1)$:\\ The
supercuspidal case}

In this section, we prove the local converse theorem for irreducible generic 
supercuspidal representations of $\SO(2n+1)$. 
The idea of the proof is to transfer the local converse theorem for 
$\GL_{2n}(k)$ to $\SO_{2n+1}(k)$ by combining various liftings. 

\proclaimtitle{Local converse theorem for $\SO(2n+1)$: The supercuspidal case}
\proclaim{Theorem} 
Let $\sigma$ and $\sigma'$ be irreducible generic supercuspidal 
representations of $\SO_{2n+1}(k)$. Assume that 
$$\gamma(\sigma\times\varrho,s,\psi)=\gamma(\sigma'\times
\varrho,s,\psi)$$
for all irreducible supercuspidal representations
$\varrho$ of $\GL_l(k)$, $l=1,2,\cdots,2n-1$. Then the
representations $\sigma$ and $\sigma'$ are equivalent.
\endproclaim

\demo{Proof}
Let $\sigma$ and $\sigma'$ be irreducible generic supercuspidal 
representations of 
$\SO_{2n+1}(k)$, such that  
\begin{equation}
\gamma(\sigma\times\varrho,s,\psi)
=
\gamma(\sigma'\times\varrho,s,\psi)
\end{equation}
for all irreducible supercuspidal representations
$\varrho$ of $\GL_l(k)$ with $l=1,\cdots,\break 2n-1$.  Let $\pi$
and $\pi'$ be the local $\psi$ - Howe lifts of $\sigma$
and $\sigma'$ respectively to $\tilspn(k)$ 
(Theorem 2.2). These are irreducible, $\psi_{{\tilde{U}},1}$-generic,
supercuspidal representations of $\tilspn(k)$. 

By Theorem 3.3, there are irreducible generic representations $\tau$ and 
$\tau'$ of $\GL_{2n}(k)$ corresponding to $\sigma$ and $\sigma'$, 
respectively, such that 
\begin{equation}
\gamma(\sigma\times\varrho,s,\psi)
=
\gamma(\tau\times\varrho,s,\psi)
\end{equation}
and 
\begin{equation}
\gamma(\sigma'\times\varrho,s,\psi)
=
\gamma(\tau'\times\varrho,s,\psi)
\end{equation}
for all irreducible supercuspidal representations $\varrho$ of $\GL_l(k)$ with 
any positive integer $l$. In particular,  
$$
\gamma(\tau\times\varrho,s,\psi)
=
\gamma(\tau'\times\varrho,s,\psi)
$$ 
for all irreducible supercuspidal representations $\varrho$ of $\GL_l(k)$ with 
$l=1,2,\cdots,\break 2n-1$. By Theorem 1.1,  
$$
\tau\cong\tau'.
$$
By Theorem 3.3 again,  
\begin{equation}
\tau\cong\tau'\cong \eta_1\times\cdots\times\eta_t
\end{equation}
such that each $\eta_i$ is an irreducible unitary
self-dual supercuspidal representation of $\GL_{2n_i}(k)$,
$\eta_i\not\cong\eta_j$ if $i\neq j$, and the local
exterior square $L$-function $L(\eta_i,\Lambda^2,s)$ has 
a pole at $s=0$ for $i=1,2,\cdots,t$. By Proposition 3.4,
both gamma factors $\gamma(\pi\times \tau,s,\psi)$ and
$\gamma(\pi'\times \tau,s,\psi)$ have a pole of order $t$ at $s=1$. 
By the theory of the local backward 
lifting from $\GL_{2n}(k)$ to $\tilspn(k)$ in
\cite{GRS6}, the image $\pi_\psi(\tau)$ of $\tau$ under
the backward lifting is the unique irreducible 
$\psi_{{\tilde{U}},1}$-generic supercuspidal
representation of $\tilspn(k)$ with the
property that the twisted local gamma factor 
$\gamma(\pi_\psi(\tau)\times\tau ,s,\psi)$ has a
pole of order $t$ at $s=1$. Since  
$\pi$  and $\pi'$  are supercuspidal and
$\psi_{\tilU,1}$-generic, we conclude that 
$$\pi_\psi(\tau)\cong\pi\cong\pi'.$$
Finally, it follows from Theorem 2.2 that the
representation $\sigma$ is equivalent to the
representation $\sigma'$. 
\enddemo

\section{\bf The local converse theorem:  general case}

We shall use the information on poles of twisted gamma
factors to determine explicitly the structure of the supercuspidal support of 
irreducible generic representations.  Then the general case of
the local converse theorem follows from the supercuspidal
case. 

\specialnumber{5.1}\proclaim{Theorem}
Let $\sigma$ and $\sigma'$ be irreducible 
generic representations of $\SO_{2n+1}(k)$ with
supercuspidal supports $(Q;\tau_1,\tau_2,\cdots,\tau_r;
\sigma_0)$ and $(Q';\tau'_1,\tau'_2,\cdots,\break \tau'_{r'};
\sigma'_0)${\rm ,} and exponents $(z_1,z_2,\cdots,z_r)$ and
$(z'_1,z'_2,\cdots,z'_{r'})${\rm ,} respectively {\rm (}\/see
{\rm (3.8)).}  Assume that 
$$\gamma(\sigma\times\varrho,s,\psi)=\gamma(\sigma'
\times\varrho,s,\psi)$$
for all irreducible supercuspidal representations
$\varrho$ of $\GL_l(k)$ with $l=1,2,\cdots,\break 2n-1$. Then{\rm ,} after a possible rearrangement of $(\tau'_1,z'_1; 
\nek ; \tau'_{r'},z'_{r'})${\rm ,} without affecting the
decreasing order of $z'_1\nek z'_{r'}${\rm ,}
\begin{itemize}
\item[{\rm (1)}] $r=r'$ and $m_i=m_i'$ for $i=0,1,\cdots,r${\rm ,} 
\item[{\rm (2)}] $z_i=z'_i$ and $\tau_i\cong\tau_i'$ for
$i=1,2,\cdots,r${\rm ,}
\item[{\rm (3)}] for the representations $\sigma_0$ and
$\sigma'_0$ of $\SO_{2m_0+1}(k)${\rm ,} the twisted gamma
factors are the same{\rm ,} i.e.\ $$\gamma(\sigma_0\times\varrho,s,\psi)=\gamma(\sigma'_0
\times\varrho,s,\psi)$$
for all irreducible supercuspidal representations
$\varrho$ of $\GL_l(k)$ with $l$ as above. 
\end{itemize}

\endproclaim

\demo{Proof}
By assumption,  
\begin{equation}
\gamma(\sigma\times\varrho,s,\psi)=\gamma(\sigma'\times
\varrho,s,\psi)
\end{equation}
for all irreducible supercuspidal representations
$\varrho$ of $\GL_l(k)$ with $l=1,2,\cdots,\break 2n-1$. 
By Proposition 3.2 and Corollary 3.2, we have $z_1=z'_1$. 
Taking $\varrho=\tau_1$ we know that $\gamma(\sigma\times
\tau_1,s,\psi)$ and hence $\gamma(\sigma'\times\tau_1,s,
\psi)$ each has a pole at $s=1+z_1=1+z'_1$. 
By the proof of Proposition 3.2 applied to the local gamma 
factor $\gamma(\sigma'\times\tau_1,s,\psi)$, we find that 
either (1) $z_1=z'_i=0$ and $\tau'_i\cong\tau_1^\vee$, for some $i$,
or (2) $z_1=z'_i$ and $\tau'_i\cong\tau_1$, for some $i$, or (3)
$\tau_1\cong\tau_1^\vee$ and $z_1=0$.  We shall consider
these three cases as follows. 

If case (1) occurs, then  $z'_1=z_1=0$. This
implies that 
$$z_1=\cdots=z_r=0=z'_1=\cdots=z'_{r'}.$$
In this case we change the order of the supercuspidal
support $(Q';\tau'_1,\tau'_2,\cdots,\break\tau'_{r'};
\sigma'_0)$ by making $\tau^{'\vee}_i$ the first
representation, so that now $\tau'_1\cong\tau_1$.
It follows that 
$$\gamma(\tau_1\times\varrho,s,\psi)=\gamma(\tau'_1 
\times\varrho,s,\psi)$$ 
and 
$$\gamma(\tau_1^\vee\times\varrho,s,\psi)=
\gamma({\tau'_1}^\vee\times\varrho,s,\psi)$$ 
for all irreducible supercuspidal representations $\varrho$ of $\GL_l(k)$. 
By Theorem 3.1 (multiplicativity of gamma factors),  
$$
\gamma(\sigma\times\varrho,s,\psi)
=
\omega_{\varrho}(-1)^{rn+m_0}
\left[\prod_{i=1}^r\gamma(\tau_i\times\varrho,s,\psi)
\gamma(\tau^\vee_i\times\varrho,s,\psi)\right]
\gamma(\sigma_0\times\varrho,s,\psi)
$$
and 
$$\gamma(\sigma'\times\varrho,s,\psi)=\omega_{\varrho}
(-1)^{r'n+m'_0}\left[\prod_{i=1}^{r'}\gamma(\tau'_i\times
\varrho,s,\psi)\gamma({\tau'_i}^\vee\times\varrho,s,\psi)\right]
\gamma(\sigma'_0\times\varrho,s,\psi).$$
By cancelling the gamma factors 
$$\gamma(\tau_1\times\varrho,s,\psi)=\gamma(\tau'_1
\times\varrho,s,\psi)$$ 
and 
$$\gamma(\tau_1^\vee\times\varrho,s,\psi)=\gamma
({\tau'_1}^\vee\times\varrho,s,\psi)$$ 
from equation (5.1), we obtain 
a new identity for products of gamma factors, i.e.\ \begin{eqnarray}
& &
\omega_{\varrho}(-1)^{rn+m_0}
\left[\prod_{i=2}^r\gamma(\tau_i\times\varrho,s,\psi)
\gamma(\tau^\vee_i\times\varrho,s,\psi\right]
\gamma(\sigma_0\times\varrho,s,\psi)  
 \\
&= &
\omega_{\varrho}(-1)^{r'n+m'_0}
\left[\prod_{i=2}^{r'}\gamma(\tau'_i\times\varrho,s,\psi)
\gamma({\tau'_i}^\vee\times\varrho,s,\psi\right]
\gamma(\sigma'_0\times\varrho,s,\psi)\nonumber.
\end{eqnarray}  
It is clear that the arguments to prove Proposition 3.2 and Corollary 3.2 are 
applicable to the above identity. 
This reduces the proof to the case where the number of gamma factors in 
both sides of the identity is smaller by two. 

If case (2) occurs, then we have $z_1=z'_1$, $z_1=z'_i$ and
$\tau'_i\cong\tau_1$. 
Hence,
$$
z'_1=\cdots=z'_i.
$$
If we change the order of the supercuspidal data 
$(Q';\tau'_1,\tau'_2,\cdots,\tau'_{r'};\sigma'_0)$ by
making $\tau'_i$ the first representation, so that
$\tau'_1\cong\tau_1$, then the same argument as in case
(1) reduces the proof to (5.2).

Finally, we shall show that case (3) can reduce to cases (1) and (2), or yield 
a contradiction.  If case (3) occurs, then $\tau_1$ is self-dual,
i.e.\ $\tau_1\cong\tau_1^\vee$ and $z_1=0$.  Hence, 
$$z_1=\cdots=z_r=0=z'_1=\cdots=z'_{r'}.$$
By the conclusion of case (1) and case (2), we may assume that the 
representation 
$\tau_1$ is not isomorphic to any one of the representations $\tau'_i$, 
${\tau'_i}^\vee$, for 
$i=1,2,\cdots r'$. This implies that the pole at $s=1$ ($z_1=0$) of 
the gamma factor $\gamma(\sigma'\times\tau_1,s,\psi)$ can only be achieved by 
the gamma factor $\gamma(\sigma'_0\times\tau_1,s,\psi)$. 
Hence this gamma factor 
$\gamma(\sigma'\times\tau_1,s,\psi)$ has a simple pole at $s=1$. 
On the other hand, since $\tau_1\cong\tau_1^\vee$, the gamma factor 
$\gamma(\sigma\times\tau_1,s,\psi)$ can be written as 
$$
\gamma(\sigma\times\tau_1,s,\psi)
=
[\gamma(\tau_1\times\tau_1,s,\psi)]^2\cdot[\cdots].
$$
It has at least a second order pole at $s=1$. This is a contradiction.  

Therefore, by repeating the above argument, we can finally conclude that 
(1) $r=r'$, (2) up to permutation, $z_i=z'_i$ and 
$\tau_i\cong\tau_i'$ for $i=1,2,\cdots,r$, and 
$$
\gamma(\sigma_0\times\varrho,s,\psi)=\gamma(\sigma'_0\times\varrho,s,\psi)
$$
for all irreducible supercuspidal representations $\varrho$ of $\GL_l(k)$ with 
$l=1,\cdots,\break 2n-1$. 
\enddemo

\demo{{\rm 5.1.} The proof of the local converse theorem
{\rm (Theorem 1.2)}} 
Let $\sigma$ and $\sigma'$ be irreducible admissible
generic representations of $\SO_{2n+1}(k)$. Assume that
$\sigma$ and $\sigma'$ have supercuspidal supports
$(Q;\tau_1,\tau_2,\cdots,\tau_r;\sigma_0)$ 
and $(Q';\tau'_1,\tau'_2,\cdots,\tau'_{r'};\sigma'_0)$, and exponents 
$(z_1,z_2,\cdots,z_r)$ and $(z'_1,z'_2,\cdots,z'_{r'})$, respectively. 
More precisely, 
the representation $\sigma$ is a subquotient of the normalized induced 
representation 
$$
\sigma\prec 
\Ind^{\SO_{2n+1}}_{Q}(\tau_1|\det|^{z_1}\otimes\cdots\otimes\tau_r|\det|^{z_r}
\otimes\sigma_0)
$$
and the representation $\sigma'$ is a subquotient of the normalized induced 
representation 
$$
\sigma'\prec \Ind^{\SO_{2n+1}}_{Q'}(\tau'_1|\det|^{z'_1}\otimes\cdots\otimes
\tau'_{r'}|\det|^{z'_{r'}}\otimes\sigma'_0). 
$$
Without loss of generality, we may assume that the exponents 
$(z_1,z_2,\cdots,z_r)$ and $(z'_1,z'_2,\cdots,z'_{r'})$ satisfy the condition:
$$
z_1\geq z_2\geq\cdots\geq z_r\geq 0;\ \ \ z'_1\geq z'_2\geq\cdots\geq z'_{r'}\geq 0.
$$

We assume now that the gamma factors 
attached to $\sigma$ and $\sigma'$, twisted by any irreducible supercuspidal 
representation $\varrho$ of $\GL_l(k)$, are the same, i.e.\ $$
\gamma(\sigma\times\varrho,s,\psi)=\gamma(\sigma'\times\varrho,s,\psi)
$$
with $l=1,2,\cdots,2n-1$. By Theorem 5.1, we conclude that 
(1) $r=r'$ and $m_i=m_i'$ for $i=0,1,\cdots,r$; 
(2) $z_i=z'_i$ and $\tau_i\cong\tau_i'$ for $i=1,2,\cdots,r$ 
(after a possible reordering of the $(\tau'_i,z'_i)$,
which does not affect the decreasing order of the
exponents; and (3) as representations of $\SO_{2m_0+1}(k)$,
the twisted gamma factors attached to $\sigma_0$ and
$\sigma'_0$ are the same, i.e.\ $$\gamma(\sigma_0\times\varrho,s,\psi)=\gamma
(\sigma'_0\times\varrho,s,\psi)$$
for all irreducible supercuspidal representations $\varrho$ of $\GL_l(k)$ with 
$l=1,2,\cdots,\break 2n-1$. Since $n\geq m_0$, it follows from Theorem 4.1 that the 
representations $\sigma_0$ and $\sigma'_0$ are isomorphic. Hence 
both irreducible admissible generic representations $\sigma$ and $\sigma'$ 
have the 
same supercuspidal support $(Q;\tau_1,\tau_2,\cdots,\tau_r;\sigma_0)$ and the 
same exponent $(z_1,z_2,\cdots,z_r)$. In other words, both $\sig$ and $\sig'$ 
are irreducible generic constituents (up to equivalence) of the induced 
representation 
$$
\Ind^{\SO_{2n+1}}_{Q}(\tau_1|\det|^{z_1}\otimes\cdots\otimes\tau_r|\det|^{z_r}
\otimes\sigma_0). 
$$
By the uniqueness of the generic constituent in the induced representation 
$$
\Ind^{\SO_{2n+1}}_{Q}(\tau_1|\det|^{z_1}\otimes\cdots\otimes\tau_r|\det|^{z_r}
\otimes\sigma_0), 
$$
the irreducible admissible generic representations $\sigma$ and $\sigma'$ 
must be equivalent. This proves the local converse theorem (Theorem 1.2) 
in general. 
\hfill \qed\enddemo

5.2. {\it Some direct applications}.
We shall briefly discuss two global applications to the theory of automorphic 
representations of $\SO_{2n+1}({\Bbb A})$. First we obtain the injectivity of 
the weak Langlands functorial lifting established in \cite{CKP-SS}. 

\specialnumber{5.2}\proclaim{Theorem}
Let ${\cal SO}_{2n+1}^{\rm igca}({\Bbb A})$ be the set of all equivalence classes 
of irreducible 
generic cuspidal automorphic representations of $\SO_{2n+1}({\Bbb A})$ and 
${\cal GL}_{2n}^{\rm ia}({\Bbb A})$ be the set of all equivalence classes of 
irreducible 
automorphic representations of $\GL_{2n}({\Bbb A})$. Then the weak 
Langlands functorial lifting from ${\cal SO}_{2n+1}^{\rm igca}({\Bbb A})$ to 
${\cal GL}_{2n}^{\rm ia}({\Bbb A})${\rm ,} established in {\rm \cite{CKP-SS},}
is an injective map. 
\endproclaim

\demo{Proof}
The proof is now a straightforward consequence of
Proposition 3.3 and our local converse theorem
(Theorem 1.2). (Recall again that the weak lift of
\cite{CKP-SS} is functorial at archimedean places.) 
\enddemo 

\proclaimtitle{rigidity theorem}
\specialnumber{5.3}\proclaim{Theorem} 
Let $\Sigma=\otimes_v\Sigma_v$ and $\Sigma'=\otimes_v \Sigma'_v$ belong to 
$\SO^{\rm igca}_{2n+1}({\Bbb A})$ {\rm (}\/as defined in Theorem {\rm 5.2).} 
If $\Sigma_v$ is equivalent to $\Sigma'_v$ 
for almost all local places $v${\rm ,} then $\Sigma$ is equivalent to $\Sigma'$. 
\endproclaim

\demo{Proof}
Let ${\cal T}$ and ${\cal T}'$ be the weak Langlands functorial 
liftings of $\Sigma$ and $\Sigma'$, respectively, as constructed in 
\cite{CKP-SS}. 
This means that at all archimedean local places and unramified finite local 
places, the lifting is the local Langlands functorial lifting. 
By assumption we know that both representations ${\cal T}$ and 
${\cal T}'$ are equivalent at almost all local places. By Theorem 3.2, both 
${\cal T}$ and ${\cal T}'$ are irreducible generic self-dual unitary 
automorphic representations of $\GL_{2n}({\Bbb A})$. 
By the strong multiplicity-one property for 
$\GL(n)$ (\cite{JS}), we conclude that 
$$
{\cal T}\cong{\cal T}'
$$
as automorphic representations of $\GL_{2n}({\Bbb A})$.
By Theorem 5.2, $\Sig\cong\Sig'$.
\enddemo

{\it Remark} 5.1.
 (1)  The rigidity theorem for $\GL_n$ was proved by  
H. Jacquet and J. Shalika (\cite{JS}) for 
generic automorphic representations of $\GL_n(\AA)$. By applying the 
multiplicity-one  theorem for irreducible cuspidal automorphic representations 
of $\GL_n(\AA)$ (\cite{Shl}), we see that the rigidity theorem implies the strong 
multiplicity-one theorem for $\GL(n)$. The rigidity theorem for some lower rank
cases has also been studied in \cite{B2}, \cite{R}, and \cite{S4}. 
\vglue4pt
 (2)  In \cite{GRS5}, the local converse theorem proved in this
paper for generic representations of $\SO_{2n+1}(k)$ is
used to prove the irreducibility of the backward lift
to $\SO_{2n+1}({\Bbb A})$ of a representation of
$\GL_{2n}({\Bbb A})$ of the form described in Theorem 3.2.
This was conjectured in \cite{GRS1}.

\section{On local Langlands conjectures}

In this chapter, we establish the local Langlands
functorial lift from irreducible generic supercuspidal
representations of $\SO_{2n+1}(k)$ to $\GL_{2n}(k)$ 
and the local Langlands reciprocity law for irreducible generic supercuspidal 
representations of $\SO_{2n+1}(k)$. We shall consider the local Langlands 
conjectures for more general representations of $\SO_{2n+1}(k)$ and other 
relevant problems in a forthcoming work of ours \cite{JngS}. 

\demo{{\rm 6.1.} On local Langlands functorial lifting from $\SO(2n+1)$ to $\GL(2n)$}
Let ${\cal SO}_{2n+1}^{\rm igsc}(k)$ be the set of all equivalence classes of 
irreducible generic supercuspidal representations of $\SO_{2n+1}(k)$ and 
${\cal GL}_{2n}^{\rm ifl}(k)$ 
(`ifl' denotes the image of the functorial lifting) be the set of all 
equivalence classes of irreducible admissible generic representations of 
$\GL_{2n}(k)$ of the form 
$$
\tau=\eta_1\times\eta_2\times\cdots\times\eta_t, 
$$
where $\eta_i$ are irreducible unitary supercuspidal self-dual 
representations of\break  $\GL_{2n_j}(k)$ ($j=1,2,\cdots,t$) 
($\sum_{j=1}^tn_i=n$) such that 
\begin{itemize}
\item[(1)] $\eta_i\not\cong\eta_j$ if $i\neq j$, and 
\item[(2)] the local $L$-function $L(\eta_j,\Lambda^2,s)$ has a pole at $s=0$ 
for $j=1,2,\cdots,t$. 
\end{itemize}
For any $\sigma\in {\cal SO}_{2n+1}^{\rm igsc}(k)$, the
weak Langlands functorial lifting from $\SO_{2n+1}$
to $\GL_{2n}$ produces a map 
$$\ell\ :\ \sigma\mapsto \tau=\ell(\sigma)$$
such that 
$$\gamma(\sigma\times\varrho,s,\psi)=\gamma(\tau\times
\varrho,s,\psi)$$
for all irreducible supercuspidal representations $\varrho$ of $\GL_l(k)$ where 
$l$ is any positive integer (Proposition 3.3). 
Then by Theorem 3.3, we know that the image $\tau=\ell(\sigma)$ belongs to 
${\cal GL}_{2n}^{\rm ifl}(k)$. So there is a map 
\begin{equation}
\ell\ :\ \sigma\mapsto \tau=\ell(\sigma)
\end{equation}
from ${\cal SO}_{2n+1}^{\rm igsc}(k)$ to ${\cal GL}_{2n}^{\rm ifl}(k)$, which 
preserves the twisted local gamma factors; i.e.,
 \begin{equation}
\gamma(\sigma\times\varrho,s,\psi)
=
\gamma(\tau\times\varrho,s,\psi)
\end{equation}
for all irreducible supercuspidal representations
$\varrho$ of $\GL_l(k)$, where $l$ is any positive integer.
In fact, we can prove that the map $\ell$ is a bijection. 
\enddemo

\specialnumber{6.1}\proclaim{Proposition}
Let $\ell$ be the map defined above with   property  {\rm (6.2).}
Then the map $\ell$ 
is a bijection from ${\cal SO}_{2n+1}^{\rm igsc}(k)$ onto 
${\cal GL}_{2n}^{\rm ifl}(k)$. Moreover{\rm ,} such a map is unique. 
\endproclaim

\demo{Proof}
The injectivity of the map $\ell$ follows from the local converse theorem for 
$\SO_{2n+1}$ (Theorem 1.2). The uniqueness of such a map $\ell$ follows from 
the local converse theorem for $\GL_n$ (Theorem 1.1). 
It remains to show that the map $\ell$ is surjective. 

For any $\tau'$ in ${\cal GL}_{2n}^{\rm ifl}(k)$, and 
$\tau'=\eta'_1\times\eta'_2\times\cdots\times\eta'_{t'}$, we have 
to construct a  representation $\sigma$ 
in ${\cal SO}_{2n+1}^{\rm igsc}(k)$ such that 
$$
\gamma(\sigma\times\varrho,s,\psi)
=
\gamma(\tau'\times\varrho,s,\psi)
$$
for all irreducible supercuspidal representations $\varrho$ of $\GL_l(k)$ where 
$l$ is any positive integer. 
By the local backward lifting from $\GL_{2n}(k)$ to 
$\tilspn(k)$ (\cite{GRS6}), 
there exists a unique (up to isomorphism) irreducible 
$\psi_{{\tilde{U}},1}$-generic supercuspidal representation 
$\pi_\psi(\tau')$ of $\tilspn(k)$ such that the twisted local gamma 
factor $\gamma(\pi_\psi(\tau')\times\tau',s,\psi)$ has a pole of order 
$t'$ at $s=1$. By Part (2) of Theorem 2.2, there exists a unique nontrivial 
irreducible generic supercuspidal representation $\sigma$ of $\SO_{2n+1}(k)$, 
which is the local $\psi$-Howe lift of $\pi_\psi(\tau')$. 
By Theorem 3.3, there exists an irreducible 
admissible generic representation $\tau$ of $\GL_{2n}(k)$ such that 
$$
\gamma(\tau\times\varrho,s,\psi)
=
\gamma(\sigma\times\varrho,s,\psi)
$$
for all irreducible supercuspidal representations $\varrho$ of $\GL_l(k)$ where 
$l$ is any positive integer, and the representation $\tau$ has the 
following properties: 
\begin{itemize}
\item[(1)] $\tau=\eta_1\times\cdots\times\eta_t$ with 
$\eta_i\neq \eta_j$ 
if $i\neq j$, where $\eta_i$ is an irreducible unitary self-dual 
supercuspidal representation of $\GL_{2n_i}(k)$;
\item[(2)] the local $L$-function $L(\eta_i,\Lambda^2,s)$ has a pole at $s=0$ 
for $i=1,2,\cdots,t$.
\end{itemize}
In particular, $\gamma(\sigma\times\tau,s,\psi)$ has a pole of order $t$
at $s=1$. The proof of Proposition 3.4, with $\pi=\pi_{\psi}(\tau')$ can
be repeated with simple modifications to conclude (from the fact that
$\gamma(\pi\times\tau',s,\psi)$ has a pole of order $t'$ at $s=1$) that
$\gamma(\sigma\times\tau',s,\psi)$ has a pole of order $t'$ at $s=1$. All 
we need to do is to take ${\cal T}_{\nu_0}\cong\tau'$, and then in the paragraph
after (3.13), conclude first, using Lemma 3.2 that the right-hand side of (3.13) has
a pole of order $t'$ at $s=1$, and hence conclude that the left-hand side of (3.13)
has a pole of order $t'$ at $s=1$. Now, by Lemma 3.1, we see that
$\gamma(\sigma\times\tau',s,\psi)$ has a pole of order $t'$ at $s=1$. Since
$$
\gamma(\sigma\times\tau',s,\psi)=\gamma(\tau\times\tau',s,\psi)
$$
we conclude, looking at the form of $\tau$ and $\tau'$ (both are 
in ${\cal GL}^{\rm ifl}_{2n}(k)$), that $\tau\cong\tau'$. 
Therefore the representation $\sigma$ just constructed has the 
property that 
$$
\ell(\sigma)=\tau
$$
and we have proved the map $\ell$ is surjective. 
\enddemo

\specialnumber{6.2}\proclaim{Proposition}
The map 
$$
\ell\ :\ \sigma\mapsto \tau=\ell(\sigma)
$$
from ${\cal SO}_{2n+1}^{\rm igsc}(k)$ to ${\cal GL}_{2n}^{\rm ifl}(k)$ 
preserves the twisted local $\epsilon$\/{\rm -}\/factors and local\break $L$-factors\/{\rm ;} i.e.{\rm ,}
 $$
\epsilon(\sigma\times\varrho,s,\psi)
=
\epsilon(\tau\times\varrho,s,\psi)
$$ 
and 
$$
L(\sigma\times\varrho,s)
=
L(\tau\times\varrho,s)
$$ 
for all irreducible supercuspidal representations $\varrho$ of $\GL_l(k)$ where 
$l$ is any positive integer.
\endproclaim

\demo{Proof}
Let $\tau=\ell(\sigma)$. Then we know that 
$$
\tau=\eta_1\times\eta_2\times\cdots\times\eta_t, 
$$
where $\eta_i$ are irreducible unitary self-dual supercuspidal 
representations of\break $\GL_{2n_j}(k)$ ($j=1,2,\cdots,t$) ($\sum_{j=1}^tn_i=n$) 
such that 
\begin{itemize}
\item[(1)] $\eta_i\not\cong\eta_j$ if $i\neq j$, and 
\item[(2)] the local $L$-function $L(\eta_j,\Lambda^2,s)$ has a pole at $s=0$ 
for $j=1,2,\cdots,t$. 
\end{itemize}
More importantly, we have 
$$
\gamma(\sigma\times\varrho,s,\psi)
=
\gamma(\tau\times\varrho,s,\psi)
$$
for all irreducible supercuspidal representations $\varrho$ of $\GL_l(k)$ 
where $l$ is any positive integer (Theorem 3.3 and Proposition 6.1). It follows 
from \cite{Sh1}, \cite{Sh2} and \cite{JP-SS} that 
\begin{equation}
\gamma(\sigma\times\varrho,s,\psi)
=
\epsilon(\sigma\times\varrho,s,\psi)\cdot
{L(\sigma\times\varrho^\vee,1-s)\over L(\sigma\times\varrho,s)}
\end{equation}
and
\begin{equation}
\gamma(\tau\times\varrho,s,\psi)
=
\epsilon(\tau\times\varrho,s,\psi)\cdot
{L(\tau\times\varrho^\vee,1-s)\over L(\tau\times\varrho,s)}. 
\end{equation}
By assumption, we have 
\begin{equation}
\epsilon(\sigma\times\varrho,s,\psi)\cdot
{L(\sigma\times\varrho^\vee,1-s)\over L(\sigma\times\varrho,s)}
=
\epsilon(\tau\times\varrho,s,\psi)\cdot
{L(\tau\times\varrho^\vee,1-s)\over L(\tau\times\varrho,s)}. \hskip.25in
\end{equation}
If the supercuspidal representation $\varrho$ is not equivalent to any one 
of the $\eta_i$'s, up to unramified unitary twisting, then  
$$
L(\tau\times\varrho^\vee,1-s)=1=L(\tau\times\varrho,s). 
$$ 
Hence, by equation (6.5), we have, 
$$
{L(\sigma\times\varrho^\vee,1-s)\over L(\sigma\times\varrho,s)}
=
\epsilon(\tau\times\varrho,s,\psi)\cdot
\epsilon(\sigma\times\varrho,s,\psi)^{-1},
$$
which is an exponential function in $s$. 

We first assume that $\varrho$ is unitary. Since both $\sigma$ and $\varrho$ 
are supercuspidal, by Proposition 7.2 in \cite{Sh2}, 
the possible poles of $L(\sigma\times\varrho^\vee,1-s)$ lie in ${\rm Re}(s)\geq 1$, 
but the possible poles of $L(\sigma\times\varrho,s)$ lie in ${\rm Re}(s)\leq 0$. 
Hence, 
$$
L(\sigma\times\varrho^\vee,1-s)=1=L(\sigma\times\varrho,s). 
$$
Therefore,  \def\nsp{\hskip-8pt}
\begin{eqnarray*}
L(\sigma\times\varrho,s)
&\nsp=\nsp&
L(\tau\times\varrho,s) \\[-5pt]
\noalign{\noindent and} 
\epsilon(\sigma\times\varrho,s,\psi)
&\nsp=\nsp&
\epsilon(\tau\times\varrho,s,\psi).
\end{eqnarray*}
It is clear that the same argument works when $\varrho$ is not necessarily \pagebreak unitary. 

If the supercuspidal representation $\varrho$ is isomorphic to one of the 
$\eta_i$'s, up to unramified unitary twisting, say 
$$
\varrho\cong\eta_i\cdot |\cdot|^y 
$$
where $y$ is purely imaginary, then we know again that there are no 
cancellations between the poles of $L(\sigma\times\varrho^\vee,1-s)$ and 
the poles of $L(\sigma\times\varrho,s)$, and 
the same thing happens with $L(\tau\times\varrho^\vee,1-s)$ and 
$L(\tau\times\varrho,s)$. 
Hence, from equation (6.5), the set of the poles of 
$L(\sigma\times\varrho,s)$ 
is equal to the set of poles of $L(\tau\times\varrho,s)$.  
Thus, the polynomials $L(\sigma\times\varrho,s)^{-1}$ and 
$L(\tau\times\varrho,s)^{-1}$ are equal. Therefore  
$$
L(\sigma\times\varrho,s)
=
L(\tau\times\varrho,s).
$$
It follows that the $\epsilon$-factors are also equal, i.e.\ 
\vglue12pt \hfill $\epsilon(\sigma\times\varrho,s,\psi)
=
\epsilon(\tau\times\varrho,s,\psi).
$ 
\enddemo

\vglue9pt
The following theorem on   local Langlands functoriality follows from 
Propositions 6.1 and 6.2. 

\proclaimtitle{local Langlands functoriality}
\proclaim{Theorem} 
There exists a unique bijective map 
$$
\ell\ :\ \sigma\mapsto \tau=\ell(\sigma)
$$
from ${\cal SO}_{2n+1}^{\rm igsc}(k)$ to ${\cal GL}_{2n}^{\rm ifl}(k)${\rm ,} which 
preserves the twisted local $L$\/{\rm -}\/factors{\rm ,}\break $\epsilon$\/{\rm -}\/factors and gamma factors{\rm ,} 
i.e.\ 
\begin{eqnarray*}
\noalign{\vskip-4pt}
L(\sigma\times\varrho,s)
&\nsp=\nsp&
L(\tau\times\varrho,s),
\\ 
\epsilon(\sigma\times\varrho,s,\psi)
&\nsp=\nsp&
\epsilon(\tau\times\varrho,s,\psi)
\\ 
\noalign{\noindent 
and}
\gamma(\sigma\times\varrho,s,\psi)
&\nsp=\nsp&
\gamma(\tau\times\varrho,s,\psi)
\end{eqnarray*}
for all irreducible supercuspidal representations $\varrho$ of $\GL_l(k)$ where 
$l$ is any positive integer.
\endproclaim  

6.2. {\it On the local Langlands reciprocity law for} $\SO(2n+1)$.
We shall establish the local Langlands reciprocity law for 
${\cal SO}_{2n+1}^{\rm igsc}(k)$ 
by using the local Langlands reciprocity law for $\GL(n)$ established by 
M. Harris and R. Taylor \cite{HT} and by G. Henniart \cite{Hn3}. 

Let $W_k$ be the Weil group associated to the the local field $k$. We take  
$$
W_k\times \SL_2({\Bbb C})
$$
as the Weil-Deligne group (\cite{M} and \cite{Kn}). 
Let ${\cal G}_n^{\rm ah}(k)$ be the set of conjugacy classes of admissible 
homomorphisms $\rho$ from $W_k\times{\rm SL}_2({\Bbb C})$ to $\GL_n({\Bbb C})$. 
If we write 
$$
\rho=\mathbold{\oplus}_i\rho^0_i\otimes\lam^0_i,
$$
then the admissibility of $\rho$ means that $\rho^0_i$'s are continuous complex 
representations of $W_k$ with $\rho^0_i(W_k)$ semi-simple and $\lam^0_i$'s 
are algebraic complex representations of $\SL_2(\CC)$. Let ${\cal GL}_n^{\rm is}(k)$ 
be the set of equivalence classes of irreducible smooth 
representations of $\GL_n(k)$. Then the local Langlands conjecture (or local 
Langlands reciprocity law), now a theorem of Harris-Taylor \cite{HT} and 
Henniart \cite{Hn3}, is the following. 

\proclaimtitle{Harris-Taylor \cite{HT} and Henniart \cite{Hn3}}
\proclaim{Theorem} 
There exists a (unique) bijection 
$$
{\frak r}_n\ :\ \rho\mapsto \tau={\frak r}_n(\rho) 
$$
from ${\cal G}_n^{\rm ah}(k)$ onto ${\cal \GL}_n^{\rm is}(k)$ satisfying the following 
conditions.
\begin{itemize}
\item[{\rm (1)}] For any $\rho\in {\cal G}_n^{\rm ah}(k)${\rm ,} $\det(\rho)$ corresponds to 
$\omega_{{\frak r}_n(\rho)}${\rm ,} the central character\/{\rm ;}\/
\item[{\rm (2)}] For any $\rho\in {\cal G}_n^{\rm ah}(k)$ and any quasi\/{\rm -}\/character 
$\chi$ of $k^\times${\rm ,} one has 
${\frak r}_n(\chi\otimes\rho)=(\chi\circ\det)\otimes{\frak r}_n(\rho)${\rm ;} 
\item[{\rm (3)}] For any $\rho\in {\cal G}_n^{\rm ah}(k)${\rm ,} one has 
${\frak r}_n(\rho)^\vee={\frak r}_n(\rho^\vee)${\rm ;}
\item[{\rm (4)}] For any $\rho\in{\cal G}_n^{\rm ah}(k)$ and 
$\rho'\in{\cal G}_{n'}^{\rm ah}(k)${\rm ,} one has 
\begin{itemize}
\item[$(4L)$] $L(\rho\otimes\rho',s)=
L({\frak r}_n(\rho)\times{\frak r}_{n'}(\rho'),s)${\rm ,}
\item[(4$\epsilon$)]
$\epsilon(\rho\otimes\rho',s,\psi)=
\epsilon({\frak r}_n(\rho)\times{\frak r}_{n'}(\rho'),s,\psi)${\rm ,}
\item[(4$\gamma$)] 
$\gamma(\rho\otimes\rho',s,\psi)=
\gamma({\frak r}_n(\rho)\times{\frak r}_{n'}(\rho'),s,\psi)${\rm ;}
\end{itemize}
\item[{\rm (5)}] If  $\rho=(\rho^0,\delta)${\rm ,} then $(\rho^0,1)$ with 
$\rho^0$ irreducible corresponds to ${\frak r}_n(\rho)${\rm ,} which is 
irreducible and supercuspidal. 
\end{itemize}
\endproclaim

{\it Remark} 6.1.
 (1)\ The uniqueness of the reciprocity map in Theorem 6.2 follows from Henniart's 
local converse theorem (Theorem 1.1) and an induction argument on $n$. 

\vglue4pt
(2)\ Theorem 6.2 has been proved for supercuspidal representations by Harris and 
Taylor (\cite{HT}) and by Henniart (\cite{Hn3}). The reduction of the general case 
to the supercuspidal case was given by A. Zelevinsky (\cite{Z}). Various special 
cases of Theorem 6.2 were proved before by several authors. 
We refer to \cite{H} and \cite{K2} for detailed comments. 
\vglue8pt

In order to establish the local Langlands reciprocity conjecture for\break 
${\cal SO}_{2n+1}^{\rm igsc}(k)$, the set of equivalence classes of irreducible 
generic supercuspidal representations of $\SO_{2n+1}(k)$, it is sufficient to 
figure out the subset of the local Langlands parameters for 
${\cal GL}_{2n}^{\rm ifl}(k)$ by using the explicit local Langlands functorial lift 
from $\SO_{2n+1}(k)$ to $\GL_{2n}(k)$ (Theorem 6.1) and the local 
Langlands reciprocity law for $\GL(n)$ (Theorem 6.2). 

Recall that the set 
${\cal GL}_{2n}^{\rm ifl}(k)$ consists of equivalence classes of representations of 
$\GL_{2n}(k)$ of the form: 
$$
\tau=\eta_1\times\eta_2\times\cdots\times\eta_t, 
$$
where $\eta_i$ are irreducible unitary supercuspidal self-dual 
representations of\break $\GL_{2n_j}(k)$ ($j=1,2,\cdots,t$) 
($\sum_{j=1}^tn_i=n$) such that 
\begin{itemize}
\item[(1)] $\eta_i\not\cong\eta_j$ if $i\neq j$, and 
\item[(2)] the local $L$-function $L(\eta_j,\Lambda^2,s)$ has a pole at $s=0$ 
for $j=1,2,\cdots,t$. 
\end{itemize}
Each irreducible unitary supercuspidal self-dual representation $\eta_j$ of 
$\GL_{2n_j}(k)$ has the local Langlands parameter $\rho^0_j$, which 
is an irreducible, $2n_j$-dimensional, admissible, complex representation of $W_k$, 
by Theorem 6.2. Further, we have $\rho^0_i\not\cong\rho^0_j$ if $i\neq j$. Hence 
the representation 
$$
\tau=\eta_1\times\eta_2\times\cdots\times\eta_t
$$
has the local Langlands parameter 
$$
\rho^0=\rho^0_1\oplus\rho^0_2\oplus\cdots\oplus\rho^0_t,
$$
which is a $2n$-dimensional, admissible, completely reducible, multiplicity-free,  
complex representation of $W_k$. 

Recently, G. Henniart communicated to us (\cite{Hn1}) 
that he can prove the following results 
among others satisfied by the local Langlands reciprocity map. 

 \proclaimtitle{Henniart \cite{Hn1}}
\proclaim{Theorem}
The local Langlands reciprocity map has the following property\/{\rm :}
the gamma factor $\gamma(\rho,\Lambda^2,s,\psi)$ has the same poles as 
the local gamma factor $\gamma({\frak r}_n(\rho),\Lambda^2,s,\psi)$ for 
any irreducible $\rho$ {\rm (}\/i.e.\ for any ${\frak r}_n(\rho)$
supercuspidal\/{\rm ).}\/ 
\endproclaim

By using Theorem 6.3, we can prove the following proposition. 

\specialnumber{6.3}\proclaim{Proposition}
$(1)$\ Let $\rho^0$ be an irreducible{\rm ,} $2m$\/{\rm -}\/dimensional{\rm ,} admissible{\rm ,} complex 
representation of $W_k$ and $\tau$ be an irreducible unitary supercuspidal 
representation of $\GL_{2m}(k)$ with the properties that $({\rm i})$\ $\tau$ has the 
local Langlands parameter $\rho^0$ and $({\rm ii})$\ the local exterior square $L$\/{\rm -}\/function 
$L(\tau,\Lambda^2,s)$ has a pole at $s=0$. Then $\rho^0$ is symplectic{\rm ,} i.e.\ $$
\rho^0(W_k)\subset\Sp_{2m}(\CC).
$$

$(2)$\ Let $\rho^0=\rho^0_1\oplus\rho^0_2$ be a $2m$\/{\rm -}\/dimensional{\rm ,} admissible{\rm ,} 
completely reducible{\rm ,} complex representation of $W_k$ with the property that 
$$
\Hom_{W_k}(\rho^0_1\otimes\rho^0_2,1)=0.
$$
Then $\rho^0$ is symplectic if and only if $\rho^0_1$ and $\rho^0_2$ are both 
symplectic.
\endproclaim

\demo{Proof}
The proof of Part (1) follows from Theorem 6.3. More precisely it goes as follows. 
Since 
$$
\gamma(\tau^\vee,\Lambda^2,s,\psi)=\epsilon(\tau^\vee,\Lambda^2,s,\psi)\cdot
{L(\tau,\Lambda^2, 1-s)\over L(\tau^\vee,\Lambda^2, s)}, 
$$
and by the assumption, the local $L$-function $L(\tau,\Lambda^2, s)$ has a pole at 
$s=0$, we obtain that the gamma factor $\gamma(\tau^\vee,\Lambda^2,s,\psi)$ has 
a pole at $s=1$. By Theorem 6.3, the gamma factor 
$\gamma((\rho^0)^\vee,\Lambda^2,s,\psi)$ has a pole at $s=1$. Because we also have 
$$
\gamma((\rho^0)^\vee,\Lambda^2,s,\psi)
=\epsilon((\rho^0)^\vee,\Lambda^2,s,\psi)\cdot
{L(\rho^0,\Lambda^2, 1-s)\over L((\rho^0)^\vee,\Lambda^2,s)}, 
$$
we get that the $L$-function $L(\rho^0,\Lambda^2,s)$ has a pole at $s=0$. 
Now it is an elementary fact that if $L(\rho^0,\Lambda^2,s)$ has a pole at $s=0$, 
then the image $\rho^0(W_k)$ is included in $\Sp_{2m}({\Bbb C})$, i.e.\ the 
parameter $\rho^0$ is symplectic.  
This proves Part~(1). 

Part (2) is basically proved by linear algebra. It is clear that if both 
$\rho^0_1$ and $\rho^0_2$ are symplectic, then $\rho^0$ is itself symplectic. 
Conversely, we use a basic fact from linear algebra that $\rho^0$ is symplectic 
if and only if $\Lam^2(\rho^0)$ has $W_k$-invariant functionals 
([GW, \S 5.1.7]). Since 
$$
\Lam^2(\rho^0)=\Lam^2(\rho^0_1)\oplus\Lam^2(\rho^0_2)\oplus
[\rho^0_1\otimes\rho^0_2],
$$
the $W_k$-invariant functionals will not vanish on at least one of 
$\Lam^2(\rho^0_1)$, $\Lam^2(\rho^0_2)$, since we assume that 
$$
\Hom_{W_k}(\rho^0_1\otimes\rho^0_2,1)=0.
$$
Without loss of generality, we assume that $\Lam^2(\rho^0_1)$ supports a  
$W_k$-invariant functional. Hence $\rho^0_1$ is symplectic. Because 
$\rho^0$ is nondegenerate and $\rho^0_2$ is the complement of $\rho^0_1$, we 
conclude that $\rho^0_2$ is also symplectic. 
\enddemo

Let ${\cal G}^0_{2n}(k)$ be the set of conjugacy classes of all $2n$-dimensional, 
admissible, completely reducible, multiplicity-free, symplectic 
complex representations $\rho^0$ of $W_k$. Then we have the following theorem. 

 \proclaimtitle{local Langlands reciprocity law}
\proclaim{Theorem}
There exists a unique bijection 
$$
{\frak R}_{2n}\ :\ \rho^0_{2n}\mapsto {\frak R}_{2n}(\rho_{2n}^0)
$$
from the set ${\cal G}_{2n}^0(k)$ onto the set 
${\cal SO}_{2n+1}^{\rm igsc}(k)$ such that 
\vglue6pt
 $(L)$\ $L(\rho^0_{2n}\otimes\rho^0_l,s,)
=
L({\frak R}_{2n}(\rho_{2n})\times {\frak r}_l(\rho^0_l),s)${\rm ,} 
\vglue4pt 
$\, (\epsilon)$\ $\epsilon(\rho^0_{2n}\otimes\rho^0_l,s,\psi)
=
\epsilon({\frak R}_{2n}(\rho_{2n})\times {\frak r}_l(\rho^0_l),s,\psi)${\rm ,} and 
\vglue4pt $\, (\gamma)$\ $\gamma(\rho^0_{2n}\otimes\rho^0_l,s,\psi)
=
\gamma({\frak R}_{2n}(\rho_{2n})\times {\frak r}_l(\rho^0_l),s,\psi)$ 
\vglue6pt\noindent 
for all irreducible continuous representations of $W_k$ of dimension $l$. 
\endproclaim

\demo{Proof}
The theorem is a direct consequence of Theorem 6.1 and Proposition 6.3.
\enddemo

{\it Remark} 6.2.
The complete local Langlands reciprocity conjecture in this case 
states that the reciprocity 
map ${\frak R}_{2n}$ takes a local Langlands parameter $\rho^0$ in 
${\cal G}_{2n}^0(k)$ to a finite set $\Pi(\rho^0)$ (local $L$-packet) of irreducible 
admissible representations of $\SO_{2n+1}(k)$. By our local converse theorem 
(Theorem 1.2), we know that in each local $L$-packet $\Pi(\rho^0)$, there is 
at most one generic member (i.e.\ with Whittaker model). It is a very interesting 
and difficult problem to give an explicit construction of the local $L$-packets.
\enddemo

\section{Appendix:  On gamma factors for
$\tilsp_{2n}\times \GL_l$}

7.1. {\it Review of the global theory}
In \cite{GRS3}, $L$-functions for generic, automorphic,
cuspidal representations on $\tilsp^{\phantom{1}}_{2n}\times \GL_l$,
are represented via global integrals of Shimura type.
We review this construction briefly.  It yields local
gamma factors at each place.

Let $F$ be a number field.  Denote by $\AA$  its adele
ring.  Fix a nontrivial character $\psi_0$ of
$F\bks\AA$.  Let $\Pi$  (resp.\ $\rho$)  be an
irreducible, automorphic, cuspidal representation of
$\tilsp_{2n}^{\phantom{l}}(\AA)$  (resp.\ $\GL_l (\AA)$).  Assume
that $\Pi$  is globally
$(\psi_0)_{\tilU_{\AA:,1}}$-generic.  In the sequel, the
cuspidality of $\rho$ is not important. What we need is
that $\rho$ is automorphic, realized in an irreducible
subspace of automorphic forms on $\GL_l(\AA)$,  and
that $\rho$ is globally generic (i.e.\ the Whittaker
coefficient is nontrivial on the space of $\rho$). 
Although this is not pointed out in \cite{GRS3}, the
proofs there do not use at all the cuspidality of
$\rho$.  Thus, we may take $\rho$  to be an Eisenstein
series induced from irreducible, automorphic, cuspidal
representations at a point of holomorphy.

Let $\xi_{\rho,s}$ be a holomorphic section for
$J_{\rho,s}$ -- the representation of
$\Sp_{2l}(\AA)$, induced from $\rho_s=\rho|\det\cdot
|^{s-1/2}$  on the Siegel parabolic subgroup $Q_l
(\AA)$, and denote by $E(\xi_{\rho,s,\cdot})$  the
corresponding Eisenstein series on $\Sp_{2l}(\AA)$.
We distinguish two cases according to whether $n\ge l$
or $n<l$.  In the first case, a Fourier-Jacobi
coefficient of a cusp form in $\Pi$  is paired against
the Eisenstein series above, and in the second case, a
cusp form in $\Pi$ is paired against a Fourier-Jacobi
coefficient of the Eisenstein series.  As we need in
this paper  only the case $n<l$   (as a matter of fact,
we need just the case $l =2n$,) we assume from now on
that $n<l$.

Let $w_{\psi_0^{-1}}$  be the Weil representation of
$\tilsp_{2n}(\AA)$, corresponding to $\psi_0^{-1}$.
Realize it in ${\cal S}(\AA^n)$, and denote, for $\phi\in
{\cal S}(\AA^n)$, by $\tet^\phi_{\psi_0^{-1}}$ the
corresponding theta series. Extend $w_{\psi_0^{-1}}$ and
$\tet^\phi_{\psi_0^{-1}}$  to $\calH_n(\AA)$ -- the
Heisenberg group in $2n+1$  variables.  Let
$N_{l,n+1}$  be the unipotent radical of the standard
parabolic subgroup $Q_{l,n+1}$  of $\Sp_{2l}$,
whose Levi part is isomorphic to
$\GL_1^{l-n-1}\times\Sp_{2(n+1)}$. Let $(\psi_0)_{N_{l ,n+1}}$ be 
the restriction to $N_{l,n+1}(\AA)$ of the standard
nondegenerate character defined by $\psi_0$.  Note that
$\calH_n$ embeds naturally in $N_{l,n}$ so that
$N_{l,n} =N_{l,n+1}\rtimes\calH_n$.
Extend $(\psi_0)_{N_{l,n+1}}$  to $N_{l,n}(\AA)$, so
that it is trivial on $\calH_n(\AA)$. Denote this
extension by $\chi_{\psi_{0;l,n}}$.  Denote by $j$
the projection of $N_{l,n}$  to $\calH_n$.  Let $j$
denote also the embedding of $\Sp_{2n}$ into the Levi 
part of $Q_{l,n}$.
Note that $j$  embeds $\Sp_{2n}\rtimes\calH_n$  into the
Levi part of $Q_{l,n+1}$.

Let $f$ be an automorphic function on $\Sp_{2l}(\AA)$.
A Fourier-Jacobi coefficient of $f$  of type
$(\psi_0,n)$, is a function on $\tilsp_{2n}(\AA)$ of the
form
$$f^\phi_{\psi_{0,n}}(\tilg)=\int_{N_{l,n}(F)\bks 
N_{l,n}(\AA)}f(uj(g))\tet^\phi_{\psi_0^{-1}}(j(u)\tilg)
\chi_{\psi^{-1}_{0,l,n}}(u)du\ .$$
Here $\tilg\in\tilsp_{2n}(\AA)$ projects to
$g\in\Sp_{2n}(\AA)$, and $\phi\in {\cal S}(\AA^n)$. Let
$\varphi$  be a cusp form in (the space of) $\Pi$.
Define 
$$\calL(\varphi,\phi,\xi_{\rho,s})=\int_{\Sp_{2n}(F)\bks
\Sp_{2n}(\AA)}\
\varphi(g)E^\phi_{\psi_{0,n}}(\xi_{\rho,s},g)dg\ .$$
(In case $n>l$, one takes
$\varphi^\phi_{\psi_{0,l}}$ and pairs it with
$E(\xi_{\tau,s},\cdot)$  along
$\Sp_{2l}(F)\bks\hskip-1pt\Sp_{2l}(\AA)$,  and in case
$n\hskip-1pt =\hskip-1pt l$, one integrates
$\varphi(g)\tet^\phi_{\psi_0^{-1}}(g)E(\xi_{\rho,s},g)$
along $\Sp_{2n}(F)\bks\hskip-1pt\Sp_{2n}(\AA)$.) We have an Euler
product decomposition, for decomposable data and\break ${\rm Re}
(s)\gg 0$
$$\calL(\varphi,\phi,\xi_{\rho,s})=\prod_\nu\calB
(W_{\varphi,\nu},\phi_\nu,\xi_{\rho_{\nu,s}})$$
where $W=\prod_\nu W_{\varphi,\nu}$  is the Whittaker
function of $\varphi$  with respect to $\psi_0$, and at
each place $\nu$,
\begin{eqnarray}\label{eq7.1}
&&\calB(W_{\varphi,\nu},\phi_\nu,\xi_{\rho_{\nu,s}})
 \\[0.3em]
&&\qquad =\int_{U_n(F_\nu)\bks\Sp_{2n}(F_\nu)}\
\int_{N_{l,n}^{\gam_{l,n}}(F_\nu)\bks
N_{l,n}(F_\nu)}\nonumber\\[0.5em]
&&\qquad\quad W_{\varphi,\nu}(g)w_{\psi^{-1}_0}(j(u)g)
\phi_\nu(e_0) f_{\xi_{\rho_{\nu,s}}}(\gam_{l,n}uj(g))
\chi_{\psi^{-1}_{0,\nu;l,n}}(u)dudg.\nonumber
\end{eqnarray}
Here $U_n$  is the standard maximal unipotent subgroup
of $\Sp_{2n}$, $\gam_{l,n}$ is a certain Weyl
element,
$N^{\gam_{l,n}}_{l,n}=\gam^{-1}_{l,n}Q_l
\gam_{l,n}\cap N_{l,n}$; and $e_0=(0\nek 0,1)$.  We obtain
$f_{\xi_{\rho,s}}=\Pi f_{\xi_{\rho_{\nu,s}}}$  from $\xi_{\rho,s}$ after taking a certain
Whittaker coefficient in the\break ``$\rho$--variable".  Thus,
we consider sections $f_{\xi_{\rho_{\nu,s}}}$ for
$J_{\rho_{\nu,s}}$, which take values in a 
certain Whittaker model $\rho_\nu$.

For decomposable data $\varphi,\phi,\xi_{\rho,s}$,  let
$S$  be a finite set of places, including those at
infinity, those above 2, and such that outside $S$  all
data are unramified, and $\psi_0$  is normalized.  Then
(normalizing $W_{\varphi,\nu}(I)=1$ outside $S$), we have  
\begin{equation}\label{eq7.2}
\calL(\varphi,\phi,\xi_{\rho,s})=\prod_{\nu\in
S}\calB(W_{\varphi,\nu},\phi_\nu,\xi_{\rho_{\nu,s}})
\frac{L^S_{\psi_0}(\Pi\times\rho,s)}{L^S(\rho,s+\frac{1}
{2})L^S(\rho,\Lam^2,2s)}.\hskip.4in
\end{equation}
This implies that $L^S_{\psi_0}(\Pi\times\rho,s)$ is
meromorphic.  Indeed $\calL(\varphi,\phi,\xi_{\rho,s})$
is clearly meromorphic, and $L^S(\rho,s +\frac{1}{2}),
L^S(\rho,\Lam^2,2s)$  are known to be meromorphic.   
For finite $\nu$ in $S$, we can choose data, such that
$\calB(W_{\varphi,\nu},\phi_\nu,\xi_{\rho_{\nu,s}})=1$,
for all~$s$ (see \cite[Prop.~6.6]{GRS3}) and given
$s_0\in\CC$  and $\nu\in S$  which is archimedean, we
can find a combination
$\sum^N_{i=1}\calB\Big(W^{(i)}_\nu,\phi^{(i)}_\nu,
\xi^{(i)}_{\rho_{\nu,s}}\Big)$ which is holomorphic and
nonzero at $s_0$  (see \cite[Prop.~6.7]{GRS3}).  From
this we conclude, choosing data in the same way, at all
places of $S$, but one place $\nu_0$, that
$\calB(W_{\varphi,\nu_0}, \phi_{\nu_0},
\xi_{\rho_{\nu_0,s}})$ is meromorphic.
(This can be shown in general, without the assumption that
the data are coming from global cusp forms.  See 
\cite[\S 1.1]{GRS2} for the case $\nu <\infty$, where
it follows that
$\calB(W_\nu,\phi_\nu,\xi_{\rho_{\nu,s}})$ is rational
in $q^{-s}_\nu$.  The case where $\nu$ is infinite can
be done exactly as in \cite{S2}.)
\def\eqref#1{(\ref{#1})}

Applying in \eqref{eq7.1} the intertwining operator 
$M$, with respect to $\left( \begin{array}{cc}
&I_l\\
-I_l
\end{array}\right)$ on $J_{\rho,s}$,  we get
\begin{eqnarray}\label{eq7.3}
&&\calL(\varphi,\phi,M(\xi_{\rho,s}))\\[0.5em]
&&\qquad =\prod_{\nu\in
S}\calB(W_{\varphi,\nu},\phi_\nu,M_\nu
(\xi_{\rho_{\nu,s}}))\frac{L^S(\rho,s-\frac{1}{2})L^S
(\rho,\Lam^2,2s-1)}{L^S(\rho,s+\frac{1}{2})L^S
(\rho,\Lam^2,2s)}\nonumber\\[0.5em]
&&\qquad \quad\cdot\ \frac{L^S_{\psi_0}(\Pi\times\hatrho,1-s)}
{L^S(\hatrho,\frac{3}{2}-s)L^S(\hatrho,\Lam^2,2-2s)}.\nonumber
\end{eqnarray}
Using the functional equation satisfied by
$E(\xi_{\rho,s,\cdot})$, we can equate \eqref{eq7.2} and
\eqref{eq7.3} to get 
\begin{eqnarray}\label{eq7.4}
&&  \frac{L^S_{\psi_0}(\Pi\times\rho,s)L^S(\hatrho,\frac{3}
{2}-s)L^S(\hatrho,\Lam^2,2-2s)}{L^S_{\psi_0}
(\Pi\times\hatrho,1-s)L^S(\rho, s-\frac{1}{2})
L^S(\rho,\Lam^2, 2s-1)}\prod_{\nu\in S}\calB
(W_{\varphi,\nu},\phi_\nu,\xi_{\rho_{\nu,s}})
\\[0.5em] 
&&\hskip1.5in =\prod_{\nu\in
S}\calB\Big(W_{\varphi,\nu},\phi_\nu, M_\nu
(\xi_{\rho_{\nu,s}})\Big).\nonumber
\end{eqnarray}
Fixing data at all places in $S$  except one place
$\nu_0$, we conclude from \eqref{eq7.4} that there is a
meromorphic function $\Gam
(\Pi_{\nu_0}\times\rho_{\nu_0},s,\psi_{0,\nu_0})$,
which is rational in $q^{-s}_{\nu_0}$, in case $\nu_0$
is finite, such that 
$$\Gam(\Pi_{\nu_0}\times\rho_{\nu_0},s,\psi_{0,\nu_0})
\calB(W_{\varphi,\nu_0},\phi_{\nu_0},\xi_{\rho_{\nu_0,s}})
=\calB(W_{\varphi,\nu_0},\phi_{\nu_0},M_{\nu_0}
(\xi_{\rho_{\nu_0,s}}))$$
for all
$W_{\varphi,\nu_0},\phi_{\nu_0},\xi_{\rho_{\nu_0,s}}$.  
We define the local gamma factor
$\gam(\Pi_{\nu_0}\times\rho_{\nu_0},s,\psi_{0,\nu_0})$
by
$$\Gam(\Pi_{\nu_0}\times\rho_{\nu_0},s,\psi_{0,\nu_0})
=\frac{\gam(\Pi_{\nu_0}\times\rho_{\nu_0},s,
\psi_{0,\nu_0})}{\gam(\rho_{\nu_0},s-\frac{1}{2},
\psi_{0,\nu_0})\gam(\rho_{\nu_0},\Lam^2,2s-1,
\psi_{0,\nu_0})}.$$
Thus
\begin{equation}\label{eq7.5}
\gam(\Pi_{\nu_0}\times\rho_{\nu_0},s,\psi_{0,\nu_0})
\calB(W_{\varphi,\nu_0},\phi_{\nu_0},\xi_{\rho_{\nu_0,s}}
)=\tilcalb (W_{\varphi,\nu_0},\phi_{\nu_0},
\xi_{\rho_{\nu_0,s}})\hskip.4in
\end{equation}
where $\tilcalb(W_{\varphi,\nu_0},\phi_{\nu_0},
\xi_{\rho_{\nu_0,s}})=\calB(W_{\varphi,\nu_0},
\phi_{\nu_0},M^*_{\nu_0}(\xi_{\rho_{\nu_0,s}})$, and 
$$M^*_{\nu_0}(\xi_{\rho_{\nu_0,s}})=\gam\left(\rho_{\nu_0},s-
\frac{1}{2},\psi_{0,\nu_0}\right)\gam(\rho_{\nu_0},\Lam^2,2s-1,
\psi_{0,\nu_0})M_{\nu_0}(\xi_{\rho_{\nu_0,s}})\ .$$
Note that at a finite place $\nu_0$  where
$\psi_{0,\nu_0}$ is normalized; and
$\Pi_{\nu_0}$ and $\rho_{\nu_0}$  are unramified, we have 
$$
\gam(\Pi_{\nu_0}\times\rho_{\nu_0},s,\psi_{0,\nu_0})
=
\frac{L_{\psi_{0,\nu_0}}(\Pi_{\nu_0}\times\hatrho_{\nu_0},1-s)}
{L_{\psi_{0,\nu_0}}(\Pi_{\nu_0}\times\rho_{\nu_0},s)}.
$$
For such $\nu_0, L_{\psi_0,\nu_0}$
$(\Pi_{\nu_0}\times\rho_{\nu_0},s)$
is nothing but $L(\tet_{\psi_0,\nu_0}(\Pi_{\nu_0})
\times\rho_{\nu_0},s)$, where
$\tet_{\psi_0,\nu_0}(\Pi_{\nu_0})$ is the unramified
representation of $\SO_{2n+1}(F_{\nu_0})$  corresponding
to $\Pi_{\nu_0}$  by the local $\psi_{0,\nu_0}$--Howe
lift.

\demo{{\rm 7.2.} A result on gamma factors at archimedean places} 
Let $\nu_0$  be an archimedean place of $F$.  Put
$\Pi_{\nu_0}=\pi,\rho_{\nu_0}=\tau, F_{\nu_0}=k$,
$W_{\varphi,\nu_0}=W$, $\phi_{\nu_0}=\phi$,
$\psi_{0,\nu_0}=\psi$. 
In this subsection, we denote $V_\pi$ and $V_\tau$ to be 
the canonical models of the Harish-Chandra modules of $\pi$ and $\tau$, respectively. 
The canonical extension of a Harish-Chandra module has the $C^\infty$-topology as 
given in [C]. From [C], such canonical extensions are unique, up to equivalence. 
Our goal in this subsection is to 
show that there is $A>0$, such that for $|{\rm Im}(s)|>A$,
$\gam (\pi\times \tau,s,\psi)$  is holomorphic and
nonzero.  The global integral $\calL(\varphi,\phi,\xi_{\rho,s})$ defined in 
subsection 7.1 is separately continuous in the $C^\infty$-topology at the place $\nu_0$. 
Hence it remains continuous after the extension to the canonical models. 
Therefore, we may regard (the analytic continuation of)
$\calB (W,\phi,\xi_{\tau,s})$  as a continuous linear
form $T$ on $V_\pi\otimes {\cal S}(k^n)\otimes V_{J_{\tau,s}}$. 
Here the notions of separate continuity and of continuity
coincide, and the two notions of tensor products
(inductive $\overline{\otimes}$, projective
$\widehat{\otimes}$) coincide. The proof of this is as
follows.  

First we note that $V_\pi$ and
$V_{J_{\tau,s}}$  are nuclear Fr{\'e}chet spaces.
Indeed, both representations are quotients of
(differentiably) induced representations coming off
Borel subgroups and quasi-characters, and since the
spaces of such induced representations are images of
surjective maps from spaces of the form $C^\infty (K)$,
where $K$  is compact, our spaces $V_\pi$  and
$V_{J_{\tau,s}}$  are quotients of such spaces.  Since
$C^\infty(K)$ is Fr{\'e}chet and nuclear, so are $V_\pi$
and $V_{J_{\tau,s}}$.  (See [Tr, pp.\ 85, 94, 514, 530].)
In particular, the two notions of tensor product for
$V_\pi\otimes V_{J_{\tau,s}}$ coincide; i.e., $V_\pi\overline{\otimes}V_{J_{\tau,s}}\cong
V_\pi\widehat{\otimes} V_{J_{\tau,s}}$,  which (by  
\cite[p.~514]{Tr}) is nuclear.  We conclude that the
two notions of tensor product for $(V_\pi\otimes
V_{J_{\tau,s}})\otimes {\cal S}(k^n)$  coincide.  Actually,
${\cal S}(k^n)$ is nuclear as well \cite[p.~530]{Tr}.

We add some more related remarks.  Note that the same
proofs work for $V_\pi\otimes V_\tau$  and
$(V_\pi\otimes V_\tau)\otimes {\cal S}(k^n)$  as well (i.e.\ $\otimes$  can be replaced by either
$\overline{\otimes}$  or $\widehat{\otimes}$).  Note
also that $V_\pi\otimes V_\tau$ is a Fr{\'e}chet space
as well, since it is a quotient of $C^\infty (K_1)\otimes
C^\infty (K_2)\cong C^\infty (K_1\times K_2)$, where
$K_1,K_2$  are compact, such that $V_\pi$ is a quotient
of $C^\infty (K_1)$ and $V_\tau$ is a quotient of
$C^\infty (K_2)$.  We conclude that $V_\pi\otimes
V_\tau\otimes {\cal S}(k^n)$  is a quotient of a Fr{\'e}chet
space, and hence is itself a Fr{\'e}chet space.
Indeed, $V_\pi\otimes V_\tau\otimes {\cal S}(k^n)$ is a
quotient of $C^\infty (K_1\times K_2)\otimes {\cal S}(k^n)$ 
which is isomorphic to ${\cal S}(k^n;C^\infty(K_1\times K_2))$
\cite[p.~533]{Tr}, and the last space is isomorphic to
${\cal S}(k^n\times K_1\times K_2)$, which is  a Fr{\'e}chet space
\cite[p.~92]{Tr}.  Note also that $V_\pi\otimes
V_\tau\otimes {\cal S}(k^n)$ is nuclear as a tensor product of
two nuclear spaces $V_\pi\otimes V_\tau$ and ${\cal S}(k^n)$.
As a corollary, we obtain (see \cite[p.~485]{Wr}):
\enddemo

\specialnumber{7.1}\proclaim{Proposition}\label{prop7.1}
Let $M$  be a $C^\infty$\/{\rm -}\/manifold{\rm ,} countable at
infinity. Then
\begin{equation}\label{eq7.6}
C^\infty_c(M)\otimes (V_\pi\otimes V_\tau\otimes
{\cal S}(k^n))\cong C^\infty_c(M; V_\pi\otimes V_\tau\otimes
{\cal S}(k^n)).
\end{equation}
\endproclaim
Again, in Proposition \ref{prop7.1}, each
$\otimes$ can be replaced by either $\overline{\otimes}$
or $\widehat{\otimes}$. 

The linear form $T$  (on $V_\pi\otimes {\cal S}(k^n)\otimes
V_{J_{\tau,s}}$) is equivariant with respect to the
subgroup $R=j(\Sp_{2n}(k))N_{l,n}(k)$.
\begin{equation}\label{eq7.7}
T((\pi(\tilg)\otimes w_{\psi^{-1}}(j(u)\tilg)\otimes
J_{\tau,s}(uj(g)))v)=\chi_{\psi;l,n}(u)T(v)
\end{equation}
where $\tilg$ is (any) inverse image in $\tilspn(k)$ of
$g$  in $\Sp_{2n}(k)$, and $u\in N_{l,n}(k)$.  This
follows easily from \eqref{eq7.1} (or even from the
structure of the global integrals).  We have a surjection
\begin{eqnarray*}
&&\tau^1_s:C^\infty_c(\Sp_{2l}(k);
V_\tau)\longrightarrow V_{J_{\tau,s}}\\
&&\tau^1_s(\varphi)(h)=\int_{Q_l(k)}
\del^{-1/2}_{Q_l}(p)\tau_s(p^{-1})(\varphi(ph))d_rp\ ,
\end{eqnarray*}
where $d_rp$  is a right invariant measure on
$Q_l(k)$. Composing $T$ with $\tau^1_s$  yields a
(continuous) linear form $t$ on 
\begin{eqnarray*}
 C^\infty_c(\Sp_{2l}(k); V_\tau)\otimes
{\cal S}(k^n)\otimes V_\pi&\cong&
(C^\infty_c(\Sp_{2l}(k))\otimes V_\tau\otimes
{\cal S}(k^n)\otimes V_\pi\\
& \cong& C^\infty_c(\Sp_{2l}(k); V_\tau\otimes
{\cal S}(k^n)\otimes V_\pi)
\end{eqnarray*}
(We used Proposition \ref{prop7.1}.)  By  
\eqref{eq7.7} and
$$\tau^1_s(\lam(p)\varphi)=\del^{1/2}_{Q_l}(p)
\tau^1_s(\tau_s(p^{-1}\circ\varphi))\ ,\quad p\in Q_l
(k)$$
($\lam(p)$ denotes the left translation by $p^{-1}$), we
conclude that $t$,  when regarded as a $V_\tau\otimes
{\cal S}(k^n)\otimes V_\pi$-distribution on
$\Sp_{2l}(k)$, satisfies 
\begin{eqnarray}\label{eq7.8}
\quad t(r(uj(g))f)&=&\chi_{\psi;l,n}(u)t((1\otimes
w_{\psi^{-1}}(\tilg^{-1}j(u^{-1}))\otimes\pi(\tilg^{-1}))
\circ f),\\ 
 t(\lam(p)f)&=&\del^{1/2}_{Q_l}(p)t((\tau_s(p^{-1})
\otimes 1\otimes 1)\circ f)\ . \nonumber
\end{eqnarray}
Here $u\in N_{l,n}(k)$, $g\in\Sp_{2n}(k)$, $p\in
Q_l (k)$  and $r$  denotes right translations.  We
are now at the situation of \cite[5.2.4]{Wr}.  (Note
that in the notation \cite[5.2.4]{Wr},
$M=\Sp_{2l}(k)$ and $G=Q_l (k)\times R$  acts on
$M$  by $(p,r)\cdot h=phr^{-1}$.  See also 
\cite[p.~408]{Wr}.)  We have a nice description of the
set $Q_l (k)\bks\Sp_{2l}(k)/R$.  This is a finite
set, and it has one open orbit $Q_l
(k)\gam_{l,n}R$.  See \cite[Sec.~4]{GRS3}.  The
reference for the following is \cite[5.2.3,5.2.4]{Wr}.

Let us show that the map
$$b:t\mapsto t\Big|_{C^\infty_c(Q_l(k)\gam_{l,n}R;
V_\tau\otimes {\cal S}(k^n)\otimes V_\pi)}$$
is injective on the space of $V_\tau\otimes
{\cal S}(k^n)\otimes V_\pi$-distributions on
$\Sp_{2l}(k)$, which satisfy \eqref{eq7.8}.  Indeed,
if $b(t)=0$, then by Bruhat theory (see the proof of 
\cite[Lemma 5.2.4.4]{Wr}), $t$ is supported in the
complement of the open orbit $Q_l (k)\gam_{l,n}R$.
The dimension of the space of such distributions is
majorized by
\begin{eqnarray*} \sum_{\gam_{l,n}\not=\gam\in Q_l
(k)\bks\Sp_{2l}(k)/R}\ \sum^\infty_{m=1}\dim\Big[&&\hskip-18pt
\Bil_{Q_l(k)\cap\gam R\gam^{-1}}\\
&& \hskip-18pt\cdot \Big(
\del^{-1/2}_{Q_l}\tau_s\otimes(\chi^{-1}_{\psi;
l,n}\otimes\pi\otimes
w_{\psi^{-1}})^\gam,\Lam_{\gam,m}\Big)\Big]
\end{eqnarray*}
where $\Lam_{\gam,m}$  are certain algebraic finite-dimensional representations, coming from derivatives.
$\Bil_H$  denotes $H$-equivariant bilinear forms.  Let
$V_l$  denote the unipotent radical of $Q_l$.  An
element of $$\Bil_{Q_l(k)\cap\gam R\gam^{-1}}\Big(
\del^{-1/2}_{Q_l}\tau_s\otimes (\chi^{-1}_{\psi;
l,n}\otimes\pi\otimes
w_{\psi^{-1}})^\gam,\Lam_{\gam,m}\Big)$$ when regarded as
a $V_l(k)\cap\gam N_{l,n}(k)\gam^{-1}$-equivariant
form embeds
$$\Big(\chi^{-1}_{\psi;l,n}\cdot\psi^{-1}\Big|_{{\rm Center}
(\calH_n(k))}\Big)^\gam$$  on $E=V_l (k)\cap\gam
N_{l,n+1}(k)\cdot {\rm Center} (\calH_n(k))\gam^{-1}$ into
the dual of $\Lam_{\gam,m}$.  Since $\Lam_{\gam,m}$  is
an algebraic finite-dimensional representation, it
cannot have nontrivial eigenvalues on the unipotent
subgroup $E$.  Thus, we must have
$$\Big(\chi_{\psi;l,n}\cdot \psi\Big|_{{\rm Center} (\calH_n
(k))}\Big)^\gam=1.$$  By \cite[p.~212]{GRS3}, this is
impossible, unless $Q_l (F)\gam R$  is the open
orbit.  This is a contradiction, and so the map $b$ is
injective.

Returning to
$T(W\otimes\phi\otimes\xi_{\tau,s})=
\calB(W,\phi,\xi_{\tau,s})$, let $s_0$  be a pole of
order $e$  of $\calB(W,\phi,\xi_{\tau,s})$  in the sense
that $(s-s_0)^e\calB(W,\phi,\xi_{\tau,s})$ is
holomorphic and not identically zero at $s_0$.  Let $t$
be the linear form on 
$C^\infty_c(\Sp_{2l}(k);V_\tau)\otimes {\cal S}(k^n)\otimes
V_\pi$ defined by
\begin{equation}\label{eq7.9}
t(\varphi\otimes\phi\otimes W)=\lim_{s\to
s_0}(s-s_0)^e\calB(W,\phi,\tau^1_s(\varphi))\ .
\end{equation}
When viewed as a $V_\tau\otimes {\cal S}(k^n)\otimes V_\pi$
distribution on $\Sp_{2l}(F)$, $t$ clearly satisfies  
\eqref{eq7.8}.  Then what we have just shown is that for
$\varphi$  supported in the open orbit $Q_l
(k)\gam_{l,n}R$, $t(\varphi\otimes\phi\otimes W)$ is
not identically zero.  This means that all poles of
$\calB(W,\phi,\xi_{\tau,s})$  are detected (with their
orders) on the open orbit.  Thus, in order to locate the
poles, it is enough to take $\xi_{\tau,s}$  with compact
support modulo $Q_l(k)$  (independent of $s$) inside
$Q_l (k)\gam_{l,n}R$.  We may take
$\xi_{\tau,s}=\tau^1_s(\varphi)$, with $\varphi$
supported in $Q_l(k)\gam_{l,n}R$.  For such
$\xi_{\tau,s}$, the unipotent inner integration in
\eqref{eq7.1} converges absolutely.  Rewrite
\eqref{eq7.1}, for ${\rm Re}(s)\gg 0$, following the Iwasawa
decomposition $\Sp_{2n}(k)=U_n(k)AK$
\begin{eqnarray}\label{eq7.10}
&&\\[-4pt]
\calB(W,\phi,\xi_{\tau,s})
&=&
\int_{h\in K}\int_{V_0}\int_A
\del^{-1}(\hata)W(\hata\tilh)w_{\psi^{-1}}(\hata
j(u)\tilh)\phi(e_0)\nonumber\\
& &\cdot \xi_{\tau,s}(\gam_{l,n}uj(h);
\left( \begin{array}{cc}
a&\\
&I_{l -n}
\end{array}\right))|\det a|^{s+s_{l, n}}\chi_{\psi;
l,n}(u)dadudh.\nonumber
\end{eqnarray}
Here $V_0$ is a compact subset of
$N^{\gam_{l,n}}_{l,n}(k)\bks N_{l,n}(k)$ (the
projection onto\break $N^{\gam_{l,n}}_{l,n}(k)\bks
N_{l,n}(k)$ of the compact support modulo $Q_l (k)$
of $\xi_{\tau,s}$  inside $Q_l(k)\gam_{l,n}R$),
$\hata$  denotes an inverse image in $\tilspn(k)$ of
$\left( \begin{array}{cc}
a&\\
&a^*
\end{array}\right)$, for a diagonal matrix $a$ in $\GL_n(k)$, 
and $s_{l,n}$  is a certain fixed translation of
$s$.  Denote the inner $da$ integration, for fixed
$(h,u)$, in \eqref{eq7.10} by $B(\pi(\tilh)W,
w_{\psi^{-1}}(j(u)\tilh)\phi,
J_{\tau,s}(uj(h))\xi_{\tau,s})$.  For a fixed $(h,u)\in
K\times V_0$, $W'_\pi =\pi(\tilh)W, \phi'=w_{\psi^{-1}}
(j(u)\tilh)\phi$, $\xi'_{\tau,s}=J_{\tau,s}
(uj(h))\xi_{\tau,s}$ and $W'_\tau(m)=\xi'_{\tau,s}(I;m)$,
the above inner integration equals
\begin{equation}\label{eq7.11}
\int_A\del^{-1}(\hata)W'_\pi
(\hata)w_{\psi^{-1}}(\hata)\phi(e_0)W'_\tau\left( \begin{array}{cc}
a&\\
&I_{l -n}
\end{array}\right)
|\det a|^{s+s_{l,n}}da\ . \hskip.4in
\end{equation}
The analytic continuation of \eqref{eq7.11} is obtained
when we replace $W'_\pi(\hata)$,\break $W'_\tau\left( \begin{array}{cc}
a&\\
&I_{l -n}
\end{array}\right)$ with their asymptotic expansions,
obtained exactly as in \cite[\S 3.3]{S1}.  (See also
\cite[\S 4]{S2} which applies here in exactly the same
way.) We conclude that \eqref{eq7.11} is a sum of
elements of the form
\begin{equation}\label{eq7.12}
\int f(a_1\nek a_n)\eta (a_1\nek
a_n)\del^{-1}(\hata)|\det a|^{s+s'_{l,n}}da
\end{equation}
where $a=\diag (a_1\nek a_n)$, $f\in {\cal S}(k^n)$  ($f$  is
independent of $s$) and $\eta$  is a finite function
varying in a finite set which depends on $\pi$  and
$\tau$  only ($s'_{l,n}$  is another fixed
translation of $s$).  Thus, there is a finite set of 
characters $X_{\pi,\tau}$ of $k^*$, and there is a
polynomial $P_{\pi,\tau}(s)$, which depend on $\pi$  and
$\tau$  only, such that $\frac{\calB
(W'_\pi,\phi',\xi'_{\tau,s})}{P_{\pi,\tau}(s)\prod_{\mu\in
X_{\pi,\tau}}L(\mu,s)}$ is holomorphic in the whole
plane. We have  
\begin{equation}\label{eq7.13}
\frac{\calB(W,\phi,\xi_{\tau,s})}{P_{\pi,\tau}(s)
\prod_{\mu\in X_{\pi,\tau}}L(\mu,s)} =\int_{K_\times
V_0}\frac{\calB(\pi(\tilh)W,w_{\psi^{-1}}(j(u)\tilh)\phi,
J_{\tau,s}(uj(h))\xi_{\tau,s})}{P_{\pi,\tau}(s)
\prod_{\mu\in X_{\pi,\tau}}L(\mu,s)}dudh.
\end{equation}
The right-hand side of \eqref{eq7.13} is holomorphic
since $K\times V_0$  is compact, the integrand is
continuous in $(u,h)$  and the convergence of the
integral is uniform in $s$,  as $s$  varies in compact
sets.  Looking at the left-hand side of \eqref{eq7.13},
we conclude:

\specialnumber{7.2}\proclaim{Proposition}\label{prop7.2}
There is $A>0${\rm ,} such that $\calB(W,\phi,\xi_{\tau,s})$
is holomorphic for $|\Im (s)|>A${\rm ,} for all data.
\endproclaim

{\it Remark} 7.1.
Since $\calB (W,\phi,M^*(\xi_{\tau,s}))$  has a similar
structure, we may take $A$ in the last proposition so
that $\widetilde{\calB}(W,\phi,\xi_{\tau,s})$  is
holomorphic for $|\Im(s)|>A$, for all data, as well.
Finally, we conclude:

\specialnumber{7.3}\proclaim{Proposition}\label{prop7.3} \hskip-8pt
There is $A>0${\rm ,} such that for $|\Im(s)|>A${\rm ,}
$\gam(\pi\times\tau,s,\psi)$ is holomorphic with no
zeroes.
\endproclaim

\demo{Proof} We have
$$\gam(\pi\times\tau,s,\psi)\calB(W,\phi,
\xi_{\tau,s})=\widetilde{\calB}(W,\phi,\xi_{\tau,s}).$$

Let $A$  be as in Proposition \ref{prop7.2} and in the remark which
follows.  Given any $s_0\in\CC$, there is a combination
$\sum^N_{i=1}\calB(W_i,\phi_i,\xi^{(i)}_{\tau,s})$
which is holomorphic and nonzero at $s_0$.  (See 
\cite[Prop.~6.7]{GRS3}.  Thus, if $s_0$ is a pole of
$\gam(\pi\times\tau,s,\psi)$, then $s_0$ is a pole of
$\widetilde{\calB}(W,\phi,\xi_{\tau,s})$.  This forces $|\Im(s_0)|\le
A$.  Similarly, since
$$\gam(\pi\times\tau,s,\psi)^{-1}\widetilde{\calB}
(W,\phi, \xi_{\tau,s})=\calB(W,\phi,\xi_{\tau,s}),$$
if $s_0$  is a zero of $\gam(\pi\times\tau,s,\psi)$, we
may assume that $M^*$  is an isomorphism between
$J_{\tau,s_0}$  and $J_{\hattau,1-s_0}$ (take $|\Im(s_0)|$
large enough), and then as before, we conclude, that
$s_0$  is a pole of ${\calB}(W,\phi,\xi_{\tau,s})$,  which
is impossible, if we assume that $|\Im(s_0)|>A$. 
\enddemo
\bye

\end{document}

%% file: Jiang.ref
\AuthorRefNames [CKP-SS]

 

%% file: Jiang.bbl
\begin{references}

\bibitem{A}
 \name{J.\ Arthur},  
The problem of classifying automorphic representations of classical
groups, {\it CRM Proc.\ and Lecture Notes\/} {\bf 11} (1997), 1--12.
 


\bibitem{AT}
\name{E.\ Artin} and \name{J.\ Tate}, 
{\it Class Field Theory}, Addison-Wesley, Reading, Mass.\ (1968).
 

\bibitem{B1}
\name{E.\ M.\ Baruch}, 
Local factors attached to representations of $p$-adic groups and strong 
multiplicity one, Ph.D.\ Thesis, Yale Univ., 1995.
  

\bibitem{B2}
\bibline,
On the gamma factors attached to representations of ${\rm U}(2,1)$
over a 
$p$-adic field, {\it Israel J.\ Math\/}.\ {\bf 102} (1997), 317--345.
 

\bibitem{B}
\name{A.\ Borel}, 
Automorphic $L$-functions, in  
{\it Automorphic Forms\/}, {\it Representations and $L$-functions\/} (Corvallis,
Ore., 1977),   Part 2, 27--61, {\it Proc.\ Sympos.\ Pure Math\/}.\ {\bf
33}, Amer.\ Math.\ Soc., Providence, RI, 1979.
 

\bibitem{BZ}
  \name{I.\ N.\ Bernstein} and \name{A.\ V.\  Zelevinsky},
Induced representations of reductive $p$-adic groups.\ I, 
{\it Ann.\ Sci.\ {\rm \'{\it E}}cole Norm.\ Sup\/}.\ {\bf 10} (1977), 441--472.
 


\bibitem{C}
 \name{W.\ Casselman},  
Canonical extensions of Harish-Chandra modules to representations of $G$,
{\it Canad.\ J.\ Math\/}.\ {\bf 41} (1989), 385--438.
  

\bibitem{CS}
 \name{W.\ Casselman} and \name{F.\  Shahidi}, 
On irreducibility of standard modules for generic representations, 
{\it Ann.\ Sci.\ {\rm \'{\it E}}cole Norm.\ Sup\/}.\ {\bf 31} (1998), 561--589.
 


\bibitem{Ch}
  \name{J.\ Chen}, 
Local factors, central characters and representations of ${\rm GL}(n)$
over  a non-archimedean local fields, Ph.D.\ Thesis, Yale University,
1996.
 

\bibitem{CKP-SS}
\name{J.\ W.\ Cogdell, H.\ H.\  Kim,  I.\ I.  Piatetski-Shapiro}, and \name{F.\
Shahidi},  
On lifting from classical groups to ${\rm GL}(n)$, {\it IHES Publ.\
Math\/}.\ {\bf 93} (2001), 5--30.
 


\bibitem{CP-S1}
\name{J.\ W.\ Cogdell} and \name{I.\ I.\  Piatetski-Shapiro},  
Converse theorems for ${\rm GL}_n$.\ II, {\it J.\ Reine Angew.\
Math\/}.\ {\bf 507} (1999), 165--188.
 


\bibitem{CP-S2}
 \bibline,  
Converse theorems for ${\rm GL}_n$, {\it IHES Publ.\ Math\/}.\ {\bf
79} (1994), 157--214.
 


\bibitem{F}
 \name{M.\ Furusawa},
On the theta lift from ${\rm SO}(2n+1)$ to $\widetilde{{\rm Sp}}(n)$,
{\it J.\ Reine Angew.\ Math\/}.\ {\bf 466} (1995), 87--110. 
 

\bibitem{GP-SR}
  \name{D.\ Ginzburg, I.\ I.\ Piatetski-Shapiro}, and \name{S.\ Rallis},
{\it $L$\/{\rm -}\/Functions for the Orthogonal Group\/}, {\it Mem.\ Amer.\ Math.\
Soc\/}.\ {\bf 128}, Providence, RI, 1977.
 

\bibitem{GRS1}
 \name{D.\ Ginzburg, S.\ Rallis}, and \name{D.\ Soudry}, 
On explicit lifts of cusp forms from ${\rm GL}(m)$ to classical 
groups, {\it Ann.\ of Math.\/} {\bf 150} (1999), 807--866.
 

\bibitem{GRS2}
 \bibline,
On a correspondence between cuspidal representations of 
${\rm GL}_{2n}$ and $\widetilde{{\rm Sp}}_{2n}$, 
{\it J.\ Amer.\ Math.\ Soc\/}.\ {\bf 12} (1999), 849--907.
 

\bibitem{GRS3}
 \bibline,
$L$-functions for symplectic groups, 
{\it Bull.\ Soc.\ Math.\ France\/} {\bf 126} (1998), 181--244.
 

\bibitem{GRS4}
\bibline,
Periods, poles of $L$-functions and symplectic-orthogonal theta lifts, 
{\it J.\ Reine Angew.\ Math\/}.\ {\bf 487} (1997), 85--114.
 

\bibitem{GRS5}
\bibline,
 Generic automorphic forms on ${\rm SO}(2n+1)$:
functorial lift to ${\rm GL}(2n)$, endoscopy and base change,
{\it Internat.\ Math.\ Res.\ Notices\/} {\bf 14} (2001), 729--764.
 

\bibitem{GRS6}
\bibline,
Endoscopic representations of $\tilde{{\rm Sp}}_{2n}$,
{\it J.\  of the Institut of  Mathematics of Jussieu\/} {\bf 1}(1) (2002), 77--123. 
 


\bibitem{GW}
 \name{R.\ Goodman} and \name{N.\ Wallach},
Representations and invariants of the classical groups, {\it Encycl.\
Math.\ and its Applications\/} {\bf 68}, Cambridge Univ.\ Press,
Cambridge, 1998.
 

\bibitem{H}
 \name{M.\ Harris},
The local Langlands conjecture for ${\rm GL}(n)$ over a $p$-adic field,
$n<p$,  {\it Invent.\ Math\/}.\ {\bf 134} (1998), 177--210.
  


\bibitem{HT}
 \name{M.\ Harris} and \name{R.\ Taylor},
{\it On the Geometry and Cohomology of Some Simple Shimura Varieties\/},
{\it Ann.\ of Math.\ Studies\/} {\bf 151}, Princeton Univ.\ Press,
Princeton, NJ (2001).
 
 

\bibitem{Hn1}
\name{G.\ Henniart}, A letter to David Soudry, dated July 31, 2000.
 


\bibitem{Hn2}
 \bibline,
Caract{\'e}risation de la correspondance de Langlands locale par les 
facteurs $\epsilon$ de paires, {\it Invent.\ Math\/}.\ {\bf 113} (1993),
339--350.
 

\bibitem{Hn3}
\bibline,
Une preuve simple des conjectures de Langlands pour ${\rm GL}(n)$
sur un corps 
$p$-adique, {\it Invent.\ Math\/}.\ {\bf 139} (2000), 439--455.
 


\bibitem{Hw}
 \name{R.\ Howe},
 $\theta $-series and invariant theory, in  
{\it Automorphic Forms\/}, {\it Representations and $L$-functions\/} (Corvallis, Ore.,
1977),  Part 1,  275--285, {\it Proc.\ Sympos.\ Pure Math\/}.\ {\bf 33}, Amer.\ Math.\ Soc.,
Providence, RI, 1979.
 
 

\bibitem{J}
 \name{H.\ Jacquet},
Sur les repr{\'e}sentations des groupes r{\'e}ductifs $p$-adiques, 
{\it C.\ R.\ Acad.\ Sci.\ Paris S\'er. A-B\/} {\bf 280} (1975), A1271--A1272.
 

\bibitem{JP-SS}
 \name{H.\ Jacquet, I.\ I.\ Piatetski-Shapiro}, and \name{J.\ Shalika},
Rankin-Selberg convolutions, {\it Amer.\ J.\ Math\/}.\ {\bf 105} (1983), 367--464.
 
 

\bibitem{JS}
\name{H.\ Jacquet} and \name{J.\ A.\ Shalika},
On Euler products and the classification of automorphic forms.\ I; II,
{\it Amer.\ J.\ Math\/}.\ {\bf 103} (1981), 499--558 and 777--815.
 
 


\bibitem{JngS}
\name{D.\ Jiang} and \name{D.\ Soudry}, 
Generic representations and the local Langlands reciprocity law for 
$p$-adic ${\rm SO}(2n+1)$, preprint, 2001.
 


\bibitem{K1}
 \name{S.\ Kudla},
On the local theta-correspondence, {\it Invent.\ Math\/}.\ {\bf 83} (1986), 229--255.
 


\bibitem{K2}
\bibline,
The local Langlands correspondence: the non-Archimedean case, 
in {\it Motives\/}  (Seattle, WA, 1991), 365--391, {\it Proc.\ Sympos.\ Pure Math\/}.\
{\bf 55}, Amer.\ Math.\ Soc., Providence, RI, 1994.
 



\bibitem{Kn}
\name{A.\ W.\ Knapp},
Introduction to the Langlands program, in {\it Representation Theory and
Automorphic Forms\/} (Edinburgh, 1996), 245--302, {\it Proc.\ Sympos.\ Pure Math\/}.\ {\bf
61}, Amer.\ Math.\ Soc., Providence, RI, 1997.
 


\bibitem{L}
\name{R.\ P.\ Langlands}, 
Problems in the theory of automorphic forms, 
{\it Lecture Notes in Math\/}.\ {\bf 170} (1970), 18--61, Springer-Verlag, New York.
 


\bibitem{M}
 \name{C.\ M\oe glin},
Representations of ${\rm GL}(n,F)$ in the non-archimedean case, 
{\it Proc.\ Sympos.\ Pure Math\/}.\ {\bf 61} (1997), 303--319, 
Amer.\ Math.\ Soc.,
Providence, RI.
 

\bibitem{MVW}
 \name{C.\ M\oe glin, M.-F.\ Vign\'eras}, and \name{J.-L.\ Waldspurger},
{\it Correspondances de Howe sur un Corps $p$-adique\/},
{\it Lecture Notes in Math\/}.\ {\bf 1291}, Springer-Verlag, New York, 1987.
 


\bibitem{Mt}
\name{H.\ Matsumoto}, 
Sur les sous-groupes arithm{\'e}tiques des groupes semi-simples 
d{\'e}ploy{\'e}s, 
{\it Ann.\ Sci.\ {\rm \'{\it E}}cole Norm.\ Sup\/}.\ {\bf 2} (1969), 1--62.
 



\bibitem{P-S}
\name{I.\ Piatetski-Shapiro}, 
Multiplicity one theorems, in {\it Automorphic Forms\/}, {\it Representations
and $L$-functions\/} 
(Corvallis, Ore., 1977), Part 1,  209--212, {\it Proc.\ Sympos.\ Pure Math\/}.\ {\bf
33}, Amer.\ Math.\ Soc., Providence, RI, 1979.
 



\bibitem{R}
 \name{J.\ Rogawski},
{\it Automorphic Representations of Unitary Groups of Three Variables\/},
{\it Ann.\ of Math.\ Studies\/} {\bf 123},  Princeton Univ.\ Press, Princeton, NJ, 1990.
 


\bibitem{Sh1}
 \name{F.\ Shahidi}, 
On multiplicativity of local factors, 
{\it Israel Math.\ Conf.\ Proc\/}.\ {\bf 3},  279--289, Weizmann, Jerusalem, 1990.
 

\bibitem{Sh2}
\bibline,
A proof of Langlands' conjecture on Plancherel measures; complementary 
series for $p$-adic groups, {\it Ann.\ of Math\/}.\ {\bf 132} (1990), 273--330.
 


\bibitem{Shl}
\name{J.\ Shalika}, 
The multiplicity one theorem for ${\rm GL}_{n}$, 
{\it Ann.\ of Math\/}.\ {\bf 100}  (1974),  171--193. 
 

\bibitem{S1}
\name{D.\ Soudry},
{\it  Rankin-Selberg Convolutions for ${\rm SO}_{2l+1}\times {\rm GL}_n$\/}: {\it  Local
theory\/} {\it Mem.\ Amer.\ Math.\ Soc\/}.\ {\bf 105}, A.\ M.\ S., Providence,
RI (1993).
 

\bibitem{S2}
\bibline,
On the Archimedean theory of Rankin-Selberg convolutions for 
${\rm SO}_{2l+1}\times{\rm GL}_n$, {\it Ann.\ Sci.\ {\rm \'{\it E}}cole Norm.\ Sup\/}.\ {\bf 28}
(1995), 161--224.
 


\bibitem{S3}
\bibline,
Full multiplicativity of gamma factors for 
${\rm SO}_{2l+1}\times {\rm GL}_n$,   
{\it Proc.\ Conf. on $p$-adic Aspects of the Theory of
Automorphic  Representations\/} (Jerusalem, 1998), {\it  Israel J.\ Math\/}.\ {\bf 120}
(2000), part B,  511--561. 

\bibitem{S4}
\bibline,
A uniqueness theorem for representations of ${\rm GSO}(6)$ and 
the strong multiplicity one theorem for generic representations
of ${\rm GSp}(4)$, {\it Israel J.\ Math\/}.\ {\bf 58} (1987), 257--287.  

\bibitem{T}
\name{J.\ Tate},
Number theoretic background, 
in {\it Automorphic Forms\/},  {\it Representations and $L$-functions\/}
(Corvallis, Ore., 1977), Part 2,  3--26, {\it Proc.\ Sympos.\ Pure Math\/}.\ {\bf 33}, 
  Amer. Math. Soc., Providence, R.I., 1979.
   

\bibitem{Tr}
\name{F.\ Tr\`eves},
{\it Topological Vector Spaces\/}, {\it Distributions and Kernels\/},
Academic Press, New York,
1967.
 


\bibitem{V}
 \name{M.-F.\ Vign\'eras}, 
Correspondances entre representations automorphes de ${\rm GL}(2)$
sur une extension quadratique de ${\rm GSp}(4)$ sur ${\scriptstyle\Bbb Q}$, Conjecture
locale de Langlands pour ${\rm GSp}(4)$,  in {\it The Selberg Trace Formula and
Related Topics\/} (D.\ Hejhal, P.\ Sarnak, A.\ Terras, eds.), {\it Contemp.\ Math\/}.\ {\bf
53}, 463--527, Amer.\ Math.\ Soc., Providence, RI, 1986.
 

 
\bibitem{W}
\name{J.-L.\ Waldspurger},
D{\'e}monstration d'une conjecture de dualit{\'e} de Howe dans le cas 
$p$-adique, $p\neq 2$, {\it Israel Math.\ Conf.\ Proc\/}.\ {\bf 2} (1990),
267--324.
 


\bibitem{Wr}
 \name{G.\ Warner},
{\it Harmonic Analysis on Semi-Simple Lie Groups \/}, {\it Grund.\ Math.\ Wiss\/}.\ {\bf
188} (1972), Springer-Verlag, New York. 
 

\bibitem{Z}
 \name{A.\ V.\ Zelevinsky}, 
Induced representations of reductive $p$-adic groups.\ II.\
On irreducible representations of ${\rm GL}(n)$, {\it Ann.\ Sci.\ {\rm \'{\it E}}cole Norm.\ Sup\/}.\
{\bf 13} (1980), 165--210.
\vglue-3pt
\end{references}
